\documentclass[12pt, a4paper]{amsart}

\usepackage[latin1]{inputenc}

\usepackage{pdfsync}
\usepackage{amsmath,amsthm,amsfonts,amscd,amssymb,eucal,latexsym,dsfont,color}
\usepackage[all]{xy}
\usepackage{latexsym, amssymb, amscd, mathrsfs, graphics, fullpage}
\usepackage{hyperref}
\usepackage{stmaryrd}

\newcommand{\ga}{\mathfrak{a}}

\newcommand{\id}{\mathrm{id}}
\newcommand{\bfG}{\mathds{G}}

\newcommand{\R}{\mathds{R}}

\newcommand{\vl}{\,\varolessthan\,}
\newcommand{\vg}{\,\varogreaterthan\,}
\swapnumbers

\theoremstyle{definition}

\def\gG{\mathfrak{G}}
\def\gH{\mathfrak{H}}

\def\gM{\mathfrak{M}}
\def\gN{\mathfrak{N}}

\def\gt{\mathfrak{t}}
\def\ga{\mathfrak{a}}
\def\gb{\mathfrak {b}}

\newcommand{\op}{\mathrm{o}}

\begin{document}

%%%%%%%%%%%%%%%%%%%%%TITLE
\title[Locally Compact Quantum groupoids] 
{Locally Compact Quantum groupoids}
\author{Michel Enock}

\address{Institut de Math\'ematiques de Jussieu\\Unit\'{e} Mixte de Recherche 7586, Sorbonne Universit\'{e} / Universit\'{e} de Paris /
CNRS  \\Case 247, 4 place Jussieu, 75252 Paris Cedex, France}
 \email{michel.enock@math.cnrs.fr}

\begin{abstract}
The theory of measured quantum groupoids, as defined in \cite{L} and \cite{E4}, was made to generalize the theory of quantum groups (\cite{KV1}, \cite{KV2}), but was only defined in a von Neumann algebra setting; Th. Timmermann constructed locally compact quantum groupoids, which is a $\bf C^*$-version of quantum groupoids \cite{Ti3}. Here, we associate to such a locally compact quantum groupoid a measured quantum groupoid in which it is weakly dense; we then associate to a measured quantum groupoid a locally compact quantum groupoid which is weakly dense in the measured quantum groupoid, but such a locally compact quantum groupoid may be not unique; we construct a duality of locally compact quantum groupoids. Following (\cite{Ti4}), we define the actions of a locally compact quantum groupoid. We give then examples of locally compact quantum groupoids.  \end{abstract}

\maketitle

\tableofcontents
\newpage
 %%%%%%%%%%intro
\section{Introduction}
\label{intro}
\subsection{Locally compact quantum groups}
The theory of locally compact quantum groups, developed by J. Kustermans and S. Vaes (\cite{KV1}, \cite{KV2}), provides a comprehensive framework for the study of quantum groups in the setting of $\bf C^*$-algebras and von Neumann algebras. It includes a far reaching generalization of the classical Pontrjagin duality of locally compact abelian groups, that covers all locally compact groups. Namely, if $G$ is a locally compact group, its von Neumann algebra is $L^\infty(G, \mu)$ (where $\mu$ is the left Haar measure on $G$), and its dual von Neumann algebra is $\mathcal L(G)$ generated by the left regular representation $\lambda_G$ of $G$ on $L^2(\mu)$, equipped with a coproduct $\Gamma_G$ from $\mathcal L(G)$ on $\mathcal L(G)\otimes\mathcal L(G)$ defined, for all $s\in G$, by $\Gamma_G(\lambda_G(s))=\lambda_G(s)\otimes\lambda_G(s)$, and with a normal semi-finite faithful weight, called the Plancherel weight $\varphi_G$, associated via the Tomita-Takesaki construction, to the left Hilbert algebra defined by the algebra $\mathcal K(G)$ of continuous functions with compact support (with convolution as product), this weight $\varphi_G$ being left- and right-invariant with respect to $\Gamma_G$ (\cite{T}, VII, 3). 

This theory builds on many preceding works, by G. Kac, G. Kac and L. Vainerman, J.-M. Schwartz and the author (\cite{ES1}, \cite{ES2}), S. Baaj and G. Skandalis (\cite{BS}), A. Van Daele (\cite{VD1}), S. Woronowicz {\cite{W1}, \cite{W5}, \cite{W6}) and many others. See the monography written by T. Timmermann for a survey of that theory (\cite{Ti1}), and the introduction of \cite{ES2} for a sketch of the historical background. It seems to have reached now a stable situation, because it fits the needs of operator algebraists for many reasons: 

First,  the axioms of this theory are very simple and elegant: they
    can be given in both $\bf C^*$-algebras and von Neumann algebras, and
    these two points of view are equivalent, as A. Weil had shown it
    was the fact for groups (namely any measurable group equipped with
    a left-invariant positive non zero measure bears a topology which makes it
    locally compact, and this measure is then the Haar measure
    (\cite{W}, Appendice I)). In a von Neumann setting, a locally compact quantum
    group is just a von Neumann algebra, equipped with a
    co-associative coproduct, and two normal faithful semi-finite
    weights, one left-invariant with respect to that coproduct, and
    one right-invariant. Then, many other data are constructed, in
    particular a multiplicative unitary (as defined in \cite{BS})
    which is manageable (as defined in \cite{W6}).

Second, all
    preceeding attemps (\cite{ES2}, \cite{W5}) appear as particular
    cases of locally compact quantum groups; and many interesting
    examples were constructed (\cite{W2}, \cite{W3}, \cite{VV}).

Third,  many constructions of harmonic analysis, or concerning
    group actions on $\bf C^*$-algebras and von Neumann algebras, were
    generalized to locally compact quantum groups (\cite{V2}).

    Finally, many constructions made by algebraists at the level of
    Hopf $*$-algebras, or multipliers Hopf $*$-algebras, can be
    generalized for locally compact quantum groups. This is the case,
    for instance, for Drinfel'd double of a quantum group (\cite{D}),
    and for Yetter-Drinfel'd algebras which were well-known in an
    algebraic approach in \cite{M}.

\subsection{Measured Quantum Groupoids}
In two articles (\cite{Val1},\cite{Val2}), J.-M. Vallin has introduced two notions (pseudo-multiplicative unitary, Hopf bimodule), in order to generalize, to the groupoid case, the classical notions of multiplicative unitary (\cite{BS}) and of a co-associative coproduct on a von Neumann algebra. Then, F. Lesieur (\cite{L}), starting from a Hopf bimodule, when there exist a left-invariant operator-valued weight and a right-invariant operator-valued weight, mimicking in that wider setting what was done in (\cite{KV1}, \cite{KV2}), obtained a pseudo-multiplicative unitary, and called ``measured quantum groupoids'' these objects. A new set of axioms had been given in an appendix of \cite{E3}. In \cite{E3} and \cite{E4}, most of the results given in \cite{V2} were generalized to measured quantum groupoids. 
Some trivial examples were given in \cite{L}. A more interesting example was constructed in \cite{ET} : here, a quantum transformation groupoid is defined, from a right action of a locally compact quantum group $\bf G$ on a von Neumann algebra $N$, which generalizes the transformation groupoid given by a locally compact group $G$ having a right action $\ga$ on a locally compact space $X$. 

This theory, up to now, had one important defect : it was only a theory in a von Neumann algebra setting. 

\subsection{Locally compact quantum groupoids}
Th. Timmermann had made many attemps in order to provide a $\bf C^*$-algebra version of it (see \cite{Ti1} for a survey); these attemps were fruitful, but not sufficient to complete a theory equivalent to the von Neumann one. 

This is the subject of this article, in order to get, for quantum groupoids, as it is for quantum groups, axioms in both $\bf C^*$-algebras and von Neumann algebras. In a $\bf C^*$-algebra setting, we first recall Timmermann's theory of \emph{locally compact quantum groupoids}; we then prove that to such an object $\bfG$, it is possible to associate a measured quantum groupoid $\gG$ such that the $\bf C^*$-algebra of $\bfG$ is a dense sub-$\bf{C^*}$ algebra of the underlying von Neumann algebra of the measured quantum groupoid. 

We must quote older articles about quantum groupoids, first a purely algebraic construction (\cite{Sc}), and second a $\bf C^*$ construction in a particular situation (\cite{KVD1} and \cite{KVD2}). 

The article is organized as follows:

 In chapter \ref{MQG} are recalled the definition of measured quantum groupoids; first the definitions of the relative tensor of Hilbert space (\ref{spatial}), the fiber product of von Neumann algebra (\ref{fiber}) and then the definitions of Hopf-bimodules (\ref{defHopf}), and of measured quantum groupoids (\ref{defMQG}).

In chapter \ref{LCQG}, we give a definition of locally compact quantum groupoids. First (\ref{C*fiber}), we recall  definitions of $\bf C^*$-relative products of Hilbert spaces and $\bf C^*$-fiber product of $\bf C^*$-algebras ([Ti2]); after, we recall (\ref{weights}) basic results about weights on $\bf C^*$-algebras, mostly due to F. Combes (\cite{C1}, \cite{C2}). We then recall Kusterman's definition of $\bf C^*$-valued weights (\cite{K2}), and define a class of $\bf C^*$-valued weights which are restrictions of operator-valued weights. Using all these notions, we then give a new definition of \emph{fiber products of $\bf C^*$-algebras} (\ref{fiber}).

 In chapter \ref{fLCQG}, we prove that, to any locally compact quantum groupoid $\bfG$, we can associate a canonical measured quantum groupoid $\gG$ such that $\bfG$ is a dense sub-$\bf C^*$-algebra of $\gG$. Then, we say that $\bfG$ is a \emph{locally compact sub-quantum groupoid} of $\gG$ (\ref{defLCQG}). 

In chapter \ref{fMQG}, starting from a measured quantum groupoid, we construct a canonical locally compact sub-quantum groupoid; such a sub $\bf C^*$-algebra, which is a locally compact quantum groupoid, may be not unique (\ref{ex}). 

In chapter \ref{Duality}, from any locally compact quantum groupoid, we construct a dual one, and prove that the bidual is isomorphic to the initial locally compact quantum groupoid (equal if we identify the canonical Hilbert spaces, on which the $\bf C^*$-algebras are constructed).  

In chapter \ref{action}, we construct the action of a locally compact quantum groupoid on a $\bf C^*$-algebra. 

In chapter \ref{Examples}, we recall several examples of locally compact quantum groupoid. 

We are indebted to Thomas Timmermann who had found several mistakes in a preliminary version of this article. 
%%%%%MQG

\section{Measured quantum groupoids }
\label{MQG}
In this chapter, we recall the definition of the relative tensor
product of Hilbert spaces, and of the fiber product of von Neumann
algebra (\ref{spatial}). Then, we recall the definition of a Hopf
bimodule (\ref{defHopf}) and a co-inverse. We then give the definition of a measured quantum groupoid, recall the construction of the pseudo-multiplicative unitary  (\ref{defW}), and all the data constructed then, including duality of measured quantum groupoids (\ref{data}).

\subsection{Relative tensor products of Hilbert spaces (\cite{C}, \cite{S}, \cite{T}, \cite{EVal})}
\label{spatial}
 Let $N$ be a von Neumann algebra, $\nu$ a normal semi-finite
 faithful weight on $N$; we shall denote by $H_\nu$, $\gN_\nu$,
 \ldots the canonical objects of the Tomita-Takesaki theory associated
 to the weight $\nu$. 

Let $\alpha$ be a non-degenerate faithful representation of $N$ on a Hilbert space $\mathcal H$. The set of $\nu$-bounded elements of the left module $_\alpha\mathcal H$ is
\[D(_\alpha\mathcal{H}, \nu)= \lbrace \xi \in \mathcal{H} : \exists C < \infty ,\| \alpha (y) \xi\|
\leq C \| \Lambda_{\nu}(y)\|,\forall y\in \gN_{\nu}\rbrace.\]
 For any $\xi$ in $D(_\alpha\mathcal{H}, \nu)$, there exists a bounded operator
$R^{\alpha,\nu}(\xi)$ from $H_\nu$ to $\mathcal{H}$ such that
\[R^{\alpha,\nu}(\xi)\Lambda_\nu(y) = \alpha (y)\xi \quad\text{for
  all } y \in \gN_{\nu},\]
and this operator exchanges the representations of $N$. 
If $\xi$ and $\eta$ are bounded vectors, we define the operator product 
\[\langle \xi|\eta\rangle _{\alpha,\nu} = R^{\alpha,\nu}(\eta)^* R^{\alpha,\nu}(\xi),\]
which belongs to $\pi_{\nu}(N)'$. Using Tomita-Takesaki theory, this last algebra will be
identified with the opposite von Neumann algebra $N^{\op}$. We shall use also the operator 
\[\theta^{\alpha, \nu}(\xi, \eta)=R^{\alpha,\nu}(\xi)R^{\alpha,\nu}(\eta)^*\]
which belongs to $\alpha (N)'$. If now $\beta$ is a non-degenerate faithful anti-representation of $N$ on a Hilbert space $\mathcal K$, the relative tensor product $\mathcal K\underset{\nu}{_\beta\otimes_\alpha}\mathcal H$ is the completion of the algebraic tensor product $K\odot D(_\alpha\mathcal{H}, \nu)$ by the scalar product defined by 
\[(\xi_1\odot\eta_1 |\xi_2\odot\eta_2 )= (\beta(\langle \eta_1|
\eta_2\rangle _{\alpha,\nu})\xi_1 |\xi_2)\]
for all $\xi_1, \xi_2\in \mathcal{K}$ and $\eta_1, \eta_2 \in D(_\alpha\mathcal{H},\nu)$.
If $\xi\in \mathcal{K}$ and $\eta\in D(_\alpha\mathcal{H},\nu)$, we
 denote by $\xi\underset{\nu}{_\beta\otimes_\alpha}\eta$ the image of $\xi\odot\eta$ into $\mathcal K\underset{\nu}{_\beta\otimes_\alpha}\mathcal H$. Writing $\rho^{\beta, \alpha}_\eta(\xi)=\xi\underset{\nu}{_\beta\otimes_\alpha}\eta$, we get a bounded linear operator from $\mathcal H$ into $\mathcal K\underset{\nu}{_\beta\otimes_\alpha}\mathcal H$, which is equal to $1_\mathcal K\otimes_\nu R^{\alpha, \psi}(\eta)$. 
 
 If $x\in\mathcal D(\sigma_{i/2}^\nu)$, then $\alpha(x)D(_\alpha\mathcal{H},\nu)\subset D(_\alpha\mathcal{H},\nu)$, and we have (\cite{S}, 2.2 b) :
 \[\xi\underset{\nu}{_\beta\otimes_\alpha}\alpha(x)\eta=\beta(\sigma_{i/2}(x))\xi\underset{\nu}{_\beta\otimes_\alpha}\eta\] 

 Changing the weight $\nu$ will give an isomorphic Hilbert space, but
 the isomorphism will not exchange elementary tensors!

We shall denote by $\sigma_\nu$ the relative flip, which is a unitary sending $\mathcal{K}\underset{\nu}{_\beta\otimes_\alpha}\mathcal{H}$ onto $\mathcal{H}\underset{\nu^{\op}}{_\alpha\otimes _\beta}\mathcal{K}$, defined by
\[\sigma_\nu
(\xi\underset{\nu}{_\beta\otimes_\alpha}\eta)=\eta\underset{\nu^{\op}}{_\alpha\otimes_\beta}\xi\]
 for all $\xi \in D(\mathcal {K}_\beta ,\nu^{\op} )$ and $\eta \in
 D(_\alpha \mathcal {H},\nu)$.

In (\cite{DC1}, chap. 11), De Commer had shown that, if $N$ is finite-dimensional, the Hilbert space $\mathcal K\underset{\nu}{_\beta\otimes_\alpha}\mathcal H$ can be isometrically imbedded into the usual Hilbert tensor product $\mathcal K\otimes\mathcal H$. 

%%%%basis
\subsubsection{\bf Definition}
\label{basis}
There exists (\cite{C}, prop.3) a family  $(e_i)_{i\in I}$ of 
$\nu$-bounded elements of $_\alpha\mathcal H$, such that
\[\sum_i\theta^{\alpha, \nu} (e_i ,e_i )=1\]
Such a family will be called an $(\alpha,\nu)$-{\it basis} of $\mathcal H$. Then, for any $i\in I$, the image of $\theta^{\alpha, \nu} (e_i ,e_i )$ is included in the closure of the subspace $\{\alpha(n)e_i, n\in N\}$. 
 \newline
 In that situation, let us consider, for all $n\in \bf N$ and finite $J\subset I$ with $|J|=n$, the $(1,n)$ matrix $R_J=(R^{\alpha, \nu}(e_i))_{i\in J}$. As $R_JR_J^*\leq1$, we get that $\|R_J\|\leq 1$, and that the $(n,n)$ matrix $(<e_i, e_j>_{\alpha, \nu})_{i,j\in J}\in M_n(N^o)$ is less than the unit matrix. 
 \newline
 It is possible (\cite{EN} 2.2) to construct 
an $(\alpha,\nu )$-basis of $\mathcal H$, $(e_i)_{i\in I}$, such that the
operators $R^{\alpha, \nu}(e_i)$ are partial isometries with final supports 
$\theta^{\alpha, \nu}(e_i ,e_i )$ 2 by 2 orthogonal, and such that, if $i\neq j$, then 
$<e_i ,e_j>_{\alpha, \nu}=0$. Such a family will be called an $(\alpha, \nu)$-{\it orthogonal basis} of $\mathcal H$. 
\newline
 We have, then :
 \[R^{\alpha, \nu}(\xi)=\sum_i\theta^{\alpha, \nu}(e_i, e_i)R^{\alpha, \nu}(\xi)=\sum_iR^{\alpha, \nu}(e_i)<\xi, e_i>_{\alpha, \nu}\]
 \[<\xi, \eta>_{\alpha, \psi}=\sum_i<\eta, e_i>_{\alpha, \psi}^*<\xi, e_i>_{\alpha, \nu}\]
 \[\xi=\sum_iR^{\alpha, \nu}(e_i)J_\psi\Lambda_\psi(<\xi, e_i>^o_{\alpha, \nu})\]
 the sums being weakly convergent. 
 \newline
 Moreover, we get that, for all $n$ in $N$, $\theta^{\alpha, \nu}(e_i, e_i)\alpha(n)e_i=\alpha(n)e_i$, and $\theta^{\alpha, \nu}(e_i, e_i)$ is the orthogonal projection on the closure of the subspace $\{\alpha(n)e_i, n\in N\}$. 

%%%%basic construction
\subsubsection{\bf Basic construction} (\cite{EN}, 3.1)
\label{basic}
Let $M_0\subset M_1$ be an inclusion of von Neumann algebras, $\psi_1$ a normal semi-finite faithful weight on $M_1$; then $M_2=J_{\psi_1} M'_0J_{\psi_1}$ is a von Neumann algebra called the basic construction from the inclusion $M_0\subset M_1$; let $\psi_0$  be a normal semi-finite faithful weight on $M_0$; if $\xi$, $\eta$ belong to $D(H_{\psi_1}, \psi_0^o)$, then $\theta^{\psi_0^o}(\xi, \eta)$ belongs to $M_2$, and the linear span of these operators is a dense ideal in $M_2$.

%%%ovw
\subsection{\bf{Operator-valued weights}}
\label{ovw}
Let $M_0\subset M_1$ be an inclusion of von Neumann algebras  (for simplification, these algebras will be supposed to be $\sigma$-finite), equipped with a normal faithful semi-finite operator-valued weight $T_1$ from $M_1$ to $M_0$ (to be more precise, from $M_1^{+}$ to the extended positive elements of $M_0$ (cf. \cite{T} IX.4.12)). Let $\psi_0$ be a normal faithful semi-finite weight on $M_0$, and $\psi_1=\psi_0\circ T_1$; for $i=0,1$, let $H_i=H_{\psi_i}$, $J_i=J_{\psi_i}$, $\Delta_i=\Delta_{\psi_i}$ be the usual objects constructed by the Tomita-Takesaki theory associated to these weights. 
 \newline
Following (\cite{EN} 10.6), for $x$ in $\gN_{T_1}$, we shall define $\Lambda_{T_1}(x)$ by the following formula, for all $z$ in $\gN_{\psi_{0}}$ :
\[\Lambda_{T_1}(x)\Lambda_{\psi_{0}}(z)=\Lambda_{\psi_1}(xz)\]
This operator belongs to $Hom_{M_{0}^o}(H_{0}, H_1)$; if $x$, $y$ belong to $\gN_{T_1}$, then $\Lambda_{T_1}(x)\Lambda_{T_1}(y)^*$ belongs to the von Neumann algebra $M_{2}=J_1M'_0J_1$, which is called the basic construction made from the inclusion $M_0\subset M_1$, and
$\Lambda_{T_1}(x)^*\Lambda_{T_1}(y)=T_1(x^*y)\in M_0$. 
\newline
By Tomita-Takesaki theory, the Hilbert space $H_1$ bears a natural structure of $M_1-M_1^o$-bimodule, and, therefore, by restriction, of $M_0-M_0^o$-bimodule. Let us write $r$ for the canonical representation of $M_0$ on $H_1$, and $s$ for the canonical antirepresentation given, for all $x$ in $M_0$, by $s(x)=J_1r(x)^*J_1$. Let us have now a closer look to the subspaces $D(H_{1s}, \psi_0^o)$ and $D(_rH_1, \psi_0)$. If $x$ belongs to $\gN_{T_1}\cap\gN_{\psi_1}$, we easily get that $J_1\Lambda_{\psi_1}(x)$ belongs to $D(_rH_1, \psi_0)$, with :
\[R^{r, \psi_0}(J_1\Lambda_{\psi_1}(x))=J_1\Lambda_{T_1}(x)J_0\]
and $\Lambda_{\psi_1}(x)$ belongs to $D(H_{1s}, \psi_0^o)$, with :
\[R^{s, \psi_0^o}(\Lambda_{\psi_1}(x))=\Lambda_{T_1}(x)\]
The subspace $D(H_{1s}, \psi_0^o)\cap D(_rH_1, \psi_0)$ is dense in $H_1$; more precisely, let $\mathcal T_{\psi_1, T_1}$ be the algebra made of elements $x$ in $\gN_{\psi_1}\cap\gN_{T_1}\cap\gN_{\psi_1}^*\cap\gN_{T_1}^*$, analytical with respect to $\psi_1$, and such that, for all $z$ in $\mathbb{C}$, $\sigma^{\psi_1}_z(x_n)$ belongs to $\gN_{\psi_1}\cap\gN_{T_1}\cap\gN_{\psi_1}^*\cap\gN_{T_1}^*$. Then (\cite{E4}, 2.2.1):
\newline
(i) the algebra $\mathcal T_{\psi_1, T_1}$ is weakly dense in $M_1$; it will be called Tomita's algebra with respect to $\psi_1$ and $T_1$; 
\newline
(ii) for any $x$ in  $\mathcal T_{\psi_1, T_1}$, $\Lambda_{\psi_1}(x)$ belongs to $D(H_{1s}, \psi_0)\cap D(_rH_1, \psi_0)$;
\newline
(iii) for any $\xi$ in $D(H_{1s}, \psi_0^o))$, there exists a sequence $x_n$ in $\mathcal T_{\psi_1, T_1}$ such that $\Lambda_{T_1}(x_n)=R^{s, \psi_0^o}(\Lambda_{\psi_1}(x_n))$ is weakly converging to $R^{s, \psi_0^o}(\xi)$ and $\Lambda_{\psi_1}(x_n)$ is converging to $\xi$. 
\newline
More precisely, in (\cite{E2}, 2.3) was constructed an increasing sequence of projections $p_n$ in $M_1$, converging to $1$, and elements $x_n$ in $\mathcal T_{\psi_1, T_1}$ such that $\Lambda_{\psi_1}(x_n)=p_n\xi$. We then get that :
\begin{align*}
T_1(x_n^*x_n)
&=
<R^{s, \psi_0^o}(\Lambda_{\psi_1}(x_n)), R^{s, \psi_0^o}(\Lambda_{\psi_1}(x_n))>_{s, \psi_0^o}\\
&=
<p_n\xi, p_n\xi>_{s, \psi_0^o}\\
&=
R^{s, \psi_0^o}(\xi)^*p_nR^{s, \psi_0^o}(\xi)
\end{align*}
which is increasing and weakly converging to $<\xi, \xi>_{s, \psi_0^o}$. 

%%%basisT
\subsubsection{\bf{Theorem}}(\cite{EN}, 10.3, 10.7, 10.11, \cite{E3}, 2.10)
\label{basisT}
{\it Let $M_1$ be a von Neumann algebra, $M_0$ be a von Neumann subalgebra of $M_1$, $\nu$ a faithful semi-finite normal weight on $M_0$ and $T$ a normal faithful semi-finite operator-valued weight from $M_1$ onto $M_0$; let $r$ be the inclusion of $M_0$ into $B(H_{\nu\circ T})$, and $s$ the anti-$*$-homomorphism from $M_0$ into $B(H_{\nu\circ T})$ defined by ($x\in M_0$) $x\mapsto J_{\nu\circ T}x^*J_{\nu\circ T}$; let us define, for $x\in\gN_T$, 
$\Lambda_{T}(x)\in B(H_\nu, H_{\nu\circ\circ T})$ by ($z\in\gN_\nu$) :
\[\Lambda_{T}(x)\Lambda_{\nu}(z)=\Lambda_{\nu\circ T}(xz)\]
(i) let $M_2$ be the basic construction made from the inclusion $M_0\subset M_1$; then, for any $x$, $y$ in $\gN_T$, then $\Lambda_T(x)\Lambda_T(y)^*$ belongs to $M_2$, and the von Neumann algebra $M_2$ is generated by these operators; moreover, there exists a normal faithful semi-finite operator-valued weight $T_2$ from $M_2$ to $M_1$ such that, $T_2(\Lambda_T(x)\Lambda_T(y)^*)=xy^*$. 
Let $X$ belong to $Hom_{M_0^o}(H_\nu, H_{\nu\circ T})$ such that $XX^*$ belongs to $\gM_{T_2}^+$; then, there exists a unique element $\Phi_1(X)$ in $M_1$ such that $T_2(X\Lambda_{T}(a)^*)=\Phi_1(X)a^*$, for all $a\in\gN_{T}$; this application $\Phi_1$ is an injective application of $(M_1, M_0)$-bimodule, and we have $\Phi_1(\Lambda_T(x))=x$, for all $x\in\gN_T$. 
\newline
(ii) there exists a family $(e_i)_{i\in I}$ in $\gN_T\cap\gN_T^*\cap\gN_{\nu\circ T}\cap\gN_{\nu\circ T}^*$ such that the operators $\Lambda_{T}(e_i)$ are partial isometries, with $T(e_j^*e_i)=0$ if $j\ne i$, and with their final supports $\Lambda_T(e_i)\Lambda_T(e_i)^*$ two by two orthogonal projections of sum $1$; moreover, for all $i\in I$, we have $e_i=e_iT(e_i^*e_i)$, and, for all $x\in\gN_T$, we have :
\[\Lambda_T(x)=\sum_i\Lambda_T(e_i)T(e_i^*x)\]
\[x=\sum_ie_iT(e_i^*x)\]
these two sums being weakly convergent. Such a family $(e_i)_{i\in I}$ will be called a basis for $(T, \nu)$. Moreover, the family $J_{\nu\circ T}\Lambda_{\nu\circ T}(e_i)$ is a basis for $(\alpha, \nu)$ (\ref{basis}), where $\alpha$ is the inclusion $M_0\subset M_1$. 

Moreover, the vectors $(\Lambda_{\nu\circ T}(e_i))_{i\in I}$ are a $(s, \nu^o)$-orthogonal basis of $H_{\nu\circ T}$, and the vectors $(J_{\nu\circ T}\Lambda_{\nu\circ T}(e_i))_{i\in I}$ are a $(r, \nu)$-orthogonal basis of $H_{\nu\circ T}$. 

(iii) for any $\xi\in D((H_{\nu\circ T})_s, \nu^o)$, there exists a sequence $x_n$ in $\gN_T\cap\gN_{\nu\circ T}$ such that $\Lambda_{\nu\circ T}(x_n)$ is converging to $\xi$, and $\Lambda_T(x_n)=R^{s, \nu^o}(\Lambda_{\nu\circ T}(x_n))$ is weakly converging to $R^{s, \nu^o}(\xi)$; equivalently, for any $\eta\in D(_rH_{\nu\circ T}, \nu)$, there exists a sequence $y_n$ in $\gN_T\cap\gN_{\nu\circ T}$ such that $\Lambda_{\nu\circ T}(y_n)$ is converging to $\eta$ and $J_{\nu\circ T}\Lambda_T(y_n)J_\nu=R^{r, \nu}(\Lambda_{\nu\circ T}(y_n))$ is weakly converging to $R^{r, \nu}(\eta)$.}

\begin{proof}
Result (i) is just (\cite{EN}, 10.3, 10.7 and 10.11). Let us simplify and clarify the proof given in (\cite{E3}, 2.10) for (ii) : let us first remark that $\Lambda_T(x)^*\Lambda_T(x)=T(x^*x)$.

By density of $\Lambda_{\nu\circ T}(\gN_T\cap\gN_T^*\cap\gN_{\nu\circ T}\cap\gN_{\nu\circ T}^*)$ into $H_{\nu\circ T}$, we can choose a family $(x_\alpha)_{\alpha\in A}$ in $\gN_T\cap\gN_T^*\cap\gN_{\nu\circ T}\cap\gN_{\nu\circ T}^*$ such that the closed subspaces $\overline{J_{\nu\circ T}M_0J_{\nu\circ T}\Lambda_{\nu\circ T}(x_\alpha)}$ are two by two orthogonal and that $H_{\nu\circ T}$ is the sum of these subspaces. Then, for any $\alpha$, we can remark that, for any $f\in L^\infty (\bf R^+)$, we have :
\[f(\Lambda_T(x_\alpha)\Lambda_T(x_\alpha)^*)\Lambda_T(x_\alpha)=\Lambda_T(x_\alpha)f(T(x_\alpha^*x_\alpha))=\Lambda_T(x_\alpha f(T(x_\alpha^*x_\alpha)))\]
 This formula is clear for any polynomial function, then by norm continuity, for any continuous function on $[0, \|T(x_\alpha^*x_\alpha)\|]$, and by weak continuity, we get the result. So, taking the function $f_{\alpha, n}(t)=\chi_{]\|T(x_\alpha^*x_\alpha\|/(n+1), \|T(x_\alpha^*x_\alpha\|/n]}(t)t^{-1/2}$, and defining $I$ as the subset of $A \times \bf N$ such that $x_\alpha f_{\alpha, n}(T(x_\alpha^*x_\alpha))\ne 0$ and $e_i=x_\alpha f_{\alpha, n}(T(x_\alpha^*x_\alpha))$, we obtain the appropriate family $(e_i)_{i\in I}$. 
 
 Result (iii) is taken from (\cite{E2}, 2.3(i)), or can be deduced from (i). \end{proof}

%%%%lemT
\subsubsection{\bf Lemma}
\label{lemT}
{\it Let $N$ be a von Neumann algebra, $\nu$ be a faithful semi-finite normal weight on $N$. Let $\alpha$ be a faithful non degenerate representation of $N$ into a von Neumann algebra $M_1$, and $T_1$ be a normal faithful semi-finite operator-valued weight from $M_1$ onto $\alpha (N)$. Let $\beta$ be a faithful non degenerate anti-representation of $N$ into a von Neumann algebra $M_2$, and $T_2$ be a normal faithful semi-finite operator-valued weight from $M_2$ onto $\beta(N)$. For $x\in\gN_{T_1}$, let us define $\Lambda_{T_1}(x)\in B(H_\nu, H_{\nu\circ\alpha^{-1}\circ T_1})$ by ($z\in\gN_\nu$) :
\[\Lambda_{T_1}(x)\Lambda_{\nu}(z)=\Lambda_{\nu\circ\alpha^{-1}\circ T_1}(x\alpha(z))\]
and, for $y\in\gN_{T_2}$, let us define as well $\Lambda_{T_2}(y)\in B(H_\nu, H_{\nu^o\circ\beta^{-1}\circ T_2^o})$ by :
\[\Lambda_{T_2}(y)J_\nu\Lambda_\nu(z)=\Lambda_{\nu^o\circ\beta^{-1}\circ T_2}(y\beta(z^*))\]
Then :

(i) for any $x\in\gN_{T_1}\cap\gN_{\nu\circ\alpha^{-1}\circ T_1}$, $J_{\nu\circ\alpha^{-1}\circ T_1}\Lambda_{\nu\circ\alpha^{-1}\circ T_1}(x)$ belongs to $D(_\alpha H_{\nu\circ\alpha^{-1}\circ T_1}, \nu)$ and $R^{\alpha, \nu}(J_{\nu\circ\alpha^{-1}\circ T_1}\Lambda_{\nu\circ\alpha^{-1}\circ T_1}(x))=J_{\nu\circ\alpha^{-1}\circ T_1}\Lambda_{T_1}(x)J_\nu$; moreover, if $x_1$ belongs to $\gN_{T_1}\cap\gN_{\nu\circ\alpha^{-1}\circ T_1}$, we have $\Lambda_{T_1}(x)^*\Lambda_{\nu\circ\alpha^{-1}\circ T_1}(x_2)=\Lambda_\nu(\alpha^{-1}T_1(x^*x_2))$; if $x'_1$ $x'_2$ belong to $\gN_{T_1}$, we have $\Lambda_{T_1}^*(x'_1)^*\Lambda_{T_1}(x'_2)=\alpha^{-1}T_1(x_1^{'*}x'_2)$. 

(ii) for any $y\in\gN_{T_2}\cap \gN_{\nu^o\circ\beta^{-1}\circ T_2}$, $\Lambda_{\nu^o\circ\beta^{-1}\circ T_2}(y)$ belongs to $D(_\beta H_{\nu^o\circ\beta^{-1}\circ T_2} \nu^o)$ and $R^{\beta, \nu^o}(\Lambda_{\nu^o\circ\beta^{-1}\circ T_2}(y))=\Lambda_{T_2}(y)$. Moreover, if $y_1$ belongs to $\gN_{T_2}\cap\gN_{\nu^o\circ\beta^{-1}\circ T_2}$, we have $\Lambda_{T_2}(y)^*\Lambda_{\nu^o\circ\beta^{-1}\circ T_2}(y_1)= J_\nu\Lambda_\nu(\beta^{-1}T_2(y_1^*y))$ and, if $y'_1$, $y'_2$ are in $\gN_{T_2}$, then $\Lambda_{T_2}(y'_1)^*\Lambda_{T_2}(y'_2)=J_\nu\beta^{-1}(T_2(y_2^{'*}y'_1))J_\nu$ (which belongs to $J_\nu N^oJ_\nu=N$). 

(iii) for any $x_1$, $x_2$ in $\gN_{T_1}\cap\gN_{\nu\circ\alpha^{-1}\circ T_1}$, $y_1$, $y_2$ in $\gN_{T_2^o}\cap \gN_{\nu^o\circ\beta^{-1}\circ T_2^o}$, the scalar product :
\[(\Lambda_{\nu^o\circ\beta^{-1}\circ T_2}(y_1)\underset{\nu}{_\beta\otimes_\alpha} J_{\nu\circ\alpha^{-1}\circ T_1}\Lambda_{\nu\circ\alpha^{-1}\circ T_1}(x_1)|\Lambda_{\nu^o\circ\beta^{-1}\circ T_2}(y_2)\underset{\nu}{_\beta\otimes_\alpha} J_{\nu\circ\alpha^{-1}\circ T_1}\Lambda_{\nu\circ\alpha^{-1}\circ T_1}(x_2))\]
is equal to $(J_\nu\Lambda_\nu(\beta^{-1}T_2(y_1^*y_2))|\Lambda_\nu(\alpha^{-1}T_1(x_1^*x_2)))$. }
\begin{proof}
It is straightforward to get that :
\[J_{\nu\circ\alpha^{-1}\circ T_1}\Lambda_{T_1}(x)J_\nu\Lambda_\nu (z)=\alpha(z)J_{\nu\circ\alpha^{-1}\circ T_1} \Lambda_{\nu\circ\alpha^{-1}\circ T_1}(x)\]
The formula about $\Lambda_{T_1}(x)^*$ is an easy calculation, and then we get the formula about $x'_1$, $x'_2$ if $x'_1\in \gN_{T_1}\cap\gN_{\nu\circ\alpha^{-1}\circ T_1}$, and then for all $x'_1$, $x'_2$ by continuity, , which finishes the proof of (i). 

The formula about $\Lambda_{\nu^o\circ\beta^{-1}\circ T_2}(y)$ is just the definition of $\Lambda_{T_2}$; the other results of (ii) are proved the same way as (i). 

Using (i), we get that the scalar product is equal to :
\[(\alpha(<\Lambda_{\nu^o\circ\beta^{-1}\circ T_2}(y_1), \Lambda_{\nu^o\circ\beta^{-1}\circ T_2}(y_2)>_{\beta, \nu^o})J_{\nu\circ\alpha^{-1}\circ T_1}\Lambda_{\nu\circ\alpha^{-1}\circ T_1}(x_1)|J_{\nu^o\circ\alpha^{-1}\circ T_1}\Lambda_{\nu^o\circ\alpha^{-1}\circ T_1}(x_2))\]
which, using (ii), is equal to 
\begin{multline*}
(\alpha\circ\beta^{-1}(T_2(y_2^*y_1))J_{\nu\circ\alpha^{-1}\circ T_1}\Lambda_{\nu\circ\alpha^{-1}\circ T_1}(x_1)|J_{\nu^o\circ\alpha^{-1}\circ T_1}\Lambda_{\nu\circ\alpha^{-1}\circ T_1}(x_2))=\\
=(J_{\nu\circ\alpha^{-1}\circ T_1}\alpha\circ\beta^{-1}(T_2(y_1^*y_2))J_{\nu\circ\alpha^{-1}\circ T_1}\Lambda_{\nu\circ\alpha^{-1}\circ T_1}(x_1)|\Lambda_{\nu\circ\alpha^{-1}\circ T_1}(x_2))
\end{multline*}
and, using (i) again, is equal to :
\begin{eqnarray*}
(\Lambda_{T_1}(x_2)J_\nu\Lambda_\nu(\beta^{-1}T_2(y_1^*y_2))|\Lambda_{\nu\circ\alpha^{-1}\circ T_1}(x_2))
&=&
(J_\nu\Lambda_\nu(\beta^{-1}T_2(y_1^*y_2))|\Lambda_{T_1}(x_2)^*\Lambda_{\nu\circ\alpha^{-1}\circ T_1}(x_2))\\
&=&(J_\nu\Lambda_\nu(\beta^{-1}T_2(y_1^*y_2))|\Lambda_\nu(\alpha^{-1}T_1(x_2^*x_1)))
\end{eqnarray*}
This finishes the proof.  \end{proof}

%%%fiber
\subsection{Fiber product of von Neumann algebras \cite{EVal}}
\label{fiber}
If $x\in \beta(N)'$ and $y\in \alpha(N)'$, it is possible to define an
operator $x\underset{\nu}{_\beta\otimes_\alpha}y$ on $\mathcal
K\underset{\nu}{_\beta\otimes_\alpha}\mathcal H$, with natural values
on the elementary tensors. As this operator does not depend upon the
weight $\nu$, it will be denoted by $x\underset{N}{_\beta\otimes_\alpha}y$. 

If $P$ is a von Neumann algebra on $\mathcal H$ with
$\alpha(N)\subset P$, and $Q$ a von Neumann algebra on $\mathcal K$
with $\beta(N)\subset Q$, then we define the {\it fiber product}
$Q\underset{N}{_\beta*_\alpha}P$ as
$\{x\underset{N}{_\beta\otimes_\alpha}y : x\in Q', y\in P'\}'$. 
 This von Neumann algebra can be defined independently of the Hilbert
 spaces on which $P$ and $Q$ are represented. If for $i=1,2$,
 $\alpha_i$ is a faithful non-degenerate homomorphism from $N$ into
 $P_i$, and $\beta_i$ is a faithful non-degenerate anti-homomorphism from $N$ into $Q_i$, and $\Phi$ (resp. $\Psi$) a homomorphism from $P_1$ to $P_2$ (resp. from $Q_1$ to $Q_2$) such that $\Phi\circ\alpha_1=\alpha_2$ (resp. $\Psi\circ\beta_1=\beta_2$), then, it is possible to define a homomorphism $\Psi\underset{N}{_{\beta_1}*_{\alpha_1}}\Phi$ from $Q_1\underset{N}{_{\beta_1}*_{\alpha_1}}P_1$ into $Q_2\underset{N}{_{\beta_2}*_{\alpha_2}}P_2$. 

We define a relative flip $\varsigma_N$ from $B(\mathcal K)\underset{N}{_\beta*_\alpha}B(\mathcal H)$ onto $B(\mathcal H)\underset{N^{\op}}{_\alpha*_\beta}B(\mathcal K)$ by $\varsigma_N(X)=\sigma_\psi X(\sigma_{\psi})^*$ for any $X\in B(\mathcal K)\underset{N}{_\beta*_\alpha}B(\mathcal H)$ and any normal semi-finite faithful weight $\psi$ on $N$.

Let now $U$ be an isometry from a Hilbert space $\mathcal K_1$ in a Hilbert space $\mathcal K_2$, which intertwines two anti-representations $\beta_1$ and $\beta_2$ of $N$, and let $V$ be an isometry from a Hilbert space $\mathcal H_1$ in a Hilbert space $\mathcal H_2$, which intertwines two representations $\alpha_1$ and $\alpha_2$ of $N$. Then, it is possible to define, on linear combinations of elementary tensors, an isometry $U\underset{\nu}{_{\beta_1}\otimes_
{\alpha_1}}V$ which can be extended to the whole Hilbert space
$\mathcal K_1\underset{\nu}{_{\beta_1}\otimes_{\alpha_1}}\mathcal
H_1$ with values in $\mathcal
K_2\underset{\nu}{_{\beta_2}\otimes_{\alpha_2}}\mathcal H_2$. One can
show that this isometry does not depend upon the weight $\nu$. It
will be denoted by $U\underset{N}{_{\beta_1}\otimes_{\alpha_1}}V$. If $U$ and $V$ are unitaries, then $U\underset{N}{_{\beta_1}\otimes_{\alpha_1}}V$ is an unitary and $(U\underset{N}{_{\beta_1}\otimes_{\alpha_1}}V)^*=U^*\underset{N}{_{\beta_2}\otimes_{\alpha_2}}V^*$.

If $x\in D(\sigma_{-i/2}^\nu)$, then it is possible to construct on elementary tensors an operator $\beta(x)\underset{\nu}{_\beta\otimes_\alpha}1=1\underset{\nu}{_\beta\otimes_\alpha}\alpha(\sigma_{-i/2}^\nu(x))$ (\cite{S}, 2.2b)
%%%%%Hopf bimodule
\subsection{Definition of a Hopf-bimodule}
\label{defHopf}
 A quintuple $(N, M, \alpha, \beta, \Gamma)$ will be called a \emph{Hopf-bimodule}, following (\cite{Val2}, \cite{EVal} 6.5), if
$N$,
$M$ are von Neumann algebras, $\alpha$ is a faithful non-degenerate
representation of $N$ into $M$, $\beta$ is a
faithful non-degenerate anti-representation of
$N$ into $M$, with commuting ranges, and $\Gamma$ is an injective $*$-homomorphism from $M$
into
$M\underset{N}{_\beta *_\alpha}M$ such that, for all $X$ in $N$ :
\begin{enumerate}
\item $\Gamma
  (\beta(X))=1\underset{N}{_\beta\otimes_\alpha}\beta(X)$,
\item  $\Gamma
  (\alpha(X))=\alpha(X)\underset{N}{_\beta\otimes_\alpha}1$,
\item  $\Gamma$ satisfies the co-associativity relation :
  \[(\Gamma \underset{N}{_\beta *_\alpha}\id)\Gamma =(\id
  \underset{N}{_\beta *_\alpha}\Gamma)\Gamma\]
\end{enumerate}
This last formula makes sense, thanks to the two preceeding ones and
\ref{spatial}. The von Neumann algebra $N$ will be called the \emph{basis} of $(N, M, \alpha, \beta, \Gamma)$. 

In (\cite{DC1}, chap. 11), De Commer has shown that, if $N$ is
finite-dimensional, the Hilbert space
$L^2(M)\underset{\nu}{_\beta\otimes_\alpha}L^2(M)$ can be
isometrically imbedded into the usual Hilbert tensor product
$L^2(M)\otimes L^2(M)$ and the projection $p$ on this closed subspace
belongs to $M\otimes M$. Moreover, the fiber product
$M\underset{N}{_\beta*_\alpha}M$ can be then identified with the
reduced von Neumann algebra $p(M\otimes M)p$ and we can consider $\Gamma$ as an usual coproduct $M\mapsto M\otimes M$, but with the condition $\Gamma(1)=p$. 

A \emph{co-inverse} $R$ for a Hopf bimodule $(N, M, \alpha, \beta,
\Gamma)$ is an involutive ($R^2=\id$) anti-$*$-isomorphism of $M$
satisfying $R\circ\alpha=\beta$ (and therefore $R\circ\beta=\alpha$)
and $\Gamma\circ R=\varsigma_{N^{\op}}\circ
(R\underset{N}{_\beta*_\alpha}R)\circ \Gamma$, where
$\varsigma_{N^{\op}}$ is the flip from
$M\underset{N^{\op}}{_\alpha*_\beta}M$ onto
$M\underset{N}{_\beta*_\alpha}M$. A Hopf bimodule is called \emph{co-commutative} if $N$ is abelian, $\beta=\alpha$, and $\Gamma=\varsigma\circ\Gamma$. 

For an example, suppose that $\mathcal G$ is a measured groupoid, with
$\mathcal G^{(0)}$ as its set of units. We denote
by $r$ and $s$ the range and source maps from $\mathcal G$ to $\mathcal G^{(0)}$, given by
$xx^{-1}=r(x)$ and $x^{-1}x=s(x)$, and  by $\mathcal G^{(2)}$ the set of composable elements, i.e.\ 
\[\mathcal G^{(2)}=\{(x,y)\in \mathcal G^2 : s(x)=r(y)\}\]
Let $(\lambda^u)_{u\in \mathcal G^{(0)}}$ be a Haar system on
$\mathcal G$ and $\nu$ a measure on $\mathcal G^{(0)}$. Let us denote by $\mu$ the measure on $\mathcal G$ given by integrating $\lambda^u$ by $\nu$,
\[\mu=\int_{{\mathcal G}^{(0)}}\lambda^ud\nu\]
By definition,  $\nu$ is called \emph{quasi-invariant} if $\mu$ is equivalent to its image under the inversion $x\mapsto x^{-1}$ of $\mathcal G$ (see \cite{R},
 \cite{C2} II.5, \cite{Pa} and \cite{AR} for more details, precise definitions and examples of groupoids). 

In \cite{Y1}, \cite{Y2}, \cite{Y3} and \cite{Val2} was associated to a measured groupoid $\mathcal G$, equipped with a Haar system $(\lambda^u)_{u\in \mathcal G ^{(0)}}$ and a quasi-invariant measure $\nu$ on $\mathcal G ^{(0)}$, a
Hopf bimodule with an abelian underlying von Neumann algebra $(L^\infty (\mathcal G^{(0)}, \nu), L^\infty (\mathcal G, \mu), r_{\mathcal G}, s_{\mathcal G}, \Gamma_{\mathcal
G})$, where 
$r_{\mathcal G}(g)=g\circ r$ and $s_{\mathcal G}(g)=g\circ s$ for all
$g$ in $L^\infty (\mathcal G^{(0)})$
 and where
$\Gamma_{\mathcal G}(f)$, for $f$ in $L^\infty (\mathcal G)$, is the
function defined on $\mathcal G^{(2)}$ by $(s,t)\mapsto f(st)$. Thus,
$\Gamma_{\mathcal G}$ is an involutive homomorphism from $L^\infty (\mathcal G)$ into $L^\infty
(\mathcal G^{(2)})$, which can be identified with
$L^\infty (\mathcal G){_s*_r}L^\infty (\mathcal G)$.

It is straightforward to get that the inversion of the groupoid gives a co-inverse for this Hopf bimodule structure.

%%%%defMQG
\subsection{Definition of measured quantum groupoids (\cite{L}, \cite{E4})}
\label{defMQG}

A \emph{measured quantum groupoid} is an octuple $\mathfrak {G}=(N, M, \alpha, \beta, \Gamma, T, T', \nu)$ such that (\cite{E4}, 3.8):

(i) $(N, M, \alpha, \beta, \Gamma)$ is a Hopf bimodule, 

(ii) $T$ is a left-invariant normal, semi-finite, faithful operator-valued weight from $M$ to $\alpha (N)$ (to be more precise, from $M^{+}$ to the extended positive elements of $\alpha(N)$ (cf. \cite{T} IX.4.12)), which means that, for any $x\in\gM_T^+$, we have $(\id\underset{\nu}{_\beta*_\alpha}T)\Gamma(x)=T(x)\underset{N}{_\beta\otimes_\alpha}1$. 

(iii) $T'$ is a right-invariant normal, semi-finite, faithful operator-valued weight from $M$ to $\beta (N)$, which means that, for any $x\in\gM_{T'}^+$, we have $(T'\underset{\nu}{_\beta*_\alpha}\id)\Gamma(x)=1\underset{N}{_\beta\otimes_\alpha}T'(x)$. 

(iv) $\nu$ is normal semi-finite faithful weight on $N$, which is relatively invariant with respect to $T$ and $T'$, which means that the modular automorphisms groups of the weights $\varphi=\nu\circ\alpha^{-1}\circ T$ and $\psi=\nu\circ\beta^{-1}\circ T'$ commute. The weight $\varphi$ will be called left-invariant, and $\psi$ right-invariant.

For example, 
let $\mathcal G$ be a measured groupoid equipped with a left Haar
system $(\lambda^u)_{u\in\mathcal G^{(0)}}$ and a quasi-invariant
measure $\nu$ on $\mathcal G^{(0)}$. Let us use the notations
introduced in \ref{defHopf}. If $f\in L^\infty(\mathcal G, \mu)^+$,
consider the function on $\mathcal G^{(0)}$, $u\mapsto \int_{\mathcal
  G}fd\lambda^u$, which belongs to $L^\infty (\mathcal G^{(0)},
\nu)$. The image of this function by the homomorphism $r_\mathcal G$
is the function on $\mathcal G$, $\gamma\mapsto \int_{\mathcal
  G}fd\lambda^{r(\gamma)}$, and the application which sends $f$ to
this function can be considered as an operator-valued weight from
$L^\infty(\mathcal G, \mu)$ to $r_{\mathcal G}(L^\infty (\mathcal
G^{(0)}, \nu))$ which is normal, semi-finite and faithful.  By
definition of the Haar system $(\lambda^u)_{u\in\mathcal G^{(0)}}$, it
is left-invariant in the sense of (ii). We shall denote this operator-valued weight from $L^\infty(\mathcal G, \mu)$ to $r_{\mathcal G}(L^\infty (\mathcal G^{(0)}, \nu))$ by $T_\mathcal G$.
 If we write $\lambda_u$ for the image of $\lambda^u$ under the inversion $x\mapsto x^{-1}$ of the groupoid $\mathcal G$, starting from the application which sends $f$ to the function on $\mathcal G^{(0)}$ defined by $u\mapsto \int_{\mathcal G}fd\lambda_u$, we define a normal semifinite faithful operator-valued weight from $L^\infty(\mathcal G, \mu)$ to $s_{\mathcal G}(L^\infty (\mathcal G^{(0)}, \nu))$, which is right-invariant in the sense of (ii),  and which we shall denote by $T_{\mathcal G}^{(-1)}$.
 
  We then get that :
  \[(L^\infty(\mathcal G^{(0)}, \nu), L^\infty(\mathcal G, \mu), r_\mathcal G, s_\mathcal G, \Gamma_\mathcal G, T_\mathcal G, T_{\mathcal G}^{(-1)}, \nu)\] is a measured quantum groupoid, which we shall denote again $\mathcal G$. 

It can be proved (\cite{E5}) that any measured quantum groupoid, whose underlying von Neumann algebra is abelian, is of that type. 

Let $\gG=(N, M, \alpha, \beta, \Gamma, T, T', \nu)$ be a measured quantum groupoid, then we denote by $\gG^o$ the octuplet $(N^o, M, \beta, \alpha, \varsigma_N\Gamma, T', T, \nu^o)$ (where $\sigma_N$ is the flip from $M\underset{N}{_\beta*_\alpha}M$ onto $M\underset{N^o}{_\alpha*_\beta}M$); it is is another measured quantum groupoid, called the \emph{opposite measured quantum groupoid} of $\gG$.

If $T$ is bounded, $\gG$ is called "of compact type". 
%%%%%defW
\subsection{Pseudo-multiplicative unitary.} 
\label{defW}
Let $\mathfrak{G}=(N, M, \alpha, \beta, \Gamma, T, T', \nu)$ be an
octuple satisfying the axioms (i), (ii) (iii) of \ref{defMQG}. With $\varphi=\nu\circ\alpha^{-1}\circ T$, we
shall write $H=H_\varphi$, $J=J_\varphi$ and $\widehat{\beta}(n)=J\alpha(n^*)J$ for
all $n\in N$.

Then (\cite{L}, 3.7.3 and 3.7.4), $\mathfrak {G}$ can be equipped with
a \emph{pseudo-multiplicative unitary} $W$ which is a unitary from
$H\underset{\nu}{_\beta\otimes_\alpha}H$ onto
$H\underset{\nu^{\op}}{_\alpha\otimes_{\widehat{\beta}}}H$ (\cite{E4}, 3.6)
that intertwines $\alpha$, $\widehat{\beta}$, $\beta$  in the following way:
for all $X\in N$, we have :
\begin{align*}
W(\alpha(X)\underset{N}{_\beta\otimes_\alpha}1)&=
(1\underset{N^{\op}}{_\alpha\otimes_{\widehat{\beta}}}\alpha(X))W, \\
W(1\underset{N}{_\beta\otimes_\alpha}\beta(X)) &=(1\underset{N^{\op}}{_\alpha\otimes_{\widehat{\beta}}}\beta (X))W, \\
W(\widehat{\beta}(X) \underset{N}{_\beta\otimes_\alpha}1) &=(\widehat{\beta}(X)\underset{N^{\op}}{_\alpha\otimes_{\widehat{\beta}}}1)W, \\
W(1\underset{N}{_\beta\otimes_\alpha}\widehat{\beta}(X))&=(\beta(X)\underset{N^{\op}}{_\alpha\otimes_{\widehat{\beta}}}1)W.
\end{align*}
Moreover, the operator $W$ satisfies the \emph{pentagonal relation} :
\[(1\underset{N^{\op}}{_\alpha\otimes_{\widehat{\beta}}}W)
(W\underset{N}{_\beta\otimes_\alpha}1_{H})
=(W\underset{N^{\op}}{_\alpha\otimes_{\widehat{\beta}}}1)
\sigma^{23}_{\alpha, \beta}(W\underset{N}{_{\widehat{\beta}}\otimes_\alpha}1)
(1\underset{N}{_\beta\otimes_\alpha}\sigma_{\nu^{\op}})
(1\underset{N}{_\beta\otimes_\alpha}W)\]
where $\sigma^{23}_{\alpha, \beta}$
goes from $(H\underset{\nu^{\op}}{_\alpha\otimes_{\widehat{\beta}}}H)\underset{\nu}{_\beta\otimes_\alpha}H$ to $(H\underset{\nu}{_\beta\otimes_\alpha}H)\underset{\nu^{\op}}{_\alpha\otimes_{\widehat{\beta}}}H$, 
and $1\underset{N}{_\beta\otimes_\alpha}\sigma_{\nu^{\op}}$ goes from $H\underset{\nu}{_\beta\otimes_\alpha}(H\underset{\nu^{\op}}{_\alpha\otimes_{\widehat{\beta}}}H)$ to $H\underset{\nu}{_\beta\otimes_\alpha}H\underset{\nu}{_{\widehat{\beta}}\otimes_\alpha}H$. 
The operators in this formula are well-defined because of the
intertwining relations listed above.

The operator $W$ is defined by the following formula, for any $a\in\gN_T\cap\gN_\varphi$, $v\in D((H_\varphi)_\beta, \nu^o)$, where $(\xi_i)_{i\in I}$ is a $(\beta, \nu^o)$ basis of $H_\varphi$ (in the sense of \ref{basis} (\cite{L}, 3.2.10)) :
\[W^*(v\underset{\nu^o}{_\alpha\otimes_{\widehat{\beta}}}\Lambda_\varphi(a))=\sum_{i\in I}\xi_i\underset{\nu}{_\beta\otimes_\alpha}\Lambda_\varphi((\omega_{v, \xi_i}\underset{\nu}{_\beta*_\alpha}id)\Gamma(a)))\]
The operator $W$ does not depend of the choice of the $(\beta, \nu^o)$ basis. 
 Moreover, $W$, $M$ and $\Gamma$ are related by the following results:

(i) $M$ is the weakly closed linear space generated by all operators $(\id*\omega_{\xi, \eta})(W)$, where $\xi\in D(_\alpha H, \nu)$ and $\eta\in D(H_{\widehat{\beta}}, \nu^{\op})$ (see \cite{E4}, 3.8(vii)).

(ii)  $\Gamma(x)=W^*(1\underset{N^{\op}}{_\alpha\otimes_{\widehat{\beta}}}x)W$
for all $x\in M$(\cite{E4}, 3.6). 

(iii) For any $x$, $y_1$, $y_2$ in $\gN_T\cap\gN_\varphi$, we have (\cite{E4}, 3.6):
\[(\id*\omega_{J\Lambda_\varphi (y_1^*y_2), \Lambda_\varphi (x)})(W)=
(\id\underset{N}{_\beta*_\alpha}\omega_{J\Lambda_\varphi(y_2), J\Lambda_\varphi(y_1)})\Gamma (x^*)\]

(iv) for any $a$ in $\gN_T\cap\gN_\varphi$, $v$ in $D(_\alpha H, \nu)\cap D(H_\beta, \nu^o)$, $w$ in $D(H_\beta, \nu^o)$, we have (\cite{L}, 3.3.3):
\[(\omega_{v,w}*id)(W^*)\Lambda_\varphi (a)=\Lambda_\varphi((\omega_{v,w}\underset{\nu}{_\beta*_\alpha}id)\Gamma(a))\]

If $N$ is finite-dimensional, using the fact that the relative tensor products can be identified with closed subspaces of the usual Hilbert tensor product (\ref{spatial}), we get that $W$ can be considered as a partial isometry on the usual Hilbert tensor product, which is multiplicative in the usual sense (i.e.\ such that $W_{23}W_{12}=W_{12}W_{13}W_{23}$.)

%%%%%%%lem1
\subsection{Lemma}
\label{lem1}
{\it Let $W$ be an $(\alpha, \hat{\beta}, \beta)$-pseudo-multiplicative unitary, $\xi_1$ in $D(\gH_\beta, \nu^o)$, $\xi_2$ in $D(_\alpha \gH, \nu)$, $\eta$ in $\gH$; let $\zeta_i$ in $D(\gH_\beta, \nu^o)$ and $\zeta'_i$ in $\gH$ such that $W^*(\xi_2\underset{\nu^o}{_\alpha\otimes_{\hat{\beta}}}\eta)=\sum_i \zeta_i\underset{\nu}{_\beta\otimes_\alpha}\zeta'_i$; then we have :}
\[\sum_i\alpha(<\zeta_i, \xi_1>_{\beta, \nu^o})\zeta'_i=(\omega_{\xi_1, \xi_2}*id)(W)^*\eta\]

\begin{proof}
Let $\theta$ in $\gH$; we have :
\begin{eqnarray*}
((\omega_{\xi_1, \xi_2}*id)(W)^*\eta|\theta)
&=&(W^*(\xi_2\underset{\nu^o}{_\alpha\otimes_{\hat{\beta}}}\eta)|\xi_1\underset{\nu}{_\beta\otimes_\alpha}\theta)\\
&=&(\sum_i \zeta_i\underset{\nu}{_\beta\otimes_\alpha}\zeta'_i|\xi_1\underset{\nu}{_\beta\otimes_\alpha}\theta)\\
&=&(\sum_i\alpha(<\zeta_i, \xi_1>_{\beta, \nu^o})\zeta'_i|\theta)
\end{eqnarray*}
from which we get the result.  \end{proof}

%%%%%%%lem2
\subsection{Lemma}
\label{lem2}
{\it Let $W$ be an $(\alpha, \hat{\beta}, \beta)$-pseudo-multiplicative unitary, $\xi_1$, $\zeta_1$ in $D(\gH_\beta, \nu^o)$, $\xi$ in $D(_\alpha\gH, \nu)$ and $\eta_1$, $\eta_2$ in $\gH$. Let us consider the flip $\sigma^{1,2}_{\hat{\beta}, \alpha}$ from $H\underset{\nu}{_\beta\otimes_\alpha}(H\underset{\nu^o}{_\alpha\otimes_{\hat{\beta}}}H)$ onto $H\underset{\nu^o}{_\alpha\otimes_{\hat{\beta}}}(H\underset{\nu}{_\beta\otimes_\alpha}H)$. Then, we have :}
\begin{multline*}
(\sigma^{1,2}_{\hat{\beta}, \alpha}
(1_{\gH}\underset{N}{_\beta\otimes_\alpha}W)(\xi_1\underset{\nu}{_\beta\otimes_\alpha}\eta_1\underset{\nu}{_\beta\otimes_\alpha}\xi)|\eta_2\underset{\nu^o}{_\alpha\otimes_{\hat{\beta}}}(\zeta_1\underset{\nu}{_\beta\otimes_\alpha}\zeta_2))=\\
(W(\eta_1\underset{\nu}{_\beta\otimes_\alpha}\xi)|\eta_2\underset{\nu^o}{_\alpha\otimes_{\hat{\beta}}}\alpha(<\zeta_1, \xi_1>_{\beta, \nu^o})\zeta_2)
\end{multline*}
\begin{proof}
The scalar product 
\[(\sigma^{1,2}_{\hat{\beta}, \alpha}
(1_{\gH}\underset{N}{_\beta\otimes_\alpha}W)(\xi_1\underset{\nu}{_\beta\otimes_\alpha}\eta_1\underset{\nu}{_\beta\otimes_\alpha}\xi)|\eta_2\underset{\nu^o}{_\alpha\otimes_{\hat{\beta}}}(\zeta_1\underset{\nu}{_\beta\otimes_\alpha}\zeta_2))\]
is equal to :
\[(\xi_1\underset{\nu}{_\beta\otimes_\alpha}W(\eta_1\underset{\nu}{_\beta\otimes_\alpha}\xi)|
\zeta_1\underset{\nu}{_\beta\otimes_\alpha}(\eta_2\underset{\nu^o}{_\alpha\otimes_{\hat{\beta}}}\zeta_2))\]
from which we get the result. \end{proof}

%%%%%prop2Gamma
\subsubsection{\bf{Proposition}}
\label{prop2Gamma}
{\it $\gG=(N, M, \alpha, \beta, \Gamma, T, T', \nu)$ is a measured quantum groupoid in the sense of \ref{defMQG}, Let $W$ be its pseudo-multiplicative unitary,, $\xi$ in $D(_\alpha\gH, \nu)$, $\eta$ in $D(\gH_{\hat{\beta}}, \nu^o)$. Let $\xi_1$, $\eta_1$ in $D(\gH_\beta, \nu^o)$, $\xi_2$, $\eta_2$ in $D(_\alpha\gH, \nu)$; then, we have :}
\[(\Gamma((id*\omega_{\xi, \eta})(W))(\xi_1\underset{\nu}{_\beta\otimes_\alpha}\eta_1)|
\xi_2\underset{\nu}{_\beta\otimes_\alpha}\eta_2)=
((\omega_{\xi_1, \xi_2}*id)(W)(\omega_{\eta_1, \eta_2}*id)(W)\xi|\eta)\]
\begin{proof}
Using the \ref{defW}(ii) we get that :
\[(\Gamma((id*\omega_{\xi, \eta})(W))(\xi_1\underset{\nu}{_\beta\otimes_\alpha}\eta_1)|
\xi_2\underset{\nu}{_\beta\otimes_\alpha}\eta_2)=
((1\underset{\nu^o}{_\alpha\otimes_{\hat{\beta}}}(id*\omega_{\xi, \eta})(W))W(\xi_1\underset{\nu}{_\beta\otimes_\alpha}\eta_1)|W(\xi_2\underset{\nu}{_\beta\otimes_\alpha}\eta_2))\]
which is equal to :
\[((1\underset{N^o}{_\alpha\otimes_{\hat{\beta}}}W)(W\underset{N}{_\beta\otimes_\alpha}1)(\xi_1\underset{\nu}{_\beta\otimes_\alpha}\eta_1\underset{\nu}{_\beta\otimes_\alpha}\xi)|(W\underset{N^o}{_\alpha\otimes_{\hat{\beta}}}1)((\xi_2\underset{\nu}{_\beta\otimes_\alpha}\eta_2)\underset{\nu^o}{_\alpha\otimes_{\hat{\beta}}}\eta)\]
which, using the pentagonal relation (\ref{defW}), is equal to :
\[(\sigma^{2,3}_{\alpha, \beta}(W\underset{N}{_{\hat{\beta}}\otimes_\alpha}1)
(1_{\gH}\underset{N}{_\beta\otimes_\alpha}\sigma_{\nu^o})
(1_{\gH}\underset{N}{_\beta\otimes_\alpha}W)(\xi_1\underset{\nu}{_\beta\otimes_\alpha}\eta_1\underset{\nu}{_\beta\otimes_\alpha}\xi)|
(\xi_2\underset{\nu}{_\beta\otimes_\alpha}\eta_2)\underset{\nu^o}{_\alpha\otimes_{\hat{\beta}}}\eta)\] 
or, to :
\[((W\underset{N}{_{\hat{\beta}}\otimes_\alpha}1)
(1_{\gH}\underset{N}{_\beta\otimes_\alpha}\sigma_{\nu^o})
(1_{\gH}\underset{N}{_\beta\otimes_\alpha}W)(\xi_1\underset{\nu}{_\beta\otimes_\alpha}\eta_1\underset{\nu}{_\beta\otimes_\alpha}\xi)|(\xi_2\underset{\nu^o}{_\alpha\otimes_{\hat{\beta}}}\eta)\underset{\nu}{_\beta\otimes_\alpha}\eta_2)\]
which is equal to :
\[(\sigma^{1,2}_{\hat{\beta}, \alpha}
(1_{\gH}\underset{N}{_\beta\otimes_\alpha}W)(\xi_1\underset{\nu}{_\beta\otimes_\alpha}\eta_1\underset{\nu}{_\beta\otimes_\alpha}\xi)|\eta_2\underset{\nu^o}{_\alpha\otimes_{\hat{\beta}}}(W^*(\xi_2\underset{\nu}{_\alpha\otimes_{\hat{\beta}}}\eta)))\]
Defining now $\zeta_i$, $\zeta'_i$ as in \ref{lem1}, we get, using \ref{lem2}, that it is equal to :
\[(W(\eta_1\underset{\nu}{_\beta\otimes_\alpha}\xi)|\eta_2\underset{\nu^o}{_\alpha\otimes_{\hat{\beta}}}\sum_i\alpha(<\zeta_i, \xi_1>_{\beta, \nu^o})\zeta'_i)\]
which, thanks to \ref{lem1}, is equal to :
\[(W(\eta_1\underset{\nu}{_\beta\otimes_\alpha}\xi)|\eta_2\underset{\nu^o}{_\alpha\otimes_{\hat{\beta}}}(\omega_{\xi_1, \xi_2}*id)(W)^*\eta)\]
and, therefore, to :
\[((\omega_{\eta_1, \eta_2}*id)(W)\xi|(\omega_{\xi_1, \xi_2}*id)(W)^*\eta)\]
which finishes the proof.  \end{proof}

%%%%%%data
\subsection{Other data associated to a measured quantum groupoid (\cite{L}, \cite{E4})}
\label{data}
 Suppose that $\gG=(N, M, \alpha, \beta, \Gamma, T, T', \nu)$ is a measured quantum groupoid in the sense of \ref{defMQG}. Let us write $\varphi=\nu\circ\alpha^{-1}\circ T$, which is a normal semi-finite faithful left-invariant weight on $M$. Then:

(i) There exists an anti-$*$-automorphism $R$ on $M$ such that :
\begin{align*}
  R^2&=\id, & R(\alpha(n))&=\beta(n) \text{ for all } n\in N, &
  \Gamma\circ
  R&=\varsigma_{N^{\op}}(R\underset{N}{_\beta*_\alpha}R)\Gamma\end{align*}
and :
\begin{align*}
  R((\id*\omega_{\xi, \eta})(W))=(\id*\omega_{J\eta, J\xi})(W) \quad \text{for all } \xi\in D(_\alpha H, \nu), \eta\in D(H_{\widehat{\beta}}, \nu^{\op}).
\end{align*}
This map $R$ will be called the \emph{co-inverse}.

(ii) There exists a one-parameter group $\tau_t$ of automorphisms of $M$ such that :
\begin{align*}
R\circ\tau_t&=\tau_t\circ R, &
\tau_t(\alpha(n))&=\alpha(\sigma^\nu_t(n)), & \tau_t(\beta(n))
&=\beta(\sigma^\nu_t(n)), &
  \Gamma\circ\sigma_t^\varphi &=(\tau_t\underset{N}{_\beta*_\alpha}\sigma_t^\varphi)\Gamma
\end{align*}
for all $t\in\R$ and and $n\in N$.
This one-parameter group will be called the \emph{scaling group}.

(iii) The weight $\nu$ is relatively invariant with respect to $T$ and
$RTR$. Moreover, $R$ and $\tau_t$ are still the co-inverse and the
scaling group of this new measured quantum groupoid, which we shall
denote  by :
\[\underline{\gG}=(N, M, \alpha, \beta, \Gamma, T, RTR, \nu), \]
and  for simplification we shall assume now that $T'=RTR$ and $\psi=\varphi\circ R$. 

(iv) There exists a one-parameter group $\gamma_t$ of automorphisms of $N$ such that :
\[\sigma_t^{T}(\beta(n))=\beta(\gamma_t(n))\]
 for all $t\in\R$ and $n\in N$.
Moreover, we get that $\nu\circ\gamma_t=\nu$. 

(v) There exist a positive non-singular operator $\lambda$ affiliated to $Z(M)$ and a positive non-singular operator $\delta$ affiliated with $M$ such that :
\[(D\varphi\circ R: D\varphi)_t=\lambda^{it^2/2}\delta^{it},\]
and therefore
\[(D\varphi\circ\sigma_s^{\varphi\circ R}:D\varphi)_t=\lambda^{ist}.\]
The operator $\lambda$ will be called the \emph{scaling operator}, and there exists a positive non-singular operator $q$ affiliated to $N$ such that $\lambda=\alpha(q)=\beta(q)$. We have $R(\lambda)=\lambda$. 

The operator $\delta$ will be called the \emph{modulus}. We have $R(\delta)=\delta^{-1}$ and $\tau_t(\delta)=\delta$ for all $t\in\R$, and we can define a one-parameter group of unitaries $\delta^{it}\underset{N}{_\beta\otimes_\alpha}\delta^{it}$ which acts naturally on elementary tensor products and satisfies for all $t\in\R$ :

\[\Gamma(\delta^{it})=\delta^{it}\underset{N}{_\beta\otimes_\alpha}\delta^{it}.\]

(vi) We have $(D\Phi\circ\tau_t: D\Phi)_s=\lambda^{-ist}$, which
proves that $\tau_t\circ\sigma_s^\varphi=\sigma_s^\varphi\circ\tau_t$  for all $s$, $t$ in $\R$ and allows to define a one-parameter group of unitaries by :
\[P^{it}\Lambda_\varphi(x)=\lambda^{t/2}\Lambda_\varphi(\tau_t(x)) \quad
\text{for all }  x\in\gN_\varphi.\]
Moreover, for any $y$ in $M$, we get that :
\[\tau_t(y)=P^{it}yP^{-it}.\]
 and it is possible to define one parameter groups of unitaries $P^{it}\underset{N}{_\beta\otimes_\alpha}P^{it}$ and $P^{it}\underset{N^o}{_\alpha\otimes_{\hat{\beta}}}P^{it}$ such that :
\[W(P^{it}\underset{N}{_\beta\otimes_\alpha}P^{it})=(P^{it}\underset{N^o}{_\alpha\otimes_{\hat{\beta}}}P^{it})W\]
Moreover, for all $v\in D(P^{-1/2})$, $w\in D(P^{1/2})$, $p$, $q$ in $D(_\alpha H_\Phi, \nu)\cap D((H_\Phi)_{\hat{\beta}}, \nu^o)$, we have
\[(W^*(v\underset{\nu^o}{_\alpha\otimes_{\hat{\beta}}}q)|w\underset{\nu}{_\beta\otimes_\alpha}p)=
(W(P^{-1/2}v\underset{\nu}{_\beta\otimes_\alpha}J_\Phi p)|P^{1/2}w\underset{\nu^o}{_\alpha\otimes_{\hat{\beta}}}J_\Phi q)\]
We shall say that the pseudo-multiplicative unitary $W$ is \emph{manageable}, with \emph{managing operator} $P$, which implies it is weakly regular in the sense of \cite{E3}, 4.1.
\newline
As $\tau_s\circ\sigma_t^\varphi=\sigma_t^\varphi\circ\tau_s$, we get that $J_\varphi PJ_\varphi=P$.
 
(vii) It is possible to construct a \emph{dual} measured quantum groupoid 
\[\widehat{\gG}=(N, \widehat{M}, \alpha, \widehat{\beta}, \widehat{\Gamma}, \widehat{T}, \widehat{T'}, \nu)\]
where $\widehat{M}$ is equal to the weakly closed linear space
generated by all operators of the form $(\omega_{\xi, \eta}*\id)(W)$,
for $\xi\in D( H_\beta, \nu^{\op})$ and $\eta\in D(_\alpha H, \nu)$,
 $\widehat{\Gamma}(y)=\sigma_{\nu^{\op}}
W(y\underset{N}{_\beta\otimes_\alpha}1)W^*\sigma_\nu$ for all
$y\in\widehat{M}$, and the dual left operator-valued weight $\widehat{T}$ is constructed in a similar way as the dual left-invariant weight of a locally compact quantum group. Namely, it is possible to construct a normal semi-finite faithful weight $\widehat{\varphi}$ on $\widehat{M}$ such that, for all $\xi\in D(H_\beta, \nu^{\op})$ and $\eta\in D(_\alpha H, \nu)$ such that $\omega_{\xi, \eta}$ belongs to $I_\varphi$ :
\[\widehat{\varphi}((\omega_{\xi, \eta}*\id)(W)^*(\omega_{\xi, \eta}*\id)(W))=\|\omega_{\xi, \eta}\|_\varphi^2.\]
We can prove  that
$\sigma_t^{\widehat{\varphi}}\circ\alpha=\alpha\circ\sigma_t^\nu$  for all $t\in\R$, which gives the existence of an operator-valued weight $\widehat{T}$, which appears then to be left-invariant. 

As the formula $y\mapsto J y^*J$ ($y\in \widehat{M}$) gives a
co-inverse for the coproduct $\widehat{\Gamma}$, we get also a
right-invariant operator-valued weight.  Moreover, the
pseudo-multiplicative unitary $\widehat{W}$ associated to
$\widehat{\gG}$ is $\widehat{W}=\sigma_\nu W^*\sigma_\nu$, its
managing operator $\widehat{P}$ is equal to $P$,  its scaling group
is given by $\widehat{\tau}_t(y)=P^{it}yP^{-it}$, its scaling operator
$\widehat{\lambda}$ is equal to $\lambda^{-1}$, and its one-parameter
group of automorphisms $\widehat{\gamma}_t$ of $N$ is equal to
$\gamma_{-t}$.

We write $\widehat{\varphi}$ for $\nu\circ\alpha^{-1}\circ\widehat{T}$,
identify $H_{\widehat{\varphi}}$ with $H$, and write
$\widehat{J}=J_{\widehat{\varphi}}$. Then
$R(x)=\widehat{J}x^*\widehat{J}$  for all $x\in M$ and $W^*=(\widehat{J}\underset{N^{\op}}{_\alpha\otimes_{\widehat{\beta}}}J)W(\widehat{J}\underset{N^{\op}}{_\alpha\otimes_{\widehat{\beta}}}J)$. 

Moreover, we have $\widehat{\widehat{\gG}}=\gG$. 

For example, let $\mathcal G$ be a measured groupoid as in \ref{defMQG}. The dual $\widehat{\mathcal G}$ of the measured quantum groupoid constructed in \ref{defMQG} (and denoted again by $\mathcal G$) is :
\[\widehat{\mathcal G}=(L^\infty(\mathcal G^{(0)}, \nu), \mathcal L(\mathcal G), r_{\mathcal G}, r_{\mathcal G}, \widehat{\Gamma}_{\mathcal G}, \widehat{T}_{\mathcal G}, \widehat{T}_{\mathcal G})\]
where $\mathcal L(\mathcal G)$ is the von Neumann algebra generated by
the convolution algebra associated to the groupoid $\mathcal G$, the
coproduct $\widehat{\Gamma}_{\mathcal G}$ had been defined in
(\cite{Val1} 3.3.2), and the operator-valued weight $\widehat{T}_{\mathcal G}$ had been defined in (\cite{Val1}, 3.3.4). The underlying Hopf-bimodule is co-commutative. 

The peudo-multiplicative unitary $W^o$ associated to the opposite measured quantum groupoid $\mathfrak G^o$ is :
\[W^o=(\widehat{J}\underset{\nu^o}{_\alpha\otimes_{\widehat{\beta}}}\widehat{J})W(\widehat{J}\underset{\nu^o}{_\alpha\otimes_\beta}\widehat{J})\]

%%%%lemW
\subsubsection{\bf Lemma}
\label{lemW}
 {\it Let $\gG=(N, M, \alpha, \beta, \Gamma, T, T', \nu)$ be a measured quantum groupoid, as defined in \ref{defMQG}, $W$ the pseudo-multiplicative unitary associated by \ref{defW} ;
 \newline
 (i) we have, using the notations of  \ref{defW} :
 \[(W^*\underset{N^o}{_\beta\otimes_\alpha}1)(1\underset{N}{_\alpha\otimes_{\widehat{\beta}}}W^*)(W\underset{N}{_\alpha\otimes_{\widehat{\beta}}}1)
 =
 (1\underset{N^o}{_\beta\otimes_\alpha}W^*)(1\underset{N}{_\beta\otimes_\alpha}\sigma_{\nu^o})^*(W^*\underset{N^o}{_{\widehat{\beta}}\otimes_\alpha}1)(\sigma^{2,3}_{\alpha, \beta})^*\]
(ii) let $a$, $b$ in $\gN_T\cap\gN_T^*\cap\gN_\varphi\cap\gN_\varphi^*$; we have :
\[\Gamma((id*\omega_{\Lambda_\varphi(b), J_\varphi\Lambda_\varphi(a)}(W^*))=
(\rho^{\beta, \alpha}_{J_\varphi\Lambda_\varphi(a)})^* (1\underset{N^o}{_\beta\otimes_\alpha}W^*)(1\underset{N}{_\beta\otimes_\alpha}\sigma_{\nu^o})^*(W^*\underset{N^o}{_{\widehat{\beta}}\otimes_\alpha}1)(\sigma^{2,3}_{\alpha, \beta})^*\rho^{\alpha, \widehat{\beta}}_{\Lambda_\varphi(b)}\]}

\begin{proof}
 Result (i) is easily obtained from \ref{defW}. Then, by \ref{defW}(ii), we get (ii). \end{proof}

%%%%AnW
\subsubsection{\bf Proposition}(\cite{E5}, 4.3).
\label{AnW}
{\it Let $\gG=(N, M, \alpha, \beta, \Gamma, T, T', \nu)$ be a measured quantum groupoid, as defined in \ref{defMQG}, $W$ the pseudo-mltiplicative unitary associated by \ref{defW}, $R$ the co-inverse associated by \ref{data}; let us define $A_n(W)$ as the norm closure of the linear span generated by all operators of the form $(id*\omega_{\xi, \eta})(W)$ for all $\xi\in D(_\alpha H, \nu), \eta\in D(H_{\widehat{\beta}}, \nu^{\op})$. Then $A_n(W)$ is an algebra, $A_n(W)\cap A_n(W)^*$ is a non degenerate sub-$\bf C^*$-algebra of $M$, weakly dense in $M$, invariant by $R$, $\sigma_t^\varphi$, $\sigma_t^{\varphi\circ R}$, $\tau_t$.}

\begin{proof} This is \cite{E3}, 3.6, \cite{E4}, 4.3 and 4.5. \end{proof}
%%%%%thL
\subsubsection{\bf Theorem}
\label{thL}
{\it Let $\gG=(N, M, \alpha, \beta, \Gamma, T, T', \nu)$ be a measured quantum groupoid; let's use all notations introduced in \ref{data}. Then, for any $\xi$, $\eta$ in $D(_\alpha H_\Phi, \nu)$, for all $t$ in $\mathbb{R}$, we have :
\newline
(i) $R((i*\omega_{\xi, J_\varphi\eta})(W))=(i*\omega_{\eta, J_\varphi\xi})(W)$; therefore $R(A_n(W))=A_n(W)$. 
\newline
(ii) $\tau_t((i*\omega_{\xi, J_\varphi\eta})(W))=(i*\omega_{\Delta_\varphi^{-it}\xi, \Delta_\varphi^{-it}J_\varphi\eta})(W)$
\newline
(iii)   $\sigma_t^\varphi((i*\omega_{\xi, J_\varphi\eta})(W))=
(i*\omega_{\delta^{it}J_\varphi\delta^{-it}J_\Phi\Delta_\varphi^{-it}\xi, P^{it}J_\varphi\eta})(W)$
\newline
$\sigma_t^{\varphi\circ R}((i*\omega_{\xi, J_\varphi\eta})(W))=
(i*\omega_{P^{it}\xi, \delta^{it}J_\varphi\delta^{-it}J_\varphi\Delta_\varphi^{-it}J_\varphi\eta})(W)$}
\begin{proof}
Results (i) and (ii) are (\cite{L} 4.6). 
\newline
Let us take $\xi=J_\varphi\Lambda_\varphi (y_1^*y_2)$, and $\eta=J_\varphi\Lambda_\varphi(x)$, with $x$, $y_1$, $y_2$ in  $\gN_T\cap\gN_\varphi$; then, using \ref{defW} and \ref{data}, we get :
\[\sigma_t^\varphi((i*\omega_{J_\varphi\Lambda_\varphi (y_1^*y_2), \Lambda_\varphi(x)})(W))
=(id\underset{N}{_\beta*_\alpha}\omega_{J_\varphi\Lambda_\varphi(y_2), J_\varphi\Lambda_\varphi(y_1)}\circ\sigma_t^{\varphi\circ R})\Gamma (\tau_t(x^*))\]
which is equal to :
\begin{multline*}
(id\underset{N}{_\beta*_\alpha}\omega_{J_\Phi\Lambda_\Phi(\lambda^{t/2}\sigma_{-t}^{\Phi\circ R}(y_2)), J_\Phi\Lambda_\Phi(\lambda^{t/2}\sigma_{-t}^{\Phi\circ R}(y_1))})\Gamma (\tau_t(x^*))\\
=(id\underset{N}{_\beta*_\alpha}\omega_{J_\Phi\Lambda_\Phi(\sigma_{-t}^{\Phi\circ R}(y_2)), J_\Phi\Lambda_\Phi(\sigma_{-t}^{\Phi\circ R}(y_1))})\Gamma (\lambda^t\tau_t(x^*))
\end{multline*}
which, using again \ref{defW} and \ref{data}, is equal to :
\[(i*\omega_{J_\varphi\Lambda_\varphi (\sigma_t^{\varphi\circ R}(y_1^*y_2)), \Lambda_\varphi(\lambda^t\tau_t(x))})(W)
=(i*\omega_{J_\varphi\delta^{-it}J_\varphi\delta^{it}J_\varphi\Delta_\varphi^{it}\Lambda_\varphi (y_1^*y_2), P^{it}\Lambda_\varphi(x)})(W)\]
which gives the first result of (iii), using \ref{defW}. 
\newline
By similar calculations, we obtain :
\[\sigma_t^{\varphi\circ R}((i*\omega_{J_\varphi\Lambda_\varphi (y_1^*y_2), \Lambda_\varphi(x)})(W))
=(id\underset{N}{_\beta*_\alpha}\omega_{J_\varphi\Lambda_\varphi(y_2), J_\varphi\Lambda_\varphi(y_1)}\circ\tau_t)\Gamma (\sigma_t^{\varphi\circ R}(x^*))\]
which is equal to :
\[(id\underset{N}{_\beta*_\alpha}\omega_{J_\varphi\Lambda_\varphi(\lambda^{t/2}\tau_t(y_2), J_\varphi\Lambda_\varphi(\lambda^{t/2}\tau_t(y_1)})\Gamma (\sigma_t^{\varphi\circ R}(x^*))
=(i*\omega_{J_\varphi\Lambda_\varphi (\lambda^t\tau_t(y_1^*y_2)), \Lambda_\varphi(\sigma^{\varphi\circ R}(x))})(W)\]
from which we obtain the second result of (iii). \end{proof}

%%%%lemW
\subsubsection{\bf Lemma}
\label{lemW}
{\it Let $\gG=(N, M, \alpha, \beta, \Gamma, T, T', \nu)$ be a measured quantum groupoid, as defined in \ref{defMQG}, $W$ the pseudo-multiplicative unitary associated by \ref{defW}, $R$ the co-inverse associated by \ref{data}; let $x$ and $y$ in $\gM_\varphi\cap\gM_T$; then, we have  :

(i) $(id*\omega_{J_\varphi\Lambda_\varphi(y), \Lambda_\varphi(x)})(W)^*=(id*\omega_{J_\varphi\Lambda_\varphi(y^*), \Lambda_\varphi (x^*)})(W)$

(ii) $R[(id*\omega_{J_\varphi\Lambda_\varphi(y), \Lambda_\varphi(x)})(W)]=(id*\omega_{J_\varphi\Lambda_\varphi(x), \Lambda_\varphi(y)})(W)$}
\begin{proof}
Using \ref{defW}(iii), we have, for any $x_1$, $x_2$, $y_1$, $y_2$ in $\gN_\varphi\cap\gN_T$ : 
\[(id*\omega_{J_\varphi\Lambda_\varphi(y_2^*y_1), \Lambda_\varphi(x_2^*x_1)})(W)=(id\underset{N}{_\beta*_\alpha}\omega_{J_\varphi\Lambda_\varphi(y_1), J_\varphi\Lambda_\varphi(y_2)})\Gamma(x_1^*x_2)\]
from which we get :
\[(id*\omega_{J_\varphi\Lambda_\varphi(y_2^*y_1), \Lambda_\varphi(x_2^*x_1)})(W)^*=(id\underset{N}{_\beta*_\alpha}\omega_{J_\varphi\Lambda_\varphi(y_2), J_\varphi\Lambda_\varphi(y_1)})\Gamma(x_2^*x_1)\]
and, using \ref{defW}(iii) again, we get (i).

Using $W^*=(\widehat{J}\underset{N^o}{_\alpha\otimes_{\widehat{\beta}}}J)W(\widehat{J}\underset{N^o}{_\alpha\otimes_{\widehat{\beta}}}J)$ (\ref{data}), we get that :
\begin{eqnarray*}
R[(id*\omega_{J_\varphi\Lambda_\varphi(y), \Lambda_\varphi(x)})(W)]
&=&
\widehat{J}(id*\omega_{J_\varphi\Lambda_\varphi(y), \Lambda_\varphi(x)})(W)^*\widehat{J}\\
&=&
\widehat{J}(id*\omega_{J_\varphi\Lambda_\varphi(x), \Lambda_\varphi(y)})(W^*)\widehat{J}\\
&=&
\widehat{J}(id*\omega_{J_\varphi\Lambda_\varphi(x), \Lambda_\varphi(y)})[(\widehat{J}\underset{N^o}{_\alpha\otimes_{\widehat{\beta}}}J)W(\widehat{J}\underset{N^o}{_\alpha\otimes_{\widehat{\beta}}}J)]\widehat{J}\\
&=&
(id*\omega_{J_\varphi\Lambda_\varphi(x), \Lambda_\varphi(y)})(W)
\end{eqnarray*}
which is (ii). 

\end{proof}

%%%PropW
\subsubsection{\bf Proposition} (\cite{E5}, 4.6).
\label{PropW}
{\it Let $\gG=(N, M, \alpha, \beta, \Gamma, T, T', \nu)$ be a measured quantum groupoid, as defined in \ref{defMQG}, $W$ the pseudo-mltiplicative unitary associated by \ref{defW}, $R$ the co-inverse associated by \ref{data}; let $x_1$, $x_2$, $y_1$, $y_2$ in $\gN_\varphi\cap\gN_T$; then :

(i)  the operator $(id*\omega_{J_\varphi\Lambda_\varphi(y_2^*y_1), \Lambda_\varphi(x_2^*x_1)})(W)$ belongs to $\gM_T\cap\gM_\varphi$ if $y_2^*y_1$ belongs to $\gM_{RTR}$, and we have :
\begin{multline*}
T((id*\omega_{J_\varphi\Lambda_\varphi(y_2^*y_1), \Lambda_\varphi(x_2^*x_1)})(W)=\alpha(<RTR(y_2^*y_1)J_\varphi\Lambda_\varphi(x_2), J_\varphi\Lambda_\varphi(x_1)>_{\alpha, \nu}^o)=\\=\alpha(\Lambda_T(x_1)^*J_\varphi RTR(y_2^*y_1)J_\varphi\Lambda_T(x_2))
\end{multline*}

(ii) the operator $(id*\omega_{J_\varphi\Lambda_\varphi(y_2^*y_1), \Lambda_\varphi(x_2^*x_1)})(W)$ belongs to $\gM_{RTR}\cap\gM_{\varphi\circ R}$ if  $x_2^*x_1$ belongs to $\gM_{RTR}$, and we have then :
\begin{multline*}
RTR((id*\omega_{J_\varphi\Lambda_\varphi(y_2^*y_1), \Lambda_\varphi(x_2^*x_1)})(W)=\beta(<RTR(x_2^*x_1)J_\varphi\Lambda_\varphi(y_2), J_\varphi\Lambda_\varphi(y_1)>_{\alpha, \nu}^o)=\\=\beta(\Lambda_T(y_1)^*J_\varphi RTR(x_2^*x_1)J_\varphi\Lambda_T(y_2))
\end{multline*}}

\begin{proof}
Let $x$, $y$ in $\gN_T\cap\gN_\varphi$; Using \ref{defW}(iii), we have 
\[(id*\omega_{J_\varphi\Lambda_\varphi(y^*y), \Lambda_\varphi(x^*x)})(W)=(id\underset{N}{_\beta*_\alpha}\omega_{J_\varphi\Lambda_\varphi(y), J_\varphi\Lambda_\varphi(y)})\Gamma(x^*x)\]
(which is positive); applying the right-invariant operator valued weight $RTR$, we get :
\[RTR((id*\omega_{J_\varphi\Lambda_\varphi(y^*y), \Lambda_\varphi(x^*x)})(W))=(id\underset{N}{_\beta*_\alpha}\omega_{J_\varphi\Lambda_\varphi(y), J_\varphi\Lambda_\varphi(y)})(1\underset{N}{_\beta\otimes_\alpha}RTR(x^*x))\]
which is a bounded positive operator if $x$ belongs to $\gN_{RTR}$; it is then equal to :
\[\beta(<RTR(x^*x)J_\varphi\Lambda_\varphi(y), J_\varphi\Lambda_\varphi(y)>_{\alpha, \nu^o})=\beta(\Lambda_T(y)^*J_\varphi RTR(x^*x)J_\varphi\Lambda_T(y))\]
 by \ref{lemT}(i); then, we get (ii) by polarization and (i), by applying \ref{data}. \end{proof}

%%%%%LCQG
\section{Locally compact quantum groupoids}
\label{LCQG}

In this chapter, we first recall (\ref{C*fiber}) the definition and basic properties of the $\bf C^*$-relative tensor product (\ref{def3}, \ref{def4}) and the definition of the $\bf C^*$-fiber product of two $\bf C^*$-algebras, as defined by T. Timmermann (\cite{Ti2}). We then recall (\ref{weights}) the definition of a weight on a $\bf C^*$-algebra and the main properties  : \emph{lower semi-continuous weights} (\ref{lsc}), and \emph{KMS weights} (\ref{KMS}). In \ref{C*valuedweight}, we recall the definition of $\bf C^*$-valued weights and lower semi-continuous $\bf C^*$-valued weights, as introduced by J. Kustermans (\cite{K2}). We introduce (\ref{KMSC*}) the notion of a KMS pair $(\nu, T)$, where $\nu$ is a lower semi-continuous weight on a $\bf C^*$-algebra $B$, and $T$ a lower semi-continuous $\bf C^*$-valued weight from a $\bf C^*$-algebra $A$ to $M(B)$. We then give another definition of a fiber product of two $\bf C^*$-algebras (\ref{fiber}, \cite{Ti2}) and of a locally compact quantum groupoid (\ref{defLCQG}).

%%%%%C*fiber

\subsection{$\bf C^*$-relative product of Hilbert spaces and $\bf C^*$-fiber product of $\bf C^*$-algebras.} 
\label{C*fiber}
%%%%def1
\subsubsection{\bf{Definition}}
\label{def1}
A $\bf C^*$-\emph{base} $
\gb$ is a triple $\gb=(\mathfrak H, B, B^\dag)$ where $\mathfrak H$ is a Hilbert space, and $B$, $B^\dag$ are commuting non degenerate sub-$\bf C^*$-algebras of $B(\mathfrak H)$. We denote $\gb^\dag=(\mathfrak H, B^\dag, B)$.

%def2
\subsubsection{\bf{Definition}}
\label{def2}
Be given a $\bf C^*$-base $\gb=(\mathfrak H, B, B^\dag)$, a $\bf C^*$-$\gb$-\emph{module} is a pair $(H, L)$, where $H$ is a Hilbert space, and $L$ is a norm closed subspace of $B(\mathfrak H, H)$, such that $[L\mathfrak H]=H$, $[LB]=L$, $[L^*L]=B$, where $[X]$ means the closed linear space generated by $X$ (for $X\subset H$, or $X\subset B(\mathfrak H, H)$, or $X\subset B(\mathfrak H)$). 

Then (\cite{Ti4}, 2.5), there exists a non-degenerate normal representation $\rho_L$ of $B'$ on $H$, such that $\rho_L(x)L_1=L_1x$, for all $x \in B'$, $L_1\in L$. It is easy to check that, if $B$ is non degenerate on $H$, then $\rho_L$ is faithful. 

Be given a $\bf C^*$-base $\gb=(\mathfrak H, B, B^\dag)$ and two $\bf C^*$-$\gb$-modules $(H, L)$ and $(\tilde{H}, \tilde{L})$, a morphism from $(H, L)$ to $(\tilde{H}, \tilde{L})$,, is an operator $S\in B(H, \tilde{H})$, such that $SL\subset \tilde{L}$ and $S^*L\subset \tilde{L}$. Then, for all $b^\dag\in B^\dag$, we have $S\rho_L(b^\dag)=\rho_{\tilde{L}}(b^\dag)S$. 

%%%def3
\subsubsection{\bf{Definition}}
\label{def3}
Be given a $\bf C^*$-base $\gb=(\mathfrak H, B, B^\dag)$, a $\bf C^*$-$\gb$-module $(H, L)$ and a $\bf C^*$-$\gb ^\dag$-module $(K, P)$, the \emph{$\bf C^*$-relative tensor product} $H\underset{\gb}{_L\otimes_P}K$ is the Hilbert space generated by elements $P_1\otimes\xi\otimes L_1$, where $P_1\in P$, $\xi\in \mathfrak H$, $L_1\in L$, with the inner product :
\[(L_1\otimes\xi_1\otimes P_1|L_2\otimes\xi_2\otimes P_2)=(\xi_1|(P_1^*P_2)(L_1^*L_2)\xi_2)\]
In this formula, note that, by \ref{def2}, $P_1^*P_2$ belongs to $B^\dag$, $L_1^*L_2$ belongs to $B$, and therefore, by \ref{def1}, commute. 

The image of $L\otimes\xi\otimes P$ in $H\underset{\gb}{_L\otimes_P}K$ will be denoted by $L\vg\xi\vl P$.

If $S$ is a morphism from the $\bf C^*$-$\gb$-module $(H, L)$ to the $\bf C^*$-$\gb$-module $(\tilde{H}, \tilde{L})$, and $T$ is a morphism from the $\bf C^*$-$\gb ^\dag$-module $(K, P)$ to the $\bf C^*$-$\gb ^\dag$-module $(\tilde{K},\tilde{P})$, then the formula
\[(S\underset{\gb}{_L\otimes_P}T)(L_1\vg\xi\vl P_1)=SL_1\vg\xi\vl TP_1\]
defines $S\underset{\gb}{_L\otimes_P}T\in B(H\underset{\gb}{_L\otimes_P}K, \tilde{H}\underset{\gb}{_{\tilde{L}}\otimes_{\tilde{P}}\tilde{K}})$.

It is clear that we can define $\sigma_{\gb} : H\underset{\gb}{_L\otimes_P}K\mapsto K\underset{\gb^\dag}{_K\otimes_L}H$ by $\sigma_{\gb} (L\vg\xi\vl P)=P\vg\xi\vl L$, and its adjoint $\sigma_{\gb^\dag}$. For $X\in B(H\underset{\gb}{_L\otimes_P}K)$, we shall write $\varsigma X=\sigma_{\gb}X\sigma_{\gb^\dag}$, which belongs to $B(K\underset{\gb^\dag}{_K\otimes_L}H)$.
%%%%%threl
\subsubsection{\bf{Theorem}}
\label{threl}
{\it Be given a $\bf C^*$-base $\gb=(\mathfrak H, B, B^\dag)$, a $\bf C^*$-$\gb$-module $(H, L)$, a $\bf C^*$-$\gb ^\dag$-module $(K, P)$, and 
the relative tensor product $K\underset{\gb}{_P\otimes_L}H$, we have :
\newline
(i) For any $L_1\in L$, there exists an element $\lambda_{L_1}\in B(H, H\underset{\gb}{_{L}\otimes_P}K)$ such that, for any $P_1\in\  P$ and $\xi\in \mathfrak H$, we have $\lambda_{L_1}(P_1\xi)=L_1\vg\xi\vl P_1$.
\newline
(ii) For any $P_1\in P$, there exits an element $\rho_{P_1}\in B(K, H\underset{\gb}{_L\otimes_P}K)$ such that, for any $L_1\in L$ and $\xi\in\mathfrak H$, we have $\rho_{P_1}(L_1\xi)=L_1\vg\xi\vl P_1$.}
\begin{proof}
We have $\|L_1\vl\xi\vl P_1\|=\||L_1||P_1|\xi\|\leq\|L_1\|\||P_1|\xi\|=\|L_1\|\|P_1\xi\|$ from which we get (i), by density (\ref{def2}). Result (ii) is obtained the very same way. \end{proof}

%%%%%def4
\subsubsection{\bf{Definition}}
\label{def4}
Be given a $\bf C^*$-base $\gb=(\mathfrak H, B, B^\dag)$ and a $\bf C^*$-$\gb$-module $(H, L)$, a $\bf C^*$-$\gb$-\emph{algebra} on $(H,L)$ is a non degenerate sub-$\bf C^*$-algebra $A$ of $B(H)$, such that $\rho_L(B^\dag)\subset M(A)$, where $\rho_L$ had been defined in \ref{def2}. 

If $(\tilde{H}, \tilde{L})$ is another $\bf C^*$-$\gb$-module, and $\tilde{A}$ a $\bf C^*$-$\gb$-algebra on $(\tilde{H}, \tilde{L})$, a $\gb$-morphism $\Phi$ from $A$ to $\tilde{A}$ is a strict morphism of $\bf C^*$-algebras from $A$ to $M(\tilde{A})$, such that, for all $b^\dag\in B^\dag$ and $a\in A$, we have $\rho_{\tilde{L}}(b^\dag)\Phi(a)=\Phi(\rho_L(b^\dag)a)$. 

%%%%%%%def5
\subsubsection{\bf{Definition}}(\cite{Ti2}, 3.3)
\label{def5}
Be given a $\bf C^*$-base $\gb=(\mathfrak H, B, B^\dag)$, a $\bf C^*$-$\gb$-module $(H, L)$, a $\bf C^*$-$\gb ^\dag$-module $(K, P)$, and 
the relative tensor product $H\underset{\gb}{_{L}\otimes_P}K$, let $A_1$ be a a $\bf C^*$-$\gb$-algebra on $(H,L)$, and $A_2$ a $\bf C^*$-$\gb^\dag$-algebra on $(K,P)$; then the set of all elements $X\in B(H\underset{\gb}{_L\otimes_P}K)$ such that, for all $L_1\in L$, all operators $X\lambda_{L_1}$ and $X^*\lambda_{L_1}$ in $B(H, H\underset{\gb}{_L\otimes_P}K)$ belong to the norm closure of the linear set generated by $\underset{L_2\in L}\bigcup \lambda_{L_2}A_2$, and such that, for all $P_1\in P$, all operators $X\rho_{P_1}$ and $X^*\rho_{P_1}$ in $B(K, H\underset{\gb}{_L\otimes_P}K)$ belong to the norm closure of the linear set generated by $\underset{P_2\in P}\bigcup \rho_{P_2}A_1$, is a $\bf C^*$-algebra, which will be denoted $A_1\underset{\gb}{_L*_P}A_2$, and called the \emph{$\bf C^*$-fiber product} of the $\bf C^*$-$\gb$-algebra $A_1$ and the $\bf C^*$-$\gb^\dag$-algebra $A_2$. With the notations of \ref{def3}, we have $\varsigma(A_1\underset{\gb}{_L*_P}A_2)=A_2\underset{\gb^\dag}{_P*_L}A_1$. 

If $X$ belongs to $A_1\underset{\gb}{_L*_P}A_2$, then, for any $L_1$, $L_2$ in $L$, $\lambda_{L_2}^*X\lambda_{L_1}$ belongs to $A_2$, and, for any $P_1$, $P_2$ in $P$, $\rho_{P_2}^*X\rho_{P_1}$ belongs to $A_1$ (\cite{Ti2}, 3.16(i)).

If $(\tilde{H}, \tilde{L})$ is another $\bf C^*$-$\gb$-module, and $\tilde{A_1}$ a $\bf C^*$-$\gb$-algebra on $(\tilde{H}, \tilde{L})$, and $\Phi$ is a $\gb$-morphism from $A_1$ to $M(\tilde{A_1})$, and if $(\tilde{K}, \tilde{P})$ is another $\bf C^*$-$\gb^\dag$-module, and $\tilde{A_2}$ a $\bf C^*$-$\gb^\dag$-algebra on $(\tilde{K}, \tilde{P})$, and $\Psi$ a $\gb^\dag$-morphism from $A_2$ to $M(\tilde{A_2})$, then, there exists a $*$-homomorphism $\Psi\underset{\gb}{_P*_L}\Phi$ from $A_1\underset{\gb}{_L*_P}A_2$ to $M(\tilde{A_1})\underset{\gb}{_{\tilde{L}}*_{\tilde{P}}}M(\tilde{A_2})\subset M(\tilde{A_1}\underset{\gb}{_{\tilde{L}}*_{\tilde{P}}}\tilde{A_2})$. (\cite{Ti2}, 3.20). This homomorphism may be degenerate. 

Let us remark that, if $L_1\in L$, and $X\in B(H\underset{\gb}{_{L}\otimes_P}K)$, then $X\lambda_{L_1}$ (and $X^*\lambda_{L_1}$) belongs to $B(H, H\underset{\gb}{_{L}\otimes_P}K)$, and, therefore, to the norm closure of the linear set generated by $\underset{L_2\in L}\bigcup \lambda_{L_2}B(H)$; we get also that $X\rho_{P_1}$ and $X^*\rho_{P_1}$ belong to the norm closure of of the linear set generated by $\underset{R_2\in R}\bigcup \lambda_{R_2}B(K)$; therefore $B(H\underset{\gb}{_{L}\otimes_P}K)=B(H)\underset{\gb}{_L*_P}B(K)$.

Let us suppose now that all operators $X\lambda_{L_1}$ and $X^*\lambda_{L_1}$ in $B(H, H\underset{\gb}{_L\otimes_P}K)$ belong to the norm closure of the linear set generated by $\underset{L_2\in L}\bigcup \lambda_{L_2}A_2$; then, as $X\rho_{P1}$ and $X^*\rho_{P_1}$ belongs to the linear set generated by $\underset{R_2\in R}\bigcup \lambda_{R_2}B(H)$, we get that $X\in B(H)\underset{\gb}{_L*_P}A_2$,

%%%%ex
\subsubsection{\bf{Example}}
\label{C}
$\gt=(\mathbb C, \mathbb C, \mathbb C)$ is a $\bf C^*$-base. Then, if $H$, $K$ are two Hilbert spaces, it is clear that the $\bf C^*$-relative tensor product $H\underset{\gt}{\otimes}K$ is just the Hilbert tensor product $H\otimes K$ (\cite{Ti2}, 2.13); moreover, if $A$ is a sub-$\bf C^*$-algebra of $B(H)$, and $K$ is a sub-$\bf C^*$-algebra of $B(K)$, then the $\bf C^*$-fiber product $A\underset{\gt}{*}B$ contains the $\bf C^*$-algebra $\widetilde{M}(A\otimes B)=\{X\in M(A\otimes B)$, such that $X(1\otimes b)\in A\otimes B$, for all $b\in B$, and $X(a\otimes 1)\in A\otimes B$, for all $a\in A$\}. (\cite{Ti2}, 3.20).

%%%%%weights

\subsection{Weights on $\bf C^*$-algebras (\cite{C1}, \cite{C2}, \cite{K1})}
\label{weights}

%%%%%aut
\subsubsection{\bf{Notations}} 
\label{aut}
Let $M$ be a von Neuman algebra, and $\alpha$ an action from a locally compact group $G$ on $M$, i.e. a homomorphism from $G$ into $Aut M$, such that, for all $x\in M$, the function $g\mapsto \alpha_g(x)$ is $\sigma$-weakly continuous. Let us denote by ${\bf C}^*(\alpha)$ the set of elements $x$ of $M$, such that this function $t\mapsto \alpha_g(x)$ is norm continuous. It is (\cite{Pe}, 7.5.1) a sub-${\bf C}^*$-algebra of $M$, invariant under the $\alpha_g$, generated by the elements ($x\in N$,$f\in L^1(G)$):
\[\alpha_f(x)=\int_{\mathbb{R}}f(s)\alpha_s(x)ds\]
More precisely, we get that, for any $x$ in $M$, $\alpha_f(x)$ is $\sigma$-weakly converging to $x$ when $f$ goes in an approximate unit of $L^1(G)$, which proves that ${\bf C}^*(\alpha)$ is $\sigma$-weakly dense in $M$, and that $x\in M$ belongs to ${\bf C}^*(\alpha)$ if and only if this file is norm converging. 
\newline
If $\alpha_t$ and $\gamma_s$ are two one-parameter automorphism groups of $M$, such that, for all $s$, $t$ in $\mathbb{R}$, we have $\alpha_t\circ\gamma_s=\gamma_s\circ\alpha_t$, by considering the action of $\mathbb{R}^2$ given by $(s,t)\mapsto \gamma_s\circ\alpha_t$, we obtain a dense sub-${\bf C}^*$-algebra of $M$, on which both $\alpha$ and $\gamma$ are norm continuous, we shall denote ${\bf C}^*(\alpha, \gamma)$. 
\newline
If $\varphi$ is a normal semi-finite faithful weight on $M$, we shall write $\bf C^*(\varphi)$ for the norm closure of $\gM_\varphi\cap\bf C^*(\sigma^\varphi)$.

%%%%%defweight
\subsubsection{\bf{Definition}}
\label{defweight}
Let $A$ be a $\bf C^*$-algebra; a \emph{weight }$\nu$ on $A$ a is function $A^+\mapsto [0, +\infty]$ such that,
for any $x$, $y$ in $A^+$, and $\lambda\in\bf R^+$, we have $\nu (x+y)=\nu(x)+\nu (y)$ and $\nu(\lambda x)=\lambda\nu(x)$. 

We note $\gM_\nu^+=\{x\in A^+, \nu(x) <\infty\}$, $\gN_\nu=\{x\in A, \nu(x^*x)<\infty\}$ and $\gM_\nu$ for the linear space generated by $\gM_\nu ^+$ (or by all products $x^*y$, where $x$, $y$ are in $\gN_\nu$). 

As $\nu$ is an increasing function, we get that $\gM_\nu^+$ is an hereditary cone, $\gN_\nu$ a left ideal (in $M(A)$), and that $\gM_\nu$ is a sub-$*$-algebra of $A$, and that $\gM_\nu^+=\gM_\nu\cap A^+$ (which justify the notation). Moreover, $\nu$ can be extended to a linear map on $\gM_\nu$, we shall denote again $\nu$. We denote $N_\nu$ the set of all $x\in A$, such that $\nu(x^*x)=0$; it is clear that $N_\nu$ is a left-ideal of $A$.

We shall say that $\nu$ is densely defined if $\gM_\nu^+$ is dense in $A^+$ (or if $\gM_\nu$ (resp. $\gN_\nu$) is dense in $A$), and that $\nu$ is faithful if, for $x\in A^+$, $\nu(x)=0$ implies that $x=0$ (and then $N_\nu=\{0\}$).

%%%%GNS
\subsubsection{\bf{Definition}}
\label{GNS}
Be given a $\bf C^*$-algebra $A$ and a weight $\nu$ on $A$, a GNS construction for $\nu$ is a triple $(H_\nu, \pi_\nu, \Lambda_\nu)$ such that :

(i) $H_\nu$ is a Hilbert space, 

(ii) $\Lambda_\nu$ is a linear map from $\gN_\nu$ in $H_\nu$, such that $\Lambda_\nu(\gN_\nu)$ is dense in $H_\nu$, and, for any $x$, $y$ in $\gN_\nu$, we have $(\Lambda_\nu(x)|\Lambda_\nu(y))=\nu(y^*x)$,

(iii) $\pi_\nu$ is a representation of $A$ on $H_\nu$, such that, for any $a\in A$, $\Lambda_\nu(ax)=\pi_\nu(a)\Lambda_\nu(x)$.

For a construction, we refer to ([C1], 2).

%%%lsc
\subsubsection{\bf{Definition}}
\label{lsc}
We shall say that a weight $\nu$ on a $\bf C^*$-algebra $A$ is \emph{lower semi-continuous} (l.s.c.) if, for all $\lambda\in\bf R^+$, the set $\{x\in A^+, \nu(x)\leq\lambda\}$ is closed. 

Then, it is proved ([C1], 1.7) that, for any $x\in A^+$, $\nu(x)=sup\{\omega(x), \omega\in A^*_+, \omega\leq \nu\}$. and ([K1], 2.3) that $\Lambda_\nu$ is closed. Moreover, $\nu$ has then a natural extension to $M(A)$. 

%%%%KMS
\subsubsection{\bf{Definition}}
\label{KMS}
Be given a $\bf C^*$-algebra $A$, a densely defined lower semi-continuous faithful weight $\nu$ on $A$, and a norm continuous one parameter group of automorphism $\sigma$ on $A$; we shall say par $\nu$ is a \emph{KMS weight} on $A$ (with respect to $\sigma$) if :

(i) for all $t\in\bf R$, $\nu\circ\sigma_t=\nu$;

(ii) for any $x$, $y$ in $\gN_\nu\cap\gN_\nu^*$, there exists a bounded function $f$ in the set $\{0\leq Im z\leq 1\}\subset \bf C$, analytic in $\{0 < Im z< 1\}$, such that, for all $t\in \bf R$, $f(t)=\nu(\sigma_t(x)y)$ and $f(t+i)=\nu(y\sigma_t(x))$ (the so-called KMS conditions). 

Then ([EVal],2.2.3; the proof is due to F. Combes), $\nu$ extends to a normal semi-finite faithful weight on the von Neumann algebra $\pi_\nu (A)''$, we shall denote by $\underline{\nu}$. Then $\sigma$ is unique and is the restriction to $A$ of the modular group $\sigma^{\underline{\nu}}$. 

With the notations of \ref{aut}, we get that, if $\varphi$ is a normal semi-finite faithful weight on a von Neumann algebra $M$, then $\varphi_{|\bf C^*(\varphi)}$ is a KMS weight on $\bf C^*(\varphi)$. 

%%%%analytic
\subsubsection{\bf{Definition}}
\label{analytic}
Be given a $\bf C^*$-algebra $A$, and a densely defined, lower semi-continuous, faithful weight $\nu$, KMS with respect to a one parameter group $\sigma$ of automorphisms, we shall say that $x$ is \emph{analytic} if the function $t\mapsto \sigma_t(x)$ extends to an analytic function. 

If $x\in\gN_\nu$, then $x_n=\frac{n}{\sqrt{\pi}}\int e^{-n^2t^2}\sigma_t(x)dt$ is analytic, belongs to $\gN_\nu$ and $x_n$ is norm converging to $x$ and $\Lambda_\nu(x_n)$ is norm converging to $\Lambda_\nu(x)$. (\cite{K1} 4.3)

%%%C*valuedweight
\subsection{$\bf C^*$-valued weights (\cite{K2})}
\label{C*valuedweight}

%%%%defT
\subsubsection{\bf{Definition}}
\label{defT}
Be given two $\bf C^*$-algebras $A$ and $B$, with $B\subset M(A)$, and a hereditary cone $P$ in $A^+$; let us write $\gN=\{a\in A, a^*a\in P\}$ and $\gM=span P=\gN^*\gN$; then a $\bf C^*$-\emph{valued weight} from $A$ into $M(B)$ is a linear map $T$ from $\gM$ into $M(B)$, such that, for any $b\in B$ and $x\in P$, $b^*xb$ belongs to $P$ and $T(b^*xb)=b^*T(x)b$. 

Then, $\gM$ will be denoted $\gM_T$, $\gN$ will be denoted $\gN_T$, and $P$ will be denoted $\gM_T^+$. 
We shall say that $T$ is densely defined if one of these sets is dense in $A$ (or $A^+$). $T$ is faithful if ($x\in\gM_T^+$), $Tx=0$ implies that $x=0$. If $T$ is faithful and densely defined, then it is easy to get that $T(\gM_T)$ is a dense ideal in $M(B)$. 

A definition of lower semi-continuity is given in (\cite{K2}, 3). More precisely, Kustermans constructs an uprising set $\mathcal G_T$ of bounded completely positive maps from $A$ to $M(B)$, and $T$ is then said lower semi-continuous if $\gM_T^+$ is the set of $x\in A^+$ such that $(\rho(x))_{\rho\in\mathcal G_T}$ is strictly convergent in $M(B)$, and if this limit is then equal to $T(x)$. So, when $\nu$ is a lower semi-continuous weight on $B$, which then extends to $M(B)$ (\ref{lsc}), then $\nu\circ T$ is a lower semi-continuous weight on $A$. More precisely, if $\nu$ and $T$ are densely defined, lower semi-continuous and faithful, so is $\nu\circ T$. 

%%%KMSC^*
\subsubsection{\bf{Definition}}
\label{KMSC*}
Be given two $\bf C^*$-algebras $A$ and $B$, with $B\subset M(A)$, a faithful densely defined lower semi-continuous $\bf C^*$-weight $T$ from $A$ to $M(B)$, and a faithful densely defined lower semi-continuous weight $\nu$ on $B$, which then extends to $M(B)$. We shall say that the pair $(\nu, T)$ is \emph{KMS} if :

(i) there exists a one parameter automorphism group $\sigma^\nu$ on $B$, such that $\nu$ is KMS with respect to $\sigma^\nu$;

 (ii) there exists a one parameter automorphism group $\sigma^{\nu\circ T}$ on $A$, such that $\nu\circ T$ is KMS with respect to $\sigma^{\nu\circ T}$;
 
 (iii) for all $t\in\bf R$, the restriction of $\sigma_t^{\nu\circ T}$ to $B$ is equal to $\sigma_t^\nu$. 
 
 In that situation, considering the GNS constructions of $\nu$ and $\nu\circ T$, we define, for any $a\in\gN_T$, the linear map $\Lambda_T(a)\in B(H_\nu, H_{\nu\circ T})$ by ($b\in\gN_\nu$) :
 \[\Lambda_T(a)\Lambda_\nu(b)=\Lambda_{\nu\circ T}(ab)\]
 
 Let us remark that $\|\Lambda_T(a)\|^2=\|T(a^*a)\|$.

 %%%%thKMS
 \subsubsection{\bf{Theorem}}
 \label{thKMS}
 {\it Be given two $\bf C^*$-algebras $A$ and $B$, with $B\subset M(A)$, a densely defined, faithful, lower semi-continuous weight $\nu$ on $B$, a densely defined, faithful, lower semi-continuous $\bf C^*$-valued weight $T$ from $A$ into $M(B)$, such that $(\nu, T)$ is KMS, in the sense of \ref{KMSC*}. Then :
 
 (i) there exists a unique injective $*$-homomorphism $\Phi$ from $\pi_\nu(B)''$ into $B(H_{\nu\circ T})$ such that $\Phi\circ\pi_\nu=\pi_{\nu\circ T}$, which allows to consider $\pi_\nu (B)''$ as a sub-algebra of  $\pi_{\nu\circ T}(A)''$. 
 
 (ii) moreover, there exists a normal faithful semi-finite operator-valued weight $\underline{T}$ from $\pi_{\nu\circ T}(A)''$ onto $\pi_\nu(B)''$ such that, for any $a$ in $\gM_T$, we have $\underline{T}(\pi_{\nu\circ T}(a))=\pi_\nu(T(a))$. }
 
 \begin{proof}
 For any $b\in B$ and $t\in \bf R$, we have :
\[\sigma_t^{\nu\circ T}\Phi\pi_\nu(b)=\sigma_t^{\nu\circ T}\pi_{\nu\circ T}(b)=\pi_{\nu\circ T}(\sigma_t^\nu(b))=\Phi\pi_\nu\sigma_t^\nu(b)\]
and, by continuity, we have $\sigma_t^{\nu\circ T}\Phi=\Phi\sigma_t^\nu$. By identifying $\pi_\nu(B)''$ with $\Phi(\pi_\nu(B)'')$, we may consider $\pi_\nu(B)''$ as a sub-algebra of $\pi_{\nu\circ T}(A)''$; then the restriction of $\sigma_t^{\nu\circ T}$ to this sub-algebra is equal to $\sigma_t^\nu$. This gives the existence of a normal semi-finite faithful operator valued weight $\underline{T}$ from $\pi_{\nu\circ T}(A)''$ onto $\pi_\nu (B)''$. The fact that the restriction of $\underline{T}$ to $\pi_{\nu\circ T}(A)$ is equal to $T$ is straightforward. \end{proof}
Note that the existence of $\Phi$ can also be deduced from (\cite{C3}, 1.7 and 1.8). Thanks to that result, we get, by restriction of $\underline{T}$, that, for all $t\in\bf R$, $\sigma_t^{\nu\circ T}=T\circ\sigma_t^{\nu\circ T}$, from which we get that $\sigma_t^{\nu\circ T}(\gN_T)\subset \gN_T$; if $x\in\gM_T$, then $x_n=\frac{n}{\sqrt{\pi}}\int e^{-n^2t^2}\sigma_t^{\nu\circ T}(x)dt$ is analytic with respect to $\sigma^{\nu\circ T}$ and belongs to $\gM_T$; moreover, if $x_1$, $x_2$ are in $\gN_T$ and analytic with respect to $\sigma^{\nu\circ T}$, then $T(x_2^*x_1)$ is analytic with respect to $\sigma^\nu$. 

Moreover, we have clearly $\underline{\nu\circ T}=\underline{\nu}\circ\underline{T}$. 

%%%rest
\subsubsection{\bf{Theorem}}
\label{rest}
{\it Let $M$ be a von Neumann algebra, $N$ a sub-von Neumann algebra of $M$, and $T$ a normal faithful semi-finite operator-valued weight from $M$ to $N$, $\nu$ a normal semi-finite faithful weight on $N$, and $\varphi=\nu\circ T$; let us define $\bf C^*$$(T, \nu)$ as the norm closure of $\gM_T\cap\bf C^*(\varphi)$; then :

(i) $\bf C^*(\sigma^\nu)$ is included into $\bf C^*(\sigma^\varphi)$;  if $x\in \gM_T\cap \bf C^*(\sigma^\varphi)$, then $T(x)$ belongs to $\bf C^*(\sigma^\nu)$; 

(ii) $\bf C^*$$(T, \nu)$ is equal to $\bf C^*(\varphi)$; 

(iii) the restriction of $T$ to $\bf C^*$$(T, \nu)$ is a lower semi-continuous densely defined $\bf C^*$-valued weight from $\bf C^*($$T, \nu)$ to $\bf C^*(\nu)$ and $(\nu_{|\bf C^*(\nu)}, T_{|\bf C^*(T,\nu)})$ is KMS. }

\begin{proof}
As ${\sigma_t^\varphi}_{|N}=\sigma_t^\nu$, it is clear that  $\bf C^*(\sigma^\nu)$ is included into $\bf C^*(\sigma^\varphi)$; let now $x\in \gM_T^+\cap \bf C^*(\sigma^\varphi)$. Using \ref{aut}, we get that $x$ is the norm limit of $\int_{\bf R}f_i(t)\sigma_t^\varphi(x)dt$, when $f_i$ is a positive continuous approximate unit of $L^1(\bf R)$.

Let $\omega\in N^{*+}$, such that $\omega_{|\bf C^*(\sigma^\nu)}=0$; as $\omega\circ T$ is norm lower semi-continuous, we get that $\omega\circ T$ is equal to $\sup\{\omega'\in M^{+*}, \omega'\leq\omega\circ T\}$; then, using the norm continuity of $\omega'$, we have :
\begin{multline*}
\omega'(x)=lim_i\int_{\bf R}f_i(t)\omega'(\sigma_t^\varphi(x))dt\leq \omega\circ T (\int_{\bf R}f_i(t)\sigma_t^\varphi(x)dt)=\\=lim_i\omega\circ\int_{\bf R}f_i(t)\sigma_t^\nu (Tx)dt=0
\end{multline*}
because, by \ref{aut} again, all elements of the form $\int_{\bf R}f(t)\sigma_t^\nu(y)dt$, for any $y\in N$ and $f\in L^1(\bf R)$, belong to $\bf C^*(\sigma^\nu)$; therefore, we get that $\omega\circ T(x)=0$, and, then, that $T(x)$ belongs to $\bf C^*(\sigma^\nu)$. 

(ii) If $x\in\gM_T\cap\bf C^*(\sigma^\varphi)$ and $y\in\gN_\nu\cap\bf C^*(\sigma^\nu)$, using (i), we get that $xy$ belongs to $\gM_T\cap\gM_\varphi\cap\bf C^*(\sigma^\varphi)$ and that such elements are dense in both $\bf C^*($$T, \nu)$ and $\bf C^*$$(\varphi)$, which is (ii). 

(iii) as $\gM_T$ is invariant under $\sigma_t^\varphi$, we get that  $\varphi_{|\bf C^*(\varphi)}$ is KMS; as $\nu_{|\bf C^*(\nu)}$ is also KMS, we get easily the result. \end{proof}

%%%basic
\subsection{Basic construction for inclusion of $\bf C^*$-algebras with KMS $\bf C^*$-valued weights.}
\label{basicC*}
Let $B$, $A$ be two $\bf C^*$-algebras, with $B\subset M(A)$. Let $\nu$ be a KMS weight on $B$, and $T$ a densely defined lower semi-continuous $\bf C^*$-valued weight from $A$ to $M(B)$, such that the pair $(\nu, T)$ is KMS. We define $\varphi=\nu\circ T$ which is a KMS weight on $A$. We shall denote $M_0=\pi_\nu(B)''$, $M_1=\pi_\varphi(A)''$, $\underline{\nu}$ the canonical extension of $\nu$ to $M_0$ (which, by \ref{KMS}, is a normal semi-finite faithful weight on $M_0$), $\underline{\varphi}$ the canonical extension of $\varphi$ to $M_1$ (which, by \ref{KMS} again, is a normal semi-finite faithful weight on $M_1$). Using again \ref{thKMS}, we get that $M_0$ can be considered as a sub-von Neumann algebra of $M_1$, and that there exists $\underline{T}$ a canonical extension of $T$ to $M_1$ (which, by \ref{thKMS}, is a normal semi-finite faithful operator-valued weight from $M_1$ to $M_0$). Let $M_2$ be the basic construction made from the inclusion $M_0\subset M_1$, and let $\underline{T_2}$ be the normal faithful semi-finite operator-valued weight constructed in \ref{basisT}(i) from $M_2$ to $M_1$. 

%%%%lemmaappunit
\subsubsection{\bf{Lemma}}
\label{lemmaappunit}
{\it Let $a$, $b$ in $A\cap\gN_T\cap\gN_\varphi$; then :

(i)  $T(b^*a)$ belongs to $M(B)$. 

(ii) The norm closure of the linear space generated by all elements of the form $\Lambda_T(a)\Lambda_T(b)^*$ is a $\bf C^*$-algebra $A_1$ and $A\subset M(A_1)$. 

(iii) There exsits elements $a_i$, $b_i$ in $A\cap\gN_T\cap\gN_\varphi\cap\gN_T^*\cap\gN_\varphi^*$ such that, for any $y\in A_1$, the sum $y\Sigma_i\Lambda_T(a_i)\Lambda_T(b_i)^*$ is norm converging to $y$. }

\begin{proof} The proof of (i) is trivial. For any $a_1$, $a_2$, $b_1$, $b_2$ in $\gN_T\cap\gN_\varphi\cap\gN_T^*\cap\gN_\varphi^*$, we have :
\[\Lambda_{T}(a_1)\Lambda_T(b_1)^*\Lambda_T(a_2)\Lambda_T(b_2)^*=\Lambda_T(a_1(T(b_1^*a_2)))\Lambda_T(b_2)^*\]
By (i) $T(b_1^*a_2)$ belongs to $M(B)\subset M(A)$, so $a_1(T(b_1^*a_2))$ belongs to $A$. Moreover, it is easy to get that it belongs also to $\gN_T\cap\gN_\varphi$, which gives that the linear space generated by all elements of the form $\Lambda_T(a)\Lambda_T(b)^*$ is an algebra, and the first part of (ii). The fact that $A\subset M(A_1)$ is trivial; which finishes the proof of (ii). Now, let us choose an approximate unit of $A_1$ of the form $\Sigma_i\Lambda_T(a_i)\Lambda_T(b_i)^*$, and we get (iii). \end{proof}
 
 %%%%xanal
 \subsubsection{\bf{Lemma}}
 \label{xanal}
 {\it Let $x\in\gN_T\cap B'$, $x$ analytic with respect to $\sigma_t^\varphi$; then $\sigma_{-i/2}^\varphi (x^*)$ belongs to $\gN_T$.}
 
 \begin{proof}
 Let $b$ in $\gN_\nu$, analytic with respect to $\sigma_t^\nu$, such that $\sigma_{-i/2}(b^*)$ belongs to $\gN_\nu$; we have :
 \begin{multline*}
  (T(x^*x)J_\nu\Lambda_\nu(b)|J_\nu\Lambda_\nu(b))=\nu(\sigma_{-i/2}^\nu(b^*)^*T(x^*x)\sigma_{-i/2}^\nu(b^*))=\varphi(\sigma_{-i/2}^\nu(b^*)^*x^*x\sigma_{-i/2}^\nu(b^*))=\\\|\Lambda_\varphi(\sigma_{-i/2}(b^*)x)\|^2=\|J_\varphi \Lambda_\varphi(\sigma_{-i/2}(b^*)x)\|^2
=\varphi(b^*\sigma_{-i/2}^\varphi(x^*)^*\sigma_{-i/2}^\varphi(x^*)b)
 \end{multline*}
 We get then that $\sigma_{-i/2}^\varphi (x^*)$ belongs to $\gN_T$ .
  \end{proof}

%%%%fiber
\subsection{\bf{Fiber product of $\bf C^*$-algebras}}
\label{fiber}
Here, using weights and $\bf C^*$-valued weights on $\bf C^*$-algebras, we give another definition of the fiber product of two $\bf C^*$-algebras; this definition is completely inspired by Timmermann's work, and, in that special situation, is equivalent to (\cite{Ti2}, 3.3) . 

%%%%%deffiber
\subsubsection{\bf{Definition}}
\label{deffiber}

Be given three $\bf C^*$-algebras $A_1$, $A_2$ and $B$, $\alpha$ an injective strict $*$-homomorphism from $B$ into $M(A_1)$, $\beta$ an injective strict $*$-antihomomorphism from $B$ into $M(A_2)$, $T_1$ a densely defined faithful lower semi-continuous $\bf C^*$-valued weight from $A_1$ to $M(B)$, $T_2$ a densely defined faithful lower semi-continuous $\bf C^*$-valued weight from $A_2$ to $M(B)$, $\nu$ a KMS weight on $B$ such that both pairs $(\nu, T_1)$ and $(\nu, T_2)$ are KMS. Let $N=\pi_\nu(B)''$, $\varphi_1=\nu\circ\alpha^{-1}\circ T_1$, $\varphi_2=\nu^o\circ\beta^{-1}\circ T_1$ $M_1=\pi_{\varphi_1}(A_1)''$, $M_2=\pi_{\varphi_2}(A_2)''$; let $\underline{\alpha}$ (resp. $\underline{\beta}$) be the imbedding of $N$ into $M_1$ (resp. $M_2$) given by \ref{thKMS}. Then $\underline{\alpha}$ (resp. $\underline{\beta}$) is an injective $*$-homomorphism (resp. anti-$*$-homomorphism) of $N$ into $M_1$ (resp. $M_2$) and let $M_2\underset{N}{_{\underline{\beta}} *_{\underline{\alpha}}}M_1$ their fiber product in the sense of \ref{fiber}. Let $\underline{\varphi_1}$ (resp. $\underline{\varphi_2}$) be the extension of $\varphi_1$ (resp. $\varphi_2$) to $M_1^+$ (resp. $M_2^+$). 
Then, we shall denote $A_2\underset{B}{_{\beta}*_{\alpha}}A_1$ the set of all elements $X\in M_2\underset{N}{_\beta*_\alpha}M_1$ such that :

- for all $y\in A_1\cap\gN_{T_1}\cap\gN_{\varphi_1}$, $X\rho^{\beta, \alpha}_{J_{\varphi_1}\Lambda_{\varphi_1}(y)}$ and $X^*\rho^{\beta, \alpha}_{J_{\varphi_1}\Lambda_{\varphi_1}(y)}$ belong to the norm closure of the linear set generated by :
\[\cup_{z\in A_1\cap\gN_{T_1}\cap\gN_{\varphi_1}}\rho^{\beta, \alpha}_{J_{\varphi_1}\Lambda_{\varphi_1}(z)}A_2\]

- for all $y'\in A_2\cap\gN_{T_2}\cap\gN_{\varphi_2}$, $X\lambda^{\beta, \alpha}_{J_{\varphi_2}\Lambda_{\varphi_2}(y')}$ and $X^*\lambda^{\beta, \alpha}_{J_{\varphi_2}\Lambda_{\varphi_2}(y')}$ belong to the norm closure of of the linear set generated by :
\[\cup_{z'\in A_2\cap\gN_{T_2}\cap\gN_{\varphi_2}}\lambda^{\beta, \alpha}_{J_{\varphi_2}\Lambda_{\varphi_2}(z')}A_1\]

It is clear that $A_2\underset{B}{_{\beta}*_{\alpha}}A_1$ is a sub-$\bf C^*$-algebra of $M_2\underset{N}{_{\underline{\beta}}*_{\underline{\alpha}}}M_1$, but it may be degenerate ([Ti2], 3.20). 

%%%%thfiber
\subsubsection{\bf{Theorem}}
\label{thfiber}
{\it Let's use the notations of \ref{deffiber}; let $X\in A_2\underset{B}{_{\beta}*_{\alpha}}A_1$, $x_1$, $x'_1$ in $\gN_{T_1}\cap \gN_{\varphi_1}$, $x_2$, $x'_2$ in $\gN_{T_2}\cap\gN_{\varphi_2}$; then :

(i) $(\omega_{\Lambda_{\varphi_2}(x'_2),\Lambda_{\varphi_2}(x_2)}\underset{\underline{\nu}}{_{\underline{\beta}}*_{\underline{\alpha}}}id)(X)$ belongs to $A_1$; 

(ii) $(id\underset{\underline{\nu}}{_{\underline{\beta}}*_{\underline{\alpha}}}\omega_{J_{\underline{\varphi_1}}\Lambda_{\varphi_1(x'_1)}, J_{\underline{\varphi_1}}\Lambda_{\varphi_1}(x_1)})(X)$ belongs to $A_2$; }
\begin{proof}
By \ref{deffiber}, we get that the operator $X\lambda_{\Lambda_{\varphi_2}(x_1)}$ belongs to the norm closure of $\cup_{x''} \lambda^{\beta, \alpha}_{\Lambda_{\varphi_2}(x'')}A_1$, for all $x''$ in $\gN_{\varphi_2}$, analytic with resect to $\varphi_2$. Then, we get that the operator $(\omega_{\Lambda_{\varphi_2}(x'_1),\Lambda_{\varphi_2}(x_1)}\underset{\nu}{_\beta*_\alpha}id)(X)$ belongs to the norm closure of all operators of the form $1\underset{\nu}{_{\widetilde{\beta}}\otimes_{\widetilde{\alpha}}}\alpha\circ\sigma_{-i/2}^\nu\circ\beta^{-1}\circ T_2(x^{'*}_1x'')A_1$, which is included in $\alpha(M(B))A_1\subset A_1$. From which we get (i) by density. Result (ii) is proved the same way.  \end{proof}

%%%%LCQG
 \subsection{Locally compact quantum groupoids}
  \label{defLCQG}
%%%threltens
\subsubsection{\bf{Theorem}}
\label{threltens}
{\it Be given three $\bf C^*$-algebras $A_1$, $A_2$ and $B$, $\alpha$ an injective strict $*$-homomorphism from $B$ into $M(A_1)$, $\beta$ an injective strict $*$-antihomomorphism from $B$ into $M(A_2)$, $T_1$ a densely defined faithful lower semi-continuous $\bf C^*$-valued weight from $A_1$ to $M(B)$, $T_2$ a densely defined faithful lower semi-continuous $\bf C^*$-valued weight from $A_2$ to $M(B)$, $\nu$ a KMS weight on $B$ such that both pairs $(\nu, T_1)$ and $(\nu, T_2)$ are KMS. Then, $(H_{\nu\circ\alpha^{-1}\circ T}, L_{T_1})$ is a $\bf C^*$-$\gb_\nu$-module, $(H_{\nu\circ\beta^{-1}\circ T_2}, L_{T_2})$ is a $\bf C^*$-$\gb_\nu^o$-module, and 
the application from $H_{\nu^o\circ \beta^{-1}\circ T_2}\underset{\nu}{_\beta\otimes_\alpha}H_{\nu\circ\alpha^{-1}\circ T_1}$ to $H_{\nu^o\circ \beta^{-1}\circ T_2}\underset{\gb_\nu}{_{L_{T_2}}\otimes_{L_{T_1}}}H_{\nu\circ\alpha^{-1}\circ T_1}$ which sends (for any $x\in \gN_{T_1}\cap\gN_{\nu\circ\alpha^{-1}\circ T_1}$, $y\in \gN_{T_2}\cap\gN_{\nu^o\circ \beta^{-1}\circ T_2}$, $z\in\gN_\nu$) :
\begin{multline*}
\Lambda_{\nu^o\circ \beta^{-1}\circ T_2}(\beta(z^*)y)\underset{\nu}{_\beta\otimes_\alpha}J_{\nu\circ\alpha^{-1}\circ T_1}\Lambda_{\nu\circ\alpha^{-1}\circ T_1}(x)=\\
=\Lambda_{\nu^o\circ \beta^{-1}\circ T_2}(y)\underset{\nu}{_\beta\otimes_\alpha}J_{\nu\circ\alpha^{-1}\circ T_1}\Lambda_{\nu\circ\alpha^{-1}\circ T_1}(x\alpha(z))
\end{multline*}
on $\Lambda_{T_2}(y)\vg J_\nu\Lambda_\nu(z)\vl\Lambda_{T_1}(x)$ is an isomorphism of Hilbert spaces. }
\begin{proof}
Note that the equality given up here is an easy corollary of \ref{spatial}. The second term can be written also as $\Lambda_{\nu^o\circ \beta^{-1}\circ T_2}(y)\underset{\nu}{_\beta\otimes_\alpha}J_{\nu\circ\alpha^{-1}\circ T_1}\Lambda_{ T_1}(x)\Lambda_\nu(z)$

Let $x_1$, $x_2$ in $\gN_{T_1}\cap\gN_{\nu\circ\alpha^{-1}\circ T_1}$, $y_1$, $y_2$ in $\gN_{T_2}\cap\gN_{\nu^o\circ \beta^{-1}\circ T_2}$, $z_1$, $z_2$ in $\gN_\nu$. Using \ref{lemT}, we get that the scalar product :
\[(\Lambda_{\nu^o\circ \beta^{-1}\circ T_2}(y_1)\underset{\nu}{_\beta\otimes_\alpha}J_{\nu\circ\alpha^{-1}\circ T_1}\Lambda_{T_1}(x_1)\Lambda_\nu(z_1)|\Lambda_{\nu^o\circ \beta^{-1}\circ T_2}(y_2)\underset{\nu}{_\beta\otimes_\alpha}J_{\nu\circ\alpha^{-1}\circ T_1}\Lambda_{T_1}(x_2)\Lambda_\nu(z_2))\]
is equal to :
\[(\alpha(<\Lambda_{\nu^o\circ \beta^{-1}\circ T_2}(y_1), \Lambda_{\nu^o\circ \beta^{-1}\circ T_2}(y_2)>_{\beta, \nu^o})J_{\nu\circ\alpha^{-1}\circ T_1}\Lambda_{T_1}(x_1)\Lambda_\nu(z_1)|J_{\nu\circ\alpha^{-1}\circ T_1}\Lambda_{T_1}(x_2)\Lambda_\nu(z_2))\]
which can be written as :
\[(\alpha\circ\beta^{-1}T_2(y_2^*y_1)J_{\nu\circ\alpha^{-1}\circ T_1}\Lambda_{T_1}(x_1)\Lambda_\nu(z_1)|J_{\nu\circ\alpha^{-1}\circ T_1}\Lambda_{T_1}(x_2)\Lambda_\nu(z_2))\]
which is equal to :
\begin{multline*}
(J_{\nu\circ\alpha^{-1}\circ T_1}\alpha\circ\beta^{-1}T_2(y_1^*y_2)J_{\nu\circ\alpha^{-1}\circ T_1}\Lambda_{T_1}(x_2)\Lambda_\nu(z_2)|\Lambda_{T_1}(x_1)\Lambda_\nu(z_1))=\\
=(\Lambda_{T_1}(x_2)J_\nu\beta^{-1}T_2(y_1^*y_2)J_\nu\Lambda_\nu(z_2)|\Lambda_{T_1}(x_1)\Lambda_\nu(z_1)=\\
=(J_\nu\beta^{-1}T_2(y_1^*y_2)J_\nu\Lambda_\nu(z_2)|\alpha^{-1}T_1(x_2^*x_1)\Lambda_\nu(z_1))\\
=(\Lambda_\nu(z_2)|\alpha^{-1}T_1(x_2^*x_1)J_\nu\beta^{-1}T_2(y_2^*y_1)J_\nu\Lambda_\nu(z_1))
\end{multline*}
which, by definition (\ref{def3}), is equal to the scalar product :
\[(\Lambda_{T_2}(y_1)\vg J_\nu\Lambda_\nu(z_1)\vl\Lambda_{T_1}(x_1)|\Lambda_{T_2}(y_2)\vg J_\nu\Lambda_\nu(z_2)\vl\Lambda_{T_1}(x_2))\]
which proves the result. \end{proof}

%%%notationsfiber
\subsubsection{\bf{Notations}}
\label{notationsfiber}
Be given three $\bf C^*$-algebras $A_1$, $A_2$ and $B$, $\alpha$ an injective strict $*$-homomorphism from $B$ into $M(A_1)$, $\beta$ an injective strict $*$-antihomomorphism from $B$ into $M(A_2)$, $T_1$ a densely defined faithful lower semi-continuous $\bf C^*$-valued weight from $A_1$ to $M(B)$, $T_2$ a densely defined faithful lower semi-continuous $\bf C^*$-valued weight from $A_2$ to $M(B)$, $\nu$ a KMS weight on $B$ such that both pairs $(\nu, T_1)$ and $(\nu, T_2)$ are KMS. Let us write $A_2\underset{B}{_{\beta}*_{\alpha}}A_1$ be the $\bf C^*$-fiber product in the sense of \ref{def5} written $A_2\underset{\gb_\nu}{_{L_{T_2}}*_{L_{T_1}}}A_1$.

%%%%thfiber
\subsubsection{\bf{Theorem}}
\label{thfiber}
{\it Let us use the notations of \ref{notationsfiber}; let $N=\pi_\nu(B)''$, $M_1=\pi_{\nu\circ\alpha^{-1}\circ T_1}(A_1)''$, $M_2=\pi_{\nu^o\circ\beta^{-1}\circ T_2}(A_2)''$; let $\underline{\alpha}$ (resp. $\underline{\beta}$) be the imbedding of $N$ into $M_1$ (resp. $M_2$) given by \ref{thKMS}. Then $\underline{\alpha}$ (resp. $\underline{\beta}$) is an injective $*$-homomorphism (resp. anti-$*$-homomorphism) of $N$ into $M_1$ (resp. $M_2$) and let $M_2\underset{N}{_{\underline{\beta}} *_{\underline{\alpha}}}M_1$ their fiber product in the sense of \ref{fiber}. Using \ref{threltens}, we shall consider that $A_2\underset{B}{_{\beta}*_{\alpha}}A_1$ is a sub-$\bf C^*$-algebra of $B(H_{\nu^o\circ\beta^{-1}\circ T_2}\underset{\underline{\nu}}{_{\underline{\beta}}\otimes_{\underline{\alpha}}}H_{\nu\circ\alpha^{-1}\circ T_1})$. Let $X\in A_2\underset{B}{_{\beta}*_{\alpha}}A_1$, $x_1$, $x'_1$ in $\gN_{T_1}\cap\gN_{\nu\circ\alpha^{-1}\circ T_1}$, $x_2$, $x'_2$ in $\gN_{T_2}\cap\gN_{\nu^o\circ\beta^{-1}\circ T_2}$; then :

(i) $(\omega_{\Lambda_{\nu ^o\circ\beta ^{-1}\circ T_2}(x'_2),\Lambda_{\nu ^o\circ\beta^{-1}\circ T_2}(x_2)}\underset{\underline{\nu}}{_{\underline{\beta}}*_{\underline{\alpha}}}id)(X)$ belongs to $A_1$; 

(ii) $(id\underset{\underline{\nu}}{_{\underline{\beta}}*_{\underline{\alpha}}}\omega_{J_{\underline{\nu\circ\alpha^{-1}\circ T_1}}\Lambda_{\nu\circ\alpha^{-1}\circ T_1(x'_1)}, J_{\underline{\nu\circ\alpha^{-1}\circ T_1}}\Lambda_{\nu\circ\alpha^{-1}\circ T_1}(x_1)})(X)$ belongs to $A_2$; 

(iii) The $\bf C^*$-fiber product $A_2\underset{B}{_{\beta}*_{\alpha}}A_1$ is a sub-$\bf C^*$-algebra of $M_2\underset{N}{_{\widetilde{\beta}} *_{\widetilde{\alpha}}}M_1$.  

}
\begin{proof}
Let us suppose first that $x''$ is analytic. By \ref{def5}, we get that the operator $X\lambda_{\Lambda_{\nu^o\circ\beta^{-1}\circ T_2}(x_1)}$ belongs to the norm closure of $\cup_{x''} \lambda_{\Lambda_{\nu^o\circ\beta^{-1}\circ T_2}(x'')}A_1$, for all $x''$ in $\gN_{T_2}\cap\gN_{\nu^o\circ\beta^{-1}\circ T_2}$, analytic with resect to $\nu^o\circ\beta^{-1}\circ T_2$. Using now \ref{threltens}, we get that the operator $(\omega_{\Lambda_{\Lambda_{\nu^o\circ\beta^{-1}\circ T_2}(x'_1)},\Lambda_{\nu^o\circ\beta^{-1}\circ T_2}(x_1)}\underset{\nu}{_\beta*_\alpha}id)(X)$ belongs to the norm closure of all operators of the form $1\underset{\nu}{_{\widetilde{\beta}}\otimes_{\widetilde{\alpha}}}\alpha\circ\sigma_{-i/2}^\nu\circ\beta^{-1}\circ T_2(x^{'*}_1x'')A_1$, which is included in $\alpha(M(B))A_1\subset A_1$. From which we get (i) by density. Result (ii) is proved the same way. 
Then, it is straightforward to get that $X$ commutes with $1\underset{N}{_{\underline{\beta}}\otimes_{\underline{\alpha}}}y$, for all $y\in\pi_{\nu\circ\alpha^{-1}\circ T_1}(A_1)'$, and with $x\underset{N}{_{\underline{\beta}}\otimes_{\underline{\alpha}}} 1$, for all $x\in \pi_{\nu^o\circ \beta^{-1}\circ T_2}(A_2)'$; so $X$ commutes with $M'_2\underset{N}{_{\underline{\beta}}\otimes_{\underline{\alpha}}}M'_1$, which gives (iii). 
 \end{proof}
 
 %%%%%not
 \subsubsection{\bf{Notations}}
 \label{not}
 Using \ref{thfiber}, we define $(id\underset{B}{_{\beta}*_{\alpha}}T_1)$ as the restriction of $(id\underset{N}{_{\underline{\beta}}*_{\underline{\alpha}}}\widetilde{T_1})$ to $A_2\underset{B}{_{\beta}*_{\alpha}}A_1$; it is a faithful lower semi-continuous $\bf C^*$-valued weight from $A_2\underset{B}{_{\beta}*_{\alpha}}A_1$ to $A_2\cap \beta(B)'\underset{\gb_\nu}{_{L_{T_2}}\otimes_{L_{T_1}}}1$. The $\bf C^*$-valued weight $(T_2\underset{B}{_{\beta}*_{\alpha}}id)$ is defined the same way. 
 
 Let us write $\varphi=\nu\circ\alpha^{-1}\circ T_1$, which is a KMS weight on $A_1$, and $\underline{\varphi}$ its extension to $M_1$. Using again \ref{thfiber}, we define also $(id\underset{B}{_{\beta}*_{\alpha}}\varphi)$ as the restriction of $(id\underset{N}{_{\underline{\beta}}*_{\underline{\alpha}}}\underline{\varphi})$ to $A_2\underset{B}{_{\beta}*_{\alpha}}A_1$. As $(id\underset{\nu}{_{\underline{\beta}}*_{\underline{\alpha}}}\underline{\nu}\circ\underline{\alpha}^{-1})$ is just the inclusion of $M_2\cap\underline{\beta}(N)'=M_2\underset{N}{_{\underline{\beta}}*_{\underline{\alpha}}}\underline{\alpha}(N)$ into $M_2$ (\cite{L}, 2.5.1), we get that, for any positive $X\in A_2\underset{B}{_{\beta}*_{\alpha}}A_1$, we have :
 \[(id\underset{B}{_{\beta}*_{\alpha}}T_1)(X)=(id\underset{B}{_{\beta}*_{\alpha}}\varphi)(X)\underset{\underline{\nu}}{_{\underline{\beta}}\otimes_{\underline{\alpha}}}1\] 
 
 \subsubsection{\bf{Definition}}
\label{defLCQG}
 An octuple $\bfG=(B, A, \alpha, \beta, \nu, T, T', \Gamma)$ will be called a \emph{locally compact quantum groupoid} if :
 
 (i) $A$ and $B$ are $\bf C^*$-algebras, $\alpha$ a strict $*$-homomorphism from $B$ to $M(A)$, $\beta$ a strict anti-$*$-homomorphism from $B$ to $M(A)$; 
 
 (ii) $\nu$ is a densely defined faithful, lower semi-continuous KMS weight on $B$, in the sense of \ref{KMS}; 
 
 (iii) $T$ is a densely defined faithful lower semi-continuous $\bf C^*$-valued weight from $A$ to $\alpha(M(B))$, $T'$ is a densely defined faithful lower semi-continuous $\bf C^*$-valued weight from $A$ to $\beta(M(B))$, such that both pairs $(T, \nu\circ\alpha^{-1})$ and $(T', \nu^o\circ\beta^{-1})$ are KMS, in the sense of \ref{KMSC*};
 
 (iv) $\Gamma$ is a strict $*$-homomorphism from $A$ to $A\underset{B}{_{\beta}*_{\alpha}}A$ such that, for all $b\in B$, $\Gamma(\alpha(b))=\alpha(b)\underset{B}{_{\beta}\otimes_{\alpha}}1$, $\Gamma(\beta(b))=1\underset{B}{_{\beta}\otimes_{\alpha}}\beta(b)$, and :
 \[(\Gamma\underset{B}{_{\beta}*_{\alpha}}id)\Gamma=(id\underset{B}{_{\beta}*_{\alpha}}\Gamma)\Gamma\]
where  $\Gamma\underset{B}{_{\beta}*_{\alpha}}id$ and $id\underset{B}{_{\beta}*_{\alpha}}\Gamma$ are defined using \ref{def5}. Such an application is called a \emph{coproduct}.  

Then, by definition of $A\underset{B}{_{\beta}*_{\alpha}}A$, for any $x\in A$, $\Lambda_1$, $\Lambda_2$ in $L_T$, $\Lambda'_1$, $\Lambda'_2$ in $L_{T'}$, $\rho_{\Lambda_2}^*\Gamma(x)\rho_{\Lambda_1}$ and $\lambda_{\Lambda'_2}^*\Gamma(x)\lambda_{\Lambda'_1}$ belong to $A$. We suppose, moreover, that the two linear subsets generated by these elements are dense in $A$. Then, $\Gamma$ is called a \emph{non degenerate coproduct}. 

 (v) the modular automorphism groups of the KMS weights $\nu\circ\alpha^{-1}\circ T$ and $\nu^o\circ\beta^{-1}\circ T'$ on $A$ commute;
 
 (vi) $T$ is left-invariant, which means that, for any $x\in\gM_T^+$, we have :
 \[(id\underset{B}{_{\beta}*_{\alpha}}T)\Gamma(x)=T(x)\underset{B}{_{\beta}\otimes_{\alpha}}1\]

 where $(id\underset{B}{_{\beta}*_{\alpha}}T)$ had been defined in \ref{not}; 
 
 (vii) $T'$ is right-invariant, which means that, for any $x'\in\gM_{T'^+}$, we have :
 \[(T'\underset{B}{_{\beta}*_{\alpha}}id)\Gamma(x)=1\underset{B}{_{\beta}*_{\alpha}}T'(x')\]
 where $(T'\underset{B}{_{\beta}\otimes_{\alpha}}id)$ had been defined in \ref{not}.
 
 %%%%%trivial
 \subsubsection{\bf{Example}}
 \label{trivial}
 Let us suppose that $\gb_\nu=\gt$; then $T$ and $T'$ are KMS weights on $A$, $\bf C^*$-tensor product over $\gb_\nu$ is the usual tensor product of Hilbert spaces, then , thanks to \ref{C}, any reduced $\bf C^*$-algebraic quantum group $(A, \Gamma, T, T')$ in the sense of (\cite{KV1}, 4.1) is a locally compact quantum groupoid $(\bf C$$, A, id, id, id, T, T', \Gamma)$.

%%%%%from LCQG to MQG
\section{From a locally compact quantum groupoid to a measured quantum groupoid}
\label{fLCQG}

In this chapter, to any locally compact quantum groupoid $\bfG=(B, A, \alpha, \beta, \nu, T, T', \Gamma)$, we associate a measured quantum groupoid $\gG=(N, M, \underline{\alpha}, \underline{\beta}, \underline{\Gamma}, \underline{T}, \underline{T'}, \underline{\nu})$ such that $B$ is weakly dense in $N$, $A$ is weakly dense in $M$, $\alpha$ (resp. $\beta$, $\Gamma$, $T$, $T'$, $\nu$) is the restriction of $\underline{\alpha}$, (resp. $\underline{\beta}$, $\underline{\Gamma}$, $\underline{T}$, $\underline{T'}$, $\underline{\nu}$). The only problem is to extend the coproduct; this is done by using \cite{L}, in which Lesieur associated a pseudo-multiplicative unitary to any quantum measured groupoid; this construction was inspired by Kustermann-Vaes (\cite{KV1}), who associated a multiplicative unitary to the $\bf C^*$-version of any locally compact group; so, it is not a surprise that Lesieur's construction can be extended to a $\bf C^*$-context. Then, we shall say that $\bfG$ is a locally compact sub-quantum groupoid of $\gG$ (\ref{def2LCQG}). 

%%%not LCQG
\subsection{Notations}
\label{not LCQG}
Let $\bfG=(B, A, \alpha, \beta, \nu, T, T', \Gamma)$ be a locally compact quantum groupoid, as defined in \ref{defLCQG}; let us denote $\underline{\nu}$ the canonical normal faithful semi-finite weight constructed on $N=\pi_\nu(B)''$ in \ref{KMS}, and $\underline{\varphi}=\underline{\nu\circ\alpha^{-1}\circ T}$ the canonical faithful semi-finite weight constructed on $M=\pi_{\nu\circ \alpha^{-1}\circ T}(A)''$. Let $\underline{\alpha}$ (resp. $\underline{\beta}$) be the canonical injective $*$-homomorphism (resp. anti-$*$-homomorphism) from $N$ into $M$ constructed in \ref{thKMS}. Let $\underline{T}$ (resp. $\underline{T'}$) be the normal faithful semi-finite operator-valued weight from $M$ to $\underline{\alpha}(N)$ (resp. $\underline{\beta}(N)$) constructed in \ref{thKMS}. Finally, let us define, for all $x\in N$, $\widehat{\beta}(x)=J_{\underline{\varphi}}\underline{\alpha}(x^*)J_{\underline{\varphi}}$. 

%%%thfond
\subsection{Theorem}
\label{thfond}
{\it Let $\bfG=(B, A, \alpha, \beta, \nu, T, T', \Gamma)$ be a locally compact quantum groupoid, as defined in \ref{defLCQG}; let's use the notations of \ref{not LCQG}; then, 

(i) for any $a\in\gN_T\cap\gN_\varphi$, and any $v$, $\xi$, in $D((H_{\underline{\varphi}})_{\underline{\beta}}, \nu^o)$, $(\omega_{v, \xi}\underset{\underline{\nu}}{_{\underline{\beta}}*_{\underline{\alpha}}}id)\Gamma(a)$ belongs to $\gN_T\cap\gN_\varphi$; 

(ii) there exists an isometry $U$ from $H_{\underline{\varphi}}\underset{\underline{\nu}^o}{_{\underline{\alpha}}\otimes_{\widehat{\beta}}}H_{\underline{\varphi}}$ to $H_{\underline{\varphi}}\underset{\underline{\nu}}{_{\underline{\beta}}\otimes_{\underline{\alpha}}}H_{\underline{\varphi}}$, such that, for any orthogonal $(\underline{\beta}, \nu^o)$-basis $(\xi_i)_{i\in I}$ of $H_{\underline{\varphi}}$, any $a\in\gN_T\cap\gN_\varphi$ and any $v\in D((H_{\underline{\varphi}})_{\underline{\beta}}, \nu^o)$, we have :}
\[U(v\underset{\underline{\nu}^o}{_{\underline{\alpha}}\otimes_{\widehat{\beta}}}\Lambda_\varphi(a))=\sum_{i\in I}\xi_i\underset{\underline{\nu}}{_{\underline{\beta}}\otimes_{\underline{\alpha}}}\Lambda_{\underline{\varphi}}((\omega_{v, \xi_i}\underset{\underline{\nu}}{_{\underline{\beta}}*_{\underline{\alpha}}}id)\Gamma(a))\]

\begin{proof} The proof of (i) is identical to (\cite{L}, 3.2.7), and the proof of (ii) is identical to (\cite{L}, 3.2.9 and 3.2.10). \end{proof}

%%%propU
\subsection{Proposition}
\label{propU}
{\it Let's use the notations of \ref{not LCQG} and \ref{thfond}; then, we have :

(i) for any $x\in \gN_\varphi$, $e\in\gN_{\underline{\varphi}}$, $\xi\in H$, we have :
\[(1\underset{N}{_{\underline{\beta}}\otimes_{\underline{\alpha}}}J_{\underline{\varphi}}eJ_{\underline{\varphi}})U(\xi\underset{\underline{\nu}}{_{\underline{\alpha}}\otimes_{\widehat{\beta}}}\Lambda_\varphi (x))=\Gamma(x)(\xi\underset{\underline{\nu}}{_{\underline{\beta}}\otimes_{\underline{\alpha}}}J_{\underline{\varphi}}\Lambda_{\underline{\varphi}}(e))\]

(ii) for any $a\in A$, we have $\Gamma(a)U=U(1\underset{N}{_{\underline{\alpha}}\otimes_{\widehat{\beta}}}a)$. }

\begin{proof} The proof of (i) is identical to (\cite{L}, 3.3.1); then (ii) is a direct corollary of (i), as in (\cite{L}, 3.4.5). \end{proof}

%%%%lemU
\subsection{Lemma}
\label{lemU}
{\it Let's use the notations of \ref{not LCQG} and \ref{thfond}; then, we have, for any $a\in\gN_\varphi$ :

(i) for any $v\in D(_{\underline{\alpha}}H, \underline{\nu})\cap D(H_{\underline{\beta}}, \underline{\nu}^o)$, $w\in D(H_{\underline{\beta}}, \underline{\nu}^o)$, we have :
\[(\omega_{v,w}*id)(U)\Lambda_{\varphi}(a)=\Lambda_{\underline{\varphi}}((\omega_{v,w}\underset{N}{_{\underline{\beta}}*_{\underline{\alpha}}}id)\Gamma(a))\]

(ii) for any $e\in\gN_{\underline{\varphi}}$, $\eta\in D(_{\underline{\alpha}}H, \underline{\nu})$, we have :
\[(id\underset{N}{_{\underline{\beta}}*_{\underline{\alpha}}}\omega_{J_{\underline{\varphi}}\Lambda_{\underline{\varphi}}(e), \eta})\Gamma(a)=
(id*\omega_{\Lambda_{\varphi}(a), J_{\underline{\varphi}}e^*J_{\underline{\varphi}}\eta})(U)\]
}

\begin{proof} The proof of (i) is identical to (\cite{L}, 3.3.3) and the proof of (ii) is identical to (\cite{L}, 3.3.4). \end{proof}

%%%thunitary
\subsection{Theorem}
\label{thunitary}
{\it Let's use the notations of \ref{not LCQG} and \ref{thfond}; the isometry $U$ is a unitary }

\begin{proof}
This is the difficult part of the result; moreover, the proof is identical to (\cite{L}, 3.5). \end{proof}

%%%%widetildeGamma
\subsection{Theorem}
\label{widetildeGamma}
{\it Let's use the notations of \ref{not LCQG} and \ref{thfond}; for any $x\in M$, let's define $\underline{\Gamma}(x)=U(1\underset{N}{_{\underline{\alpha}}\otimes_{\underline{\beta}}}x)U^*$; then, $\underline{\Gamma}(x)$ belongs to $M\underset{N}{_{\underline{\beta}}*_{\underline{\alpha}}}M$, and $(N, M, \underline{\alpha}, \underline{\beta}, \underline{\Gamma})$ is a Hopf-bimodule in the sense of \ref{defHopf}.}

\begin{proof} For any $a\in A$, we have, thanks to \ref{propU}(ii) and \ref{thunitary}, $\underline{\Gamma}(a)=\Gamma(a)$, so, by continuity and \ref{thfiber}, we get that $\underline{\Gamma}(x)$ belongs to $M\underset{N}{_{\underline{\beta}}*_{\underline{\alpha}}}M$. We get also that $\underline{\Gamma}$ is a coproduct by continuity. \end{proof}

%%%%thMQG
\subsection{Theorem}
\label{thMQG}
{\it Let $\bfG=(B, A, \alpha, \beta, \nu, T, T', \Gamma)$ be a locally compact quantum groupoid; let's use the notations of \ref{not LCQG} and \ref{widetildeGamma}. Then, $\gG=(N, M, \underline{\alpha}, \underline{\beta}, \underline{\Gamma}, \underline{T}, \underline{T'}, \underline{\nu})$ is a measured quantum groupoid, we shall denote  $\gG(\bfG)$, and we have $\alpha=\underline{\alpha}|_B$, $\beta=\underline{\beta}|_B$, $\Gamma=\underline{\Gamma}|_A$, $T=\underline{T}|_A$, $T'=\underline{T'}|_A$, $\nu=\underline{\nu}|_B$. }

\begin{proof} By normality of $\underline{T}$, we easily get that $\underline{T}$ is left-invariant with respect to $\underline{\Gamma}$; the commutation property of  $\sigma^{\underline{\nu}\circ\underline{\alpha}^{-1}\circ\underline{T}}$ and $\sigma^{\underline{\nu}^o\circ\underline{\beta}^{-1}\circ\underline{T'}}$ are clear, by continuity and \ref{defLCQG}(v). \end{proof}

Using again (\cite{L}, 3.6) we get easily that the pseudo-multiplicative unitary $W$ associated to $\gG$ is equal to $U^*$.  

%%%%thW
\subsection{Theorem}
\label{thW}
{\it Let's use the notations of \ref{thMQG}. Then :

(i) let $a$, $b$, $e$ in $\gN_\varphi\cap\gN_T$; we heve :
\[(id*\omega_{J_{\underline{\varphi}}\Lambda_{\varphi}(e^*b), \Lambda_{\varphi} (a)})(W)
=(id\underset{N}{_{\underline{\beta}}*_{\underline{\alpha}}}\omega_{J_{\underline{\varphi}}\Lambda_\varphi(b), J_{\underline{\varphi}}\Lambda_\varphi (e)})\Gamma(a^*)\]

(ii) the linear set generated by all elements of the form $(id*\omega_{J_{\underline{\varphi}}\Lambda_\varphi(b), \Lambda_\varphi(a)})(W)$, for $a$, $b$ in $\gN_\varphi\cap\gN_T$ is norm dense in $A$. Moreover, $A$ is a sub-$\bf C^*$-algebra of the $\bf C^*$-algebra $A_n(W)\cap A_n(W)^*$ introduced in \ref{AnW}. 

(iii) the co-inverse $R$ of $\gG$ is such that $R(A)=A$. 

(iv) the modulus $\delta$ and the scaling operator $\lambda$ of the measured quantum groupoid $\gG$ are affiliated to $A$, in the sense of (\cite{B}, \cite{W4}). 

(v) The scaling group $\tau_t$ of $\gG$ satisfies, for all $t\in\bf R$, $\tau_t(A)=A$, and $\tau_{|A}$ is a norm continuous one parameter group of automorphims of $A$.

}
\begin{proof}
Formula (i) is just an application of \ref{lemU}(ii) and \ref{defLCQG}(iv). Then, using the notations of \ref{AnW}, we have $A\subset A_n(W)$, from which we get (ii). 

We have, using (\cite{E4}, 3.10(v) and 3.11):
\[R((id*\omega_{J_{\underline{\varphi}}\Lambda_\varphi(b), \Lambda_\varphi(a)})(W))=J_{\widehat{(\underline{\varphi})}}(id*\omega_{J_{\underline{\varphi}}\Lambda_\varphi(b), \Lambda_\varphi(a)})(W)^*J_{\widehat{(\underline{\varphi})}}=(id*\omega_{J_{\underline{\varphi}}\Lambda_\varphi(a), \Lambda_\varphi(b)})(W)\]
from which, using (ii),  we get (iii). For the modulus, the proof of (iv) is completely similar to (\cite{KV1} 7.10); then, as $\lambda^{ist}=\sigma_t^{\varphi}(\delta^{is})\delta^{-is}$ (\cite{V1} 5.1), we get (iv). 

Using (\cite{E4}, 3.8 (ii)), and (i), we get :
\begin{eqnarray*}
\tau_t[(id*\omega_{J_{\underline{\varphi}}\Lambda_{\varphi}(e^*b), \Lambda_{\varphi} (a)})(W)]
&=&\tau_t[(id\underset{N}{_{\underline{\beta}}*_{\underline{\alpha}}}\omega_{J_{\underline{\varphi}}\Lambda_\varphi(b), J_{\underline{\varphi}}\Lambda_\varphi (e)})\Gamma(a^*)]\\
&=&(id\underset{N}{_{\underline{\beta}}*_{\underline{\alpha}}}\omega_{J_{\underline{\varphi}}\Lambda_\varphi(b), J_{\underline{\varphi}}\Lambda_\varphi (e)}\circ \sigma_{-t}^{\underline{\varphi}})\Gamma(\sigma_t^{\underline{\varphi}}(a^*))\\
&=&(id\underset{N}{_{\underline{\beta}}*_{\underline{\alpha}}}\omega_{J_{\underline{\varphi}}\Lambda_\varphi(\sigma_t^{\underline{\varphi}}(b)), J_{\underline{\varphi}}\Lambda_\varphi (\sigma_t^{\underline{\varphi}}(e))})\Gamma(\sigma_t^{\underline{\varphi}}(a^*))\\
&=&
(id*\omega_{J_{\underline{\varphi}}\Lambda_{\varphi}(\sigma_t^{\underline{\varphi}}(e^*b)), \Lambda_{\varphi}( \sigma_t^{{\underline{\varphi}}}(a))})(W)
\end{eqnarray*}
from which we get that $\tau_t[(id*\omega_{J_{\underline{\varphi}}\Lambda_{\varphi}(b), \Lambda_{\varphi} (a)})(W)]=(id*\omega_{J_{\underline{\varphi}}\Lambda_{\varphi}(\sigma_t^{\underline{\varphi}}(b)), \Lambda_{\varphi}( \sigma_t^{{\underline{\varphi}}}(a))})(W)$, which gives that $\tau_t(A)\subset A$, for any $t\in \bf R$, and, therefore $\tau_t(A)=A$. Moreover, as 
\[\|R^{\underline{\alpha}}(J_{\underline{\varphi}}\Lambda_{\underline{\varphi}}(\sigma_t^{\underline{\varphi}}(a))\|^2=\|\sigma_t^{\underline{\nu}\circ\underline{\alpha}^{-1}}(T(a^*a))\|\]
we get (iv), because $(\sigma_t^{\underline{\nu}})_{|M(B)}$ is norm continuous.  \end{proof}

%%%lemA
\subsection{Lemma}
\label{lemA}
{\it Let $\gG=(N,M,\alpha, \beta, \Gamma, T, RTR, \nu)$ be a measured quantum groupoid. Let $A$ be a $\bf C^*$-algebra, weakly dense in $M$, invariant by $R$, and $B$ a $\bf C^*$-algebra weakly dense in $N$, such that $T(\gM_T\cap A)\subset M(B)$ and let us suppose that $(\nu_{|B}, T_{|A})$ is KMS. Moreover, let us suppose that, for any $a\in A$, $b$, $c$ in $A\cap\gN_\varphi\cap\gN_T$, $(i\underset{N}{_\beta*_\alpha}\omega_{J_\varphi\Lambda_\varphi(b), J_\varphi\Lambda_\varphi(c)})\Gamma(a)$ belongs to $A$, and that $A$ is equal to the closed linear set generated by these elements. Then, for any $c_1$, $c_2$ in $A\cap\gN_\varphi\cap\gN_T$ and $c_2$ in $A\cap\gM_T\cap\gM_\varphi$, $(id*\omega_{\Lambda_{\varphi}(c_1), J_\varphi\Lambda_\varphi(c_2)})(W^*)$ belongs to $A$, and $A$ is the closed linear set generated by these elements. Moreover, if $d_1$, $d_2$ belong to $A\cap\gN_{\varphi\circ R}$, $(\omega_{J_{\varphi\circ R}\Lambda_{\varphi\circ R}(d_1), J_{\varphi\circ R}\Lambda_{\varphi\circ R}(d_2)}\underset{N}{_\beta*_\alpha}id)\Gamma(a)$ belongs to $A$; it belongs to $A\cap\gN_T\cap\gN_T^*\cap\gN_\varphi\cap\gN_{\varphi}^*$ if $a$ belongs to $A\cap\gN_T\cap\gN_T^*\cap\gN_\varphi\cap\gN_{\varphi}^*$. }

\begin{proof}
Using  \ref{defMQG}(iii), we get that $(id*\omega_{\Lambda_\varphi(x), J_\varphi\Lambda_\varphi(y_1^*y_2))}(W^*)=(id\underset{N}{_\beta*_\alpha}\omega_{J_\varphi\Lambda_\varphi(y_1), J_\varphi\Lambda_\varphi(y_2)})\Gamma(x)$ belongs to $A$, for any $x$, $y_1$, $y_2$ in $A\cap\gN_T\cap\gN_\varphi$; using now \ref{lemW}, we get the first result. The second result of is given by the hypothesis $R(A)=A$. 
\end{proof}

%%%def
\subsection{Definition}
\label{def}
Let $\gG=(N,M,\alpha, \beta, \Gamma, T, RTR, \nu)$ be a measured quantum groupoid. Using \ref{aut}, we can define $\bf C^*$$(\sigma^\varphi, \sigma^{\varphi\circ R}, \tau)$ as the weakly dense sub-$\bf C^*$-algebra of $M$ on which the one-parameter groups $\sigma_t^\varphi$, $\sigma_t^{\varphi\circ R}$, $\tau_t$ are norm continuous, and $\bf C^*(\sigma^\nu, \gamma)$ a the weakly dense sub-$\bf C^*$-algebra of $N$ on which the one parameter groups $\sigma_t^\nu$ and $\gamma_t$ are norm continuous. Then $\alpha(\bf C^*(\sigma^\nu, \gamma))$ and $\beta (\bf C^*(\sigma^\nu, \gamma))$ are included into $\bf C^*$$(\sigma^\varphi, \sigma^{\varphi\circ R}, \tau)$. With the notations of \ref{thMQG}, we get that $A\subset \bf C^*$$(\sigma^{\underline{\varphi}}, \sigma^{\underline{\varphi\circ R}}, \tau)$.

 %%%%%not
 \subsection{\bf{Notations}}
 \label{not2}
 Using \ref{thfiber}, we define $(id\underset{B}{_{\beta}*_{\alpha}}T_1)$ as the restriction of $(id\underset{N}{_{\underline{\beta}}*_{\underline{\alpha}}}\widetilde{T_1})$ to $A_2\underset{B}{_{\beta}*_{\alpha}}A_1$; it is a faithful lower semi-continuous $\bf C^*$-valued weight from $A_2\underset{B}{_{\beta}*_{\alpha}}A_1$ to $A_2\cap \beta(B)'\underset{B}{_{\beta}\otimes_{\alpha}}1$. The $\bf C^*$-valued weight $(T_2\underset{B}{_{\beta}*_{\alpha}}id)$ is defined the same way. 
 
 Using again \ref{thfiber}, we define also $(id\underset{B}{_{\beta}*_{\alpha}}\varphi)$ as the restriction of $(id\underset{N}{_{\underline{\beta}}*_{\underline{\alpha}}}\underline{\varphi})$ to $A_2\underset{B}{_{\beta}*_{\alpha}}A_1$. As $(id\underset{\nu}{_{\underline{\beta}}*_{\underline{\alpha}}}\underline{\nu}\circ\underline{\alpha}^{-1})$ is just the inclusion of $M_2\cap\underline{\beta}(N)'=M_2\underset{N}{_{\underline{\beta}}*_{\underline{\alpha}}}\underline{\alpha}(N)$ into $M_2$ (\cite{L}, 2.5.1), we get that, for any positive $X$ in $A_2\underset{B}{_{T_2}*_{T_1}}A_1$, we have 
 \[(id\underset{B}{_{\beta}*_{\alpha}}T_1)(X)=(id\underset{\nu}{_{\beta}*_{\alpha}}\varphi)(X)\underset{\underline{\nu}}{_{\underline{\beta}}\otimes_{\underline{\alpha}}}1\] 
 
 %%%%
 \subsection{\bf{Definition}}
\label{def2LCQG}
 Let $\mathfrak {G}=(N, M, \alpha, \beta, \Gamma, T, RTR, \nu)$ a measured quantum groupoid, with $\varphi=\nu\circ T$; let $B$ be a weakly dense sub-$\bf{C^*}$-algebra of $N$, $A$ a sub-$\bf{C^*}$-algebra of $M$ such that $\alpha(B)\subset M(A)$, $R(A)=A$ (and, therefore, $\beta(B)\subset M(A)$), $A\cap\gM_T$ is dense in $A$, and $T(A\cap\gM_T)\subset \alpha(M(B))$. Let us suppose that $\nu_|B$ is a KMS weight on $B$, $T_|A$ is a densely defined faithful lower semi-continuous $\bf C^*$-valued weight from $A$ to $\alpha(M(B))$, such that the pair $(T_{|A}, \varphi_{|A})$ is KMS in the sense of \ref{KMSC*}. Let us suppose that $\Gamma(A)\subset A\underset{B}{_{\beta}*_{\alpha}}A$ (which implies that $A\underset{B}{_{\beta}*_\alpha}A$ is non degenerate). When, we shall say that $(B, A, \alpha_{|B}, \beta_{|B}, \nu_{|B}, T_{|A}, RTR_{|A}, \Gamma_{|A})$ is a \emph{locally compact quantum groupoid}, or, more pricesely, a \emph{locally compact sub-quantum groupoid of $\mathfrak{G}$}.

 %%%%%trivial
 \subsubsection{\bf{Example}}
 \label{trivial}
 Let us suppose that $B=\bf{C}$; then $T$ and $T'$ are KMS weights on $A$, then any reduced $\bf C^*$-algebraic quantum group $(A, \Gamma, T, T')$ in the sense of (\cite{KV1}, 4.1) is a locally compact quantum groupoid $(\bf C$$, A, id, id, id, T, T', \Gamma)$.

%%%%Ex
\subsubsection{\bf{Example}}
\label{Ex}
Let $(\bf C$$, A, id, id, id, T, T', \Gamma)$ a locally compact quantum groupoid described in \ref{trivial}  Then, using \ref{thMQG}, we obtain that 
$(\bf C$$, \pi_T(A)'', id, id, \underline{\Gamma}, \underline{T}, \underline{T'}, id)$ is aa measured quantum groupoid. 
More precisely, $(\pi_T(A)'', \underline{\Gamma}, \underline{T}, \underline{T'})$ is (the von Neumann version) of a locally compact quantum group. Let now $(A_1, \underline{\Gamma}_{|A_1}, \underline{T}_{|A_1}, \underline{T'}_{|A_1})$ be the the $\bf C^*$ version of this locally compact quantum group. Using \ref{thW}(i) and (\cite{KV2}, 1.6), we get that $A= A_1$, and, therefore, that $\Gamma(A)\subset  \widetilde{M}(A\otimes A)$ and we get that that $(A, \Gamma, T, T')$ is (the $\bf C^*$-version of) a locally compact quantum group. So, we get the converse of \ref{trivial}. 
 
 %%%%QSQG
 \subsubsection{\bf{Example}}
 \label{QSNG}
 Let $A$ be a $\bf C^*$-algebra. Let $\nu$ be a densely defined lower semi-continuous faithful weight on $A$, which is KMS with respect to a norm continuous one parameter group of automorphisms $\sigma$ (\ref{KMS}). Let us denote $M=\pi_\nu(A)''$, and $\underline{\nu}$ the normal semi-finite faithful weight on $M$ which extends $\nu$. 
 Lesieur had introduced in (\cite{L}, 14) a quantum space quantum groupoid associated to any von Neumann algebra $M$, acting on the Hilbert space $L^2(M)$, we identify $M'\underset{Z(M)}{*}M$ with $M'\underset{Z(M)}{\otimes} M$. Let $tr$ a normal semi-finite trace on $Z(M)$, such that $tr_{|Z(M(A))}$ is a densely defined faithful lower semi-continuous KMS trace on $Z(M(A))$, and $T$ the normal semi-finite operator-valued weight from $M$ to $Z(M)$ such that $\underline{\nu}=tr\circ T$. Let $j$ be the anti-isomorphism of $\mathcal L(H_\nu)$ given by $j(x)=J_\nu x^*J_\nu$, for any $x\in\mathcal L(H_\nu)$, which will identify $M'$ with $M^o$. Let $\alpha$ (resp. $\beta$) be the representation (resp. anti-representation) of $M$ into $M^o\underset{Z(M)}{\otimes} M$, such that $\alpha(m)=1\underset{Z(M)}{\otimes}m$ (resp. $\beta(m)=j(m)\underset{Z(M)}{\otimes} 1$).  Let $I$ be the isomorphism from $[L^2(M)\underset{tr}{\otimes}L^2(M)]\underset{\nu}{_\beta\otimes_\alpha}[L^2(M)\underset{tr}{\otimes}L^2(M)]$ onto $L^2(M)\underset{tr}{\otimes}L^2(M)\underset{tr}{\otimes}L^2(M)$ defined by $I[\Lambda_\nu(y)\underset{tr}{\otimes}\eta]\underset{\nu}{_\beta\otimes_\alpha}\Xi=\alpha(y)\Xi\underset{tr}{\otimes} 1$. This isomorphism will allow us to identify the fiber product $(M^o\underset{Z(M)}{\otimes}M)\underset{M}{_\beta*_\alpha}(M^o\underset{Z(M)}{\otimes} M)$ with $M^o\underset{Z(M)}{\otimes}M$; then we define $\Gamma$ on $M^o\underset{Z(M)}{\otimes}M$ by ($n, m \in M$) :
 \[\Gamma(n^o\underset{Z(M)}{\otimes} m)=[1\underset{Z(M)}{\otimes} m]\underset{\nu}{_\beta\otimes_\alpha}[n^o\underset{Z(M}{\otimes}1]\]
 Using the previous isomorphism, $\Gamma$ is just the identity. 
 
 The co-inverse $R$ is then given by $R(n^o\underset{Z(M)}{\otimes} m)=m^o\underset{Z(M)}{\otimes}j(n)$. Then Lesieur proved that $(M, M^o\underset{Z(M)}{\otimes} M, \alpha, \beta, \Gamma, \underline{\nu}, R(id\underset{Z(M)}{*}T)R, id\underset{Z(M)}{*}T)$ is a measured quantum groupoid he called a "quantum space quantum groupoid". 
 
 Let us consider now $(A, A^o\underset{\gb}{*}A, \alpha_{|A}, \beta_{|A}, tr_{|A}, R(id\underset{Z(M)}{*}T)R_{|A}, id\underset{Z(M)}{*}T_{|A}, \Gamma_{|A})$, where $\gb$ is the $\bf C^*$-base constructed using $tr_{|Z(A)}$. It is easy to verify that it is a locally compact sub-quantum groupoid of Lesieur's quantum space quantum groupoid. 
 %%%%PQG
 \subsubsection{\bf{Example}}
 \label{PQG}
 Let use again the notations of \ref{QSNG}. Lesieur had introduced in (\cite{L}, 15) a measured quantum groupoid, called a "pair quantum groupoid" associated to any von Neumann algebra $M'\otimes M$ ($M'$ being the commutant of $M$ in $L^2(M)$; let $\underline{\nu'}$ be the canonical weight on $M'$ constructed from $\underline{\nu}$.). 
 
 The fiber product $(M'\otimes M)\underset{\underline{\nu}}{_\beta *_\alpha}(M'\otimes M)$ can be identified with $M'\otimes Z(M)\otimes M$; then, with the nomal $*$-homomrphism $\Gamma$ defined by $\Gamma(n'\otimes m)=(1\otimes m)\underset{\underline{\nu}}{_\beta\otimes_\alpha}(n'\otimes 1)$, we obtain a measured quantum groupoid $(M, M'\underset{Z(M)}\otimes M, \alpha, \beta, \Gamma, \underline{\nu}, (\underline{\nu'}\otimes id), (id\otimes\underline{\nu}))$. Then, it is clear that :
 \[(A, J_{\underline{\nu}}AJ_{\underline{\nu}}\underset{Z(M)}{\otimes} A, \alpha_{|A}, \beta_{|A}, \Gamma_{|(J_{\underline{\nu}}AJ_{\underline{\nu}}\underset{Z(M)}{\otimes} A)}, \nu, (\nu'\otimes id), (id\otimes\nu))\] is a locally compact sub-quantum groupoid of Lesieur's pairs quantum groupoid.

%%%%from MQG to LCQG
\section{From a measured quantum groupoid to a canonical locally compact sub-quantum groupoid}
\label{fMQG}

In this chapter, starting from a measured quantum groupoid, we construct a sub $\bf C^*$-algebra which has a structure of a locally compact quantum groupoid (\ref{ThLCQG}). But, the examples \ref{QSNG} and \ref{PQG} show that this object may be not unique.

%%%Notations
\subsection{Notations}
\label{Notations}
Let $\mathfrak {G}=(N, M, \alpha, \beta, \Gamma, T, T', \nu)$ be a measured quantum groupoid; we shall use all notations of chapter \ref{MQG}; using \ref{aut}, we can define $\bf C^*$$(\sigma^\varphi, \sigma^{\varphi\circ R}, \tau)$ as the weakly dense sub-$\bf C^*$-algebra of $M$ on which the one-parameter groups $\sigma_t^\varphi$, $\sigma_t^{\varphi\circ R}$, $\tau_t$ are norm continuous, and $\bf C^*(\sigma^\nu, \gamma)$ a the weakly dense sub-$\bf C^*$-algebra of $N$ on which the one parameter groups $\sigma_t^\nu$ and $\gamma_t$ are norm continuous. Then $\alpha(\bf C^*(\sigma^\nu, \gamma))$ and $\beta (\bf C^*(\sigma^\nu, \gamma))$ are included into $\bf C^*$$(\sigma^\varphi, \sigma^{\varphi\circ R}, \tau)$. Using \ref{rest}, let us define the $\bf C^*$-algebra $\bf C^*(\mathfrak G)$ as $\bf C^*$$(T, \nu)\cap $$\bf C^*$$(RTR, \nu)\cap \bf C^*$$(\sigma^\varphi, \sigma^{\varphi\circ R}, \tau)$; using \ref{rest}, we see that $\bf C^*(\mathfrak G)$  contains $\gM_\varphi\cap\gM_{\varphi\circ R}\cap \bf C^*$$(\sigma^\varphi, \sigma^{\varphi\circ R}, \tau)$. As the modular groups $\sigma^\varphi$ and $\sigma^{\varphi\circ R}$ commute, we know that $\gM_\varphi\cap\gM_{\varphi\circ R}$ is weakly dense in $M$; using now \ref{aut}, we see that $\bf C^*(\mathfrak G)$ is weakly dense in $M$. Moreover, it is clear that $R($$\bf C^*(\mathfrak G)$)=$\bf C^*(\mathfrak G)$ and that $T(\gM_T\cap \bf C^*(\mathfrak G))\subset \alpha(\bf C^*(\nu, \gamma))$ (and $RTR(\gM_{RTR}\cap \bf C^*(\mathfrak G))\subset \beta(\bf C^*(\nu, \gamma))$, where $\bf C^*(\nu, \gamma)=\bf C^*(\nu)\cap\bf C^*(\gamma)$.

%%%lemmaGamma
\subsection{Lemma}
\label{lemmaGamma}
{\it Let $\mathfrak {G}=(N, M, \alpha, \beta, \Gamma, T, T', \nu)$ be a measured quantum groupoid; we shall use all notations of chapter \ref{MQG}; for any $x$, $y$ in $\gN_\varphi\cap\gN_T$, we have :
\[\Gamma(x)\rho^{\beta, \alpha}_{J_\varphi\Lambda_\varphi(y)}=(1\underset{N}{_\beta\otimes_\alpha}J_\varphi yJ_\varphi)W^*\rho^{\alpha, \gamma}_{\Lambda_\varphi(x)}\]
}
\begin{proof}
Let $z\in\gN_\varphi\cap\gN_T$; using \ref{defW}(iii), we have :
\begin{eqnarray*}
(\rho^{\beta, \alpha}_{J_\varphi\Lambda_\varphi(z)})^*\Gamma(x)\rho^{\beta, \alpha}_{J_\varphi\Lambda_\varphi(y)}
&=&
(id\underset{N}{_\beta*_\alpha}\omega_{J_\varphi\Lambda_\varphi(y), J_\varphi\Lambda_\varphi(z)})\Gamma(x)\\
&=&
(id*\omega_{\Lambda_\varphi(x), J_\varphi\Lambda_\varphi(y^*z)})(W^*)\\
&=&
(\rho^{\beta, \alpha}_{J_\varphi\Lambda_\varphi(z)})^*(1\underset{N}{_\beta\otimes_\alpha}J_\varphi yJ_\varphi)W^*
\end{eqnarray*}
Taking now a basis $(e_i)_{i\in I}$ for $T$ (\ref{basisT}), we have :
\begin{eqnarray*}
\Gamma(x)\rho^{\beta, \alpha}_{J_\varphi\Lambda_\varphi(y)}
&=&\sum_i\rho^{\beta, \alpha}_{J_\varphi\Lambda_\varphi(e_i)}(\rho^{\beta, \alpha}_{J_\varphi\Lambda_\varphi(e_i)})^*\Gamma(x)\rho^{\beta, \alpha}_{J_\varphi\Lambda_\varphi(y)}\\
&=&\sum_i\rho^{\beta, \alpha}_{J_\varphi\Lambda_\varphi(e_i)}(\rho^{\beta, \alpha}_{J_\varphi\Lambda_\varphi(e_i)})^*(1\underset{N}{_\beta\otimes_\alpha}J_\varphi yJ_\varphi)W^*\rho^{\alpha, \widehat{\beta}}_{\Lambda_\varphi(x)}\\
&=&(1\underset{N}{_\beta\otimes_\alpha}J_\varphi yJ_\varphi)W^*\rho^{\alpha, \widehat{\beta}}_{\Lambda_\varphi(x)} \end{eqnarray*}\end{proof}

%%%%lemmaappunit2
\subsection{Lemma}
\label{lemmaappunit2}
{\it Let's use the notations of \ref{Notations}. 

(i) Let $a$, $b$ in $\gM_T\cap \bf C^*(\mathfrak G)$. The norm closure of the linear space generated by all elements of the form $\Lambda_T(a)\Lambda_T(b)^*$ is a $\bf C^*$-algebra $A$ and $\bf C^*(\mathfrak G)$$\subset M(A)$. 

(ii) we have $A\subset J_\varphi \alpha(N)'J_\varphi$. }

\begin{proof} For any $a_1$, $a_2$, $b_1$, $b_2$ in $A\cap\gM_T$, we have :
\[\Lambda_{T}(a_1)\Lambda_T(b_1)^*\Lambda_T(a_2)\Lambda_T(b_2)^*=\Lambda_T(a_1(T(b_1^*a_2)))\Lambda_T(b_2)*\]
By hypothesis, $T(b_1^*a_2)$ belongs to $M(\bf C^*(\mathfrak G))$, so $a_1T(b_1^*a_2)$ belongs also to $\bf C^*(\mathfrak G)\cap\gM_T$. Then, it is easy to get that the linear space generated by all elements of the form $\Lambda_T(a)\Lambda_T(b)^*$ is an algebra. The fact that $\bf C^*(\mathfrak G)$$\subset M(A)$ is trivial; which finishes the proof of (i).
Let now $n\in \gN_\nu$, analytic with respect to $\sigma^\nu$, such that $\sigma_{-i/2}^\nu(n^*)$ belongs to $\gN_\nu$, and $x\in \gN_\varphi$; we have $\Lambda_T(b)^*J_\varphi\alpha(n)J_\varphi\Lambda_\varphi(x)= T(b^*x)J_\nu\Lambda_\nu(n)$ and, therefore, 
\begin{eqnarray*}
\Lambda_T(a)\Lambda_T(b)^*J_\varphi\alpha(n)J_\varphi\Lambda_\varphi(x)
&=&
\Lambda_\varphi (aT(b^*x)\alpha(\sigma_{-/2}^\nu(n^*)))\\
&=&
J_\varphi\alpha(n)J_\varphi\Lambda_T(a)\Lambda_T(b)^*\Lambda_\varphi(x)
\end{eqnarray*}
from which we get, by continuity, that, for any $n\in N$
\[\Lambda_T(a)\Lambda_T(b)^*J_\varphi\alpha(n)J_\varphi=J_\varphi\alpha(n)J_\varphi\Lambda_T(a)\Lambda_T(b)^*\]
which gives (ii). \end{proof}
 
 %%%%propA
 \subsection{Proposition}
 \label{propA}
 {\it Let's use the notations of \ref{Notations}. Let $x$, $y$ in $\bf C^*(\mathfrak G)$$ \cap\gN_T\cap\gN_\varphi$; then the operator $(id*\omega_{J_\varphi\Lambda_\varphi(y), \Lambda_\varphi(x)})(W)$ belongs to $\bf C^*(\mathfrak G)$.  }
 
 \begin{proof}
 Using \ref{lemmaappunit}(i), we have $<J_\varphi\Lambda_\varphi((y), J_\varphi\Lambda_\varphi (y)>_{\alpha, \nu}=\alpha^{-1}T(y^*y)\in \bf{C^*(\nu, \gamma})$ and  
 
 $<\Lambda_\varphi(x), \Lambda_\varphi(x)>_{\gamma, \nu^o}^o=\alpha^{-1}T(x^*x)\in \bf{C^*(\nu, \gamma})$; therefore, using (\cite{E4}, 4.5 (ii), (iii) and (iv)) we get the result. \end{proof}

 %%%%proprho
 \subsection{Proposition}
 \label{proprho}
 {\it Let's use the notations of \ref{Notations} and \ref{lemmaappunit2}. Let $x$, $y$ in $\bf C^*(\mathfrak G)$$\cap\gN_T\cap\gN_\varphi$,  and $x_1$, $x_2$ in $\bf C^*(\mathfrak G)$$\cap\gM_T$; then, using again \ref{lemmaappunit2}(i), $\Lambda_T(x_1)\Lambda_T(x_2)^*y$ belongs to $A$;  let now $\Sigma_{\bf N}\Lambda_T(a_n)\Lambda_T(b_b)^*$ be an approximate unit of $A$; then $\Sigma_{\bf N}\Lambda_T(x_1)\Lambda_T(x_2)^*y\Lambda_T(a_n)\Lambda_T(b_n)^*$ is norm converging to $\Lambda_T(x_1)\Lambda_T(x_2)^*y$. Then:

(i) $(1\underset{N}{_\beta\otimes_\alpha}\Lambda_T(x_1)\Lambda_T(x_2)^*)\Gamma(x)\rho^{\beta, \alpha}_{J_\varphi\Lambda_\varphi(y)}$ is the norm limit of :
\[\sum_{\bf N} \rho^{\beta, \alpha}_{J_\varphi\Lambda_\varphi (x_1T(x_2^*ya_n))}(i*\omega_{J_\varphi\Lambda_\varphi(b_n), \Lambda_\varphi(x)})(W)^*\]}
 
 (ii) Let now $z\in \bf C^*(\mathfrak G)\cap\gM_T\cap\gM_\varphi$ and $X$, $Y$ in $A$; then $(1\underset{N}{_\beta\otimes_\alpha}X)\Gamma(x)(1\underset{N}{_\beta\otimes_\alpha}Y)$ belongs to $\bf C^*(\mathfrak G)$
 $\underset{C^*(\nu, \gamma)}{_{\beta_{|C^*(\nu, \gamma)}}*_{\alpha_{|C^*(\nu, \gamma)}}}$$B(H)$. 
 
 (iii) $\Gamma(\bf C^*(\mathfrak G))$ is included into $\bf C^*(\mathfrak G)$
 $\underset{C^*(\nu, \gamma)}{_{\beta_{|C^*(\nu, \gamma)}}*_{\alpha_{|C^*(\nu, \gamma)}}}$$B(H)$. 
 
 (iv) $\Gamma(\bf C^*(\mathfrak G))$ is included into  $\bf C^*(\mathfrak G)$
 $\underset{C^*(\nu, \gamma)}{_{\beta_{|C^*(\nu, \gamma)}}*_{\alpha_{|C^*(\nu, \gamma)}}}$$\bf C^*(\mathfrak G)$
 \begin{proof}
 Using \ref{lemmaGamma} and \ref{lemmaappunit}(iii), we get that $(1\underset{N}{_\beta\otimes_\alpha}\Lambda_T(x_1)\Lambda_T(x_2)^*)\Gamma(x)\rho^{\beta, \alpha}_{J_\varphi\Lambda_\varphi(y)}$ is the norm limit of :
 \begin{eqnarray*}
 \sum_{\bf N}(1\underset{N}{_\beta\otimes_\alpha}J_\varphi (1\underset{N}{_\beta\otimes_\alpha}\Lambda_T(x_1)\Lambda_T(x_2)^*)y\Lambda_T(a_n)\Lambda_T(b_n)^*)J_\varphi W^*\rho^{\alpha, \widehat{\beta}}_{\Lambda_\varphi(x)}=\\
 \sum_{\bf N}(1\underset{N}{_\beta\otimes_\alpha}J_\varphi\Lambda_T(x_1T(x_2^*ya_n))\Lambda_T(b_n)^*J_\varphi )W^*\rho^{\alpha, \widehat{\beta}}_{\Lambda_\varphi(x)}=\\
 \sum_{\bf N}\rho^{\beta, \alpha}_{\Lambda_\varphi(x_1T(x_2^*ya_n))}(id*\omega_{\Lambda_\varphi(x), J_\varphi\Lambda_\varphi(b_n)})(W)^*
 \end{eqnarray*}
which is (i). 
 
 Let now $x_3$, $x_4$ in $\bf C^*(\mathfrak G)\cap\gM_T$; we have :
\begin{eqnarray*}
(1\underset{N}{_\beta\otimes_\alpha}\Lambda_T(x_1)\Lambda_T(x_2)^*)\Gamma(z)(1\underset{N}{_\beta\otimes_\alpha}\Lambda_T(x_3)\Lambda_T(x_4)^*)\rho^{\beta, \alpha}_{J_\varphi\Lambda_\varphi(y)}=\\
(1\underset{N}{_\beta\otimes_\alpha}\Lambda_T(x_1)\Lambda_T(x_2)^*)\Gamma(z)\rho^{\beta, \alpha}_{J_\varphi\Lambda_\varphi(x_3T(x_2^*y))}
\end{eqnarray*}
So, applying (i) and \ref{def5}, we get that 
$(1\underset{N}{_\beta\otimes_\alpha}\Lambda_T(x_1)\Lambda_T(x_2)^*)\Gamma(z)(1\underset{N}{_\beta\otimes_\alpha}\Lambda_T(x_3)\Lambda_T(x_4)^*)$ belongs to $\bf C^*(\mathfrak G)$
 $\underset{\gb_\nu}{_{L_{RTR}}*_{L_T}}$$B(H)$. Then, by norm continuity, we get (ii). 

Let $X$ be in $A$, invertible in $\alpha(N)'$; then $1\underset{N}{_\beta\otimes_\alpha}X^{-1}$ belongs to $M(\bf C^*(\mathfrak G)$$
 \underset{B}{_{\beta}*_{\alpha}}$$B(H))$; and $\Gamma(z)=(\underset{N}{1_\beta\otimes_\alpha}X^{-1})(\underset{N}{1_\beta\otimes_\alpha}X)\Gamma(z)(\underset{N}{1_\beta\otimes_\alpha}X)(\underset{N}{1_\beta\otimes_\alpha}X^{-1})$ belongs to $\bf C^*(\mathfrak G)$$
 \underset{B}{_{\beta}*_{\alpha}}$$B(H)$. Then, by norm continuity, we get (iii). 
 
 Applying this result to $\mathfrak G^o$, we obtain that $\Gamma(\bf C^*(\mathfrak G^o)))$ is included into $B(H)$$\underset{B}{_{\beta}*_{\alpha}}$$\bf C^*(\mathfrak G^o)$; as $\bf C^*(\mathfrak G^o)=\bf C^*(\mathfrak G)$,  we get (iv). 
\end{proof}

%%%%ThLCQG
\subsection{Theorem}
\label{ThLCQG}
{\it Let's use the notations of \ref{Notations}. Then, the octuple 

\centerline{$(\bf C^*(\nu, \gamma), \bf C^*$$(\mathfrak G), \alpha_{|\bf C^*(\nu, \gamma)}, \beta_{|\bf C^*(\nu, \gamma)}, \nu_{|\bf C^*(\nu, \gamma)}, T_{|\bf C^*(\mathfrak G)}, RTR_{|\bf C^*(\mathfrak G)}, \Gamma_{|\bf C^*(\mathfrak G)})$} 
 is a locally compact quantum groupoid. }

\begin{proof}
We have $\alpha(\bf C^*(\nu, \gamma))$$\subset \bf C^*(\mathfrak G)$, and $\beta(\bf C^*(\nu, \gamma))$$\subset \bf C^*(\mathfrak G)$; by construction, $\nu_{|\bf C^* (\nu, \gamma)}$ is a KMS weight on $\bf C^*(\nu, \gamma)$, and $(T, \nu\circ\alpha^{-1})$ (resp. $(RTR, \nu^o\circ\beta^{-1})$) are KMS, in the sense of \ref{KMSC*}. And we have obtain in \ref{proprho}(iv) that 
$\Gamma_{|\bf C^*(\mathfrak G)}$ is a coproduct in the sense of \ref{LCQG}(iv). 

It is then clear also that the restriction of $T$ to $\bf C^*$$(\mathfrak G)$ is left-invariant, and the restriction of $RTR$ to $\bf C^*$$(\mathfrak G)$ is right-invariant, and that the modular groups of $\nu\circ\alpha^{-1}\circ T_{|\bf C^*(\mathfrak G)}$ and $\nu\circ\beta^{-1}\circ RTR_{|\bf C^*(\mathfrak G)}$ commute, which finishes the proof. \end{proof}

%%%%propLCQG
\subsection{Proposition}
\label{propLCQG}
{\it Let's use the notations of \ref{Notations}. The scaling operator $\lambda$ of $\gG$ (which is a non singular positive operator affiliated to $Z(M)$, of the form $\lambda=\alpha(q)=\beta(q)$, where $q$ is a non singular positive operator affiliated to $Z(N)$ (\cite{E3}, 3.8 (vi))), is affiliated to $Z(\bf C^*$$(\mathfrak G))$ (and $q$ is affiliated to $Z(\bf C^*(\nu, \gamma)))$) in the sense of (\cite{B}, \cite{W4})}

\begin{proof}
It is clear, using \ref{thW}(iv) applied to $\bf C^*$$(\mathfrak G)$. \end{proof}

%%%%%ex
\subsection{Remark}
\label{ex}
The examples \ref{QSNG} and \ref{PQG} show that, in general, a locally compact sub-quantum groupoid of a given measured quantum groupoid is not unique.

%%%%%%%%Duality
\section{Duality of locally compact groupoids}
\label{Duality}

In this chapter, starting from a locally compact quantum groupoid, and using then the construction of a measured quantum groupoid, we construct a dual locally compact quantum groupoid (\ref{thduality}), and prove that the bidual locally compact quantum groupoid is the initial one (\ref{propbidual}). 
%%%%%notd
\subsection{Notations}
\label{notd}
Let $\bfG=(B,A,\alpha, \beta, \nu, T, T', \Gamma)$ be a locally compact quantum groupoid, as defined in \ref{defLCQG}, and  $\gG=(N,M, \underline{\alpha}, \underline{\beta}, \underline{\Gamma}, \underline{T}, \underline{T'}, \underline{\nu})$ be the measured quantum groupoid constructed in \ref{thMQG}. As $B\subset N$, $A\subset M$, $\alpha=\underline{\alpha}_{|B}$, $\beta=\underline{\beta}_{|B}$, $\nu=\underline{\nu}_{|B}$, $\Gamma=\underline{\Gamma}_{|A}$, $T=\underline{T}_{|A}$, $T'=\underline{T'}_{|A}$, we shall use, for simplification of notations, $\alpha$, $\beta$, $\nu$, $\Gamma$, $T$, $T'$, instead of $\underline{\alpha}$, $\underline{\beta}$, $\underline{\nu}$, $\underline{T}$, $\underline{T'}$. Moreover, as indicated in \ref{data}(iii), we shall consider the measured quantum groupoid $\underline{\gG}$, and, therefore, assume that $T'=RTR$, where $R$ is the co-inverse of $\bfG$. 
\newline
As recalled in \ref{data}(viii), we consider now the dual measured quantum groupoid $\widehat{\gG}=(N, \widehat{M}, \alpha, \widehat{\beta}, \widehat{\Gamma}, \widehat{T}, \widehat{T'}, \nu)$, and we shall construct a dual locally compact quantum groupoid $\widehat{\bfG}=(B, \widehat{A}, \alpha_{|B}, \widehat{\beta}_{|B}, \nu_{|B}, \widehat{T}_{|\widehat{A}}, \widehat{T'}_{|\widehat{A}}, \widehat{\Gamma}_{|\widehat{A}})$, where $\widehat{A}$ is a sub-$\bf C^*$-algebra of $\widehat{M}$. 
\newline
In \ref{not2}, we had writen $\varphi=\nu\circ\alpha^{-1}\circ T$ the KMS weight on $A$, and $\underline{\varphi}=\underline{\nu\circ\alpha^{-1}\circ T}$ its extension to $M$; for simplification, let us write $\varphi$ instead of $\underline{\varphi}$, and $J$ instead of $J_{\underline{\varphi}}$. 
\newline
Let $\widehat{\varphi}$ be the canonical weight on the dual measured quantum groupoid $\widehat{\gG}$, and let us write $\widehat{J}$ for $J_{\widehat{\varphi}}$.

%%%%lembeta
\subsection{Lemma}
\label{lembeta}
{\it Let $x\in\gN_T\cap\gN_\varphi$; then $\widehat{J}J\Lambda_\varphi(x)$ belongs to $D(H_\beta, \nu^o)$.}

\begin{proof}
Let $n\in\gN_\nu$; we have :
\[\beta(n^*)\widehat{J}J\Lambda_\varphi(x)=\widehat{J}\alpha(n)J\Lambda_\varphi(x)=\widehat{J}J\Lambda_T(x)J_\nu\Lambda_\nu (n)\]
which gives the result. \end{proof}

%%%%prophat
\subsection{Proposition}
\label{prophat}
{\it (i) let $x$, $y$, in $A\cap\gN_T\cap\gN_\varphi$, and let $z=(\omega_{\widehat{J}J\Lambda_\varphi(x), J\Lambda_\varphi (y)}*id)(W)$; then $z$ belongs to $\gN_{\widehat{T}}\cap\gN_{\widehat{\varphi}}$, and we have $\Lambda_{\widehat{\varphi}}(z)=\widehat{J}JR(y)\Lambda_\varphi(x)$, and $\Lambda_{\widehat{T}}(z)=\widehat{J}JR(y)\Lambda_T(x)$. 
\newline
(ii) for $i=1,2$, let $x_i$, $y_i$, in $A\cap\gN_T\cap\gN_\varphi$, and let us write $z_i=(\omega_{\widehat{J}J\Lambda_\varphi(x_i), J\Lambda_\varphi (y_i)}*id)(W)$; then $z_i$ belongs to $\gN_{\widehat{T}}\cap\gN_{\widehat{\varphi}}$, and $\widehat{T}(z_2^*z_1)$ belongs to $\alpha(M(B))$. }

\begin{proof}
Using (\cite{E3}, 3.10 (v)); we get that $z$ belongs to $\gN_{\widehat{\varphi}}$, and :
 \[\Lambda_{\widehat{\varphi}}(z)=a_\varphi (\omega_{\widehat{J}J\Lambda_\varphi(x), J\Lambda_\varphi(y)})=Jy^*J\widehat{J}J\Lambda_\varphi(x)=\widehat{J}JR(y)\Lambda_\varphi(x)\]
Moreover, for any $n\in B\cap \gN_\nu$, analytic with respect to $\sigma^\nu$, and $x$ analytic with respect to $\sigma^\varphi$, we get :
\begin{eqnarray*}
z\alpha(n)&=&(\omega_{\beta(\sigma^\nu_{i/2}(n))\widehat{J}J\Lambda_\varphi(x), J\Lambda_\varphi(y)}*id)(W)\\
&=&(\omega_{\widehat{J}\alpha(\sigma_{-i/2}(n^*))J\Lambda_\varphi(x), J\Lambda_\varphi(y)}*id)(W)\\
&=&(\omega_{\widehat{J}\alpha(\sigma_{-i/2}(n^*))\Lambda_\varphi(\sigma_{-i/2}^\varphi(x^*), J\Lambda_\varphi(y)}*id)(W)\\
&=&(\omega_{\widehat{J}J\Lambda_\varphi(x\alpha(n)), J\Lambda_\varphi(y)}*id)(W)
\end{eqnarray*}
which remains true for any $x\in A\cap\gN_T\cap\gN_\varphi$, and any $n\in B\cap\gN_\nu$, from which we get that $z\alpha(n)$ belongs to $\gN_{\widehat{\varphi}}$, and $\Lambda_{\widehat{\varphi}}(z\alpha(n))=\widehat{J}JR(y)\Lambda_\varphi (x\alpha(n))=\widehat{J}JR(y)\Lambda_T(x)\Lambda_\nu(n)$, which finishes the proof of (i). 
Using (i),  we get that $\widehat{T}(z_2^*z_1)=T(x_2^*R(y_2^*)R(y_1)x_1)$ which belongs to $\alpha(M(B))$, which finishes the proof. 
\end{proof}
%%%prop2hat
\subsection{Proposition}
\label{prop2hat}
{\it Let $x$, $y$ in $A\cap\gN_T\cap\gN_\varphi$; we have :
\newline 
\[\widehat{R}((\omega_{\widehat{J}J\Lambda_\varphi(x), J\Lambda_\varphi(y)}*id)(W))=(\omega_{\widehat{J}J\Lambda_\varphi(y), J\Lambda_\varphi(x)}*id)(W)\]
(ii) For all $t\in \bf R$, $\delta^{it}\widehat{J}J\Lambda_\varphi(x)$ belongs to $D(H_\beta, \nu^o)$, and :
\[\delta^{it}((\omega_{\widehat{J}J\Lambda_\varphi(x), J\Lambda_\varphi(y)}*id)(W))\delta^{-it}=(\omega_{\delta^{it}\widehat{J}J\Lambda_\varphi(x), J\Lambda_\varphi (y)}*id)(W)\]
(iii) 
\[J\delta^{it}J((\omega_{\widehat{J}J\Lambda_\varphi(x), J\Lambda_\varphi(y)}*id)(W))J\delta^{-it}J=(\omega_{\widehat{J}J\Lambda_\varphi(x), \delta^{it}J\Lambda_\varphi(y)}*id)(W)\]
(iv)
\[P^{it}((\omega_{\widehat{J}J\Lambda_\varphi(x), J\Lambda_\varphi(y)}*id)(W))P^{-it}=(\omega_{\widehat{J}JP^{it}\Lambda_\varphi(x), JP^{-it}\Lambda_\varphi(y)}*id)(W)\]
(v) \[\sigma_t^{\widehat{\varphi}}((\omega_{\widehat{J}J\Lambda_\varphi(x), J\Lambda_\varphi(y)}*id)(W))=(\omega_{\widehat{J}JP^{it}\Lambda_\varphi(x), \delta^{it}JP^{-it}\Lambda_\varphi(y)}*id)(W)\]

}
\begin{proof}
We have :
\[\widehat{R}((\omega_{\widehat{J}J\Lambda_\varphi(x), J\Lambda_\varphi(y)}*id)(W))=J(\omega_{\widehat{J}J\Lambda_\varphi(x), J\Lambda_\varphi(y)}*id)(W)^*J=(\omega_{\widehat{J}J\Lambda_\varphi(y), J\Lambda_\varphi(x)}*id)(W)\]
which gives (i).
\newline
For all $n\in N$, we have :
\begin{eqnarray*}
\beta(n)\delta^{it}\widehat{J}J\Lambda_\varphi(x)
&=&\delta^{it}\sigma_t^{\varphi\circ R}\sigma_{-t}^\varphi(\beta(n))\widehat{J}J\Lambda_\varphi(x)\\
&=&\delta^{it}\beta(\sigma_t^\nu\gamma_{-t}(n))\widehat{J}J\Lambda_\varphi(x)\\
&=&\delta^{it}\widehat{J}\alpha(\sigma_t^\nu\gamma_{-t}(n))J\Lambda_\varphi(x)\\
&=&\delta^{it}\widehat{J}J\Lambda_T(x)J_\nu\Delta_\nu^{it}H^{-it}\Lambda_\nu(n)
\end{eqnarray*}
where $H^{it}$ is the canonical implementation of $\gamma_t$ on $H_\nu$. Which gives that $\delta^{it}\widehat{J}J\Lambda_\varphi(x)$ belongs to $D(H_\beta, \nu^o)$. Moreover, using now that $\Gamma(\delta^{it})=\delta^{it}\underset{N}{_\beta\otimes_\alpha}\delta^{it}$ (\cite{E3},3.8 (vi)), we finish the proof of (ii). Then (iii) is obtained easily from (ii) and (i). Using then $W(P^{it}\underset{N}{_\beta\otimes_\alpha}P^{it})=(P^{it}\underset{N^o}{_\alpha\otimes_\gamma}P^{it})(W)$ (\cite{E3}, 3.8 (vii)), we get (iv). Then (v) is obtained using (iii), (iv) and \cite{E3}, 3.10 (vii). \end{proof}
 
 %%%propA
\subsection{Proposition}
\label{propA}
{\it Let's use the notations of  \ref{Notations}. Let $\mathcal A_0$ be the set of elements $x$ in $A\cap\gM_T\cap\gM_{RTR}\cap\gM_\varphi\cap\gM_{\varphi\circ R}$, which are analytic with respect to $\sigma_t^\varphi$, $\sigma_t^{\varphi\circ R}$ and $\tau_t$, and such that, for any $z, z', z"$ in $\bf C$, $\sigma_z^\varphi\circ\sigma_{z'}^{\varphi\circ R}\circ\tau_{z"}(x)$ belongs to $A\cap\gN_T\cap\gN_{RTR}\cap\gN_\varphi\cap\gN_{\varphi\circ R}$. Then, $\mathcal A_0$ is a $*$-subalgebra of $A$, invariant by $R$, norm dense in $A$. Moreover, $\Lambda_\varphi(\mathcal A_0)$ is dense in $H_\varphi$. Moreover, if $\mathcal A_0^2$ denotes the linear span generated b products $xy$, where $x$, $y$ belong to $\mathcal A_0$, then $\mathcal A_0^2$ is a sub $*$-algebra of $\mathcal A_0$, invariant by $R$, norm dense in $A$, and such that $\Lambda_\varphi(\mathcal A_0^2)$ is dense in $H_\varphi$. Moreover, if $b\in B$ is analytic with respect to $\sigma_t^\nu$ and $\gamma_t$, we get that $\alpha(b)\mathcal A_0\subset \mathcal A_0$ and $\mathcal A_0\beta(b)\subset \mathcal A_0$. }

\begin{proof} Let's have a closer look at (\cite{L}, 6.0.10) and its proof. For any $x\in A\cap\gM_T\cap\gM_{RTR}\cap\gM_\varphi\cap\gM_{\varphi\circ R}$, we obtain a sequence $x_n$ in $\mathcal A_0$, which is norm converging to $x$, by similar arguments as in \ref{aut}. Then, we can prove that $A\cap\gM_T\cap\gM_{RTR}\cap\gM_\varphi\cap\gM_{\varphi\circ R}$ is norm dense in $A$, using \ref{defLCQG} (iii) and (vi). \end{proof}

%%%%%propA2
\subsection{Lemma}
\label{propA2}
{\it Let's use the notations of \ref{Notations} and \ref{propA}. The linear set generated by all elements of the form $(\omega_{\widehat{J}J\Lambda_\varphi(x), J\Lambda_\varphi (y)}*id)(W)$ with $x$, $y$ in $\mathcal A_0$ is invariant by $*$. }

 \begin{proof}
 We have, using \ref{prophat}, $\Lambda_{\widehat{\varphi}}((\omega_{\widehat{J}J\Lambda_\varphi(x), J\Lambda_\varphi (y)}*id)(W))=\widehat{J}JR(y)\Lambda_\varphi(x)$. Therefore, $(\omega_{\widehat{J}J\Lambda_\varphi(x), J\Lambda_\varphi (y)}*id)(W)$ belongs to $\gN_{\widehat{\varphi}}^*$ if and only if $JR(y)\Lambda_\varphi(x)$ belongs to $\mathcal D(\widehat{\Delta}^{-1/2})$. But, if $x$, $y$ belong to $\mathcal A_0$, we get, using successively \cite{E3}3.10(vii), \cite{V1}, \cite{E3} 3.8(vii), \cite{E3}3.10(v) that :
 \begin{eqnarray*}
 \widehat{\Delta}^{-1/2}JR(y)\Lambda_\varphi(x)
 &=&
 P^{-1/2}J\delta^{1/2}R(y)\Lambda_\varphi(x)\\
 &=&
 P^{-1/2}\lambda^{-i/4}\Lambda_\varphi(\sigma^\varphi_{-i/2}(x^*R(y^*))\delta^{1/2})\\
 &=&
 \Lambda_\varphi(\tau_{i/2}\sigma^\varphi_{-i/2}(x^*R(y^*))\delta^{1/2})\\
 &=&
 \widehat{J}\Lambda_\varphi(R(\tau_{i/2}\sigma^\varphi_{-i/2}(x^*R(y^*))^*)\\
 &=&
 \widehat{J}\Lambda_\varphi(\tau_{-i/2}\sigma^{\varphi\circ R}_{i/2}(R(x)y)\\
 &=&
 \widehat{J}J\Lambda_\varphi(\sigma_{-i/2}\tau_{i/2}\sigma^{\varphi\circ R}_{-i/2}(y^*R(x^*))
 \end{eqnarray*}
 which is equal to $\Lambda_{\widehat{\varphi}}(\omega_{\widehat{J}J\Lambda_\varphi(x'), J\Lambda_\varphi (y')}*id)(W))$, with $x'=\sigma_{-i/2}^\varphi\tau_{i/2}\sigma_{-i/2}^{\varphi\circ R}(R(x^*))$ and $y'=\sigma_{-i/2}^\varphi\tau_{i/2}\sigma_{-i/2}^{\varphi\circ R}(y^*)$, which both belong to $\mathcal A_0$. Therefore, we get that :
\[(\omega_{\widehat{J}J\Lambda_\varphi(x), J\Lambda_\varphi (y)}*id)(W)^*= (\omega_{\widehat{J}J\Lambda_\varphi(x'), J\Lambda_\varphi (y')}*id)(W)^*\]
which finishes the proof. \end{proof}

 %%%lemproduct
 \subsection{Lemma}
 \label{lemproduct}
 {\it Let $X$ in $\gN_{\widehat{T}}\cap\gN_{\widehat{\varphi}}$ such that there exists $Y\in A\cap\gN_T\cap \gN_{\varphi}$ with $\Lambda_{\widehat{T}}(X)=\widehat{J}J\Lambda_T(Y)$ and $\Lambda_{\widehat{\varphi}}(X)=\widehat{J}J\Lambda_{\varphi}(Y)$ (for instance, using \ref{prophat}, any linear combination of operators of the form $(\omega_{\widehat{J}J\Lambda_\varphi(x_1), J\Lambda_\varphi (y_1)}*id)(W)$ with $x_1$, $y_1$, in $A\cap\gN_T\cap\gN_\varphi)$. Then, for any $x_2$, $y_2$ in $A\cap\gN_T\cap \gN_{RTR}\cap\gN_\varphi\cap\gN_{\varphi\circ R}$, $X(\omega_{\widehat{J}J\Lambda_\varphi(x_2), J\Lambda_\varphi (y_2)}*id)(W)$ satisfies the same property. }
 
 \begin{proof}
 Let us suppose first that $x_2$, $y_2$ are analytics with respect to $\tau_t$ and that $\tau_{-i/2}(x_2)$ and $R(\tau_{i/2}(y_2)^*)$ belong to $\gN_\varphi$. 
Using (\cite{E3}3.1.(iv) and 3.10.5), \cite{V1} 5.1, we get that :
\[\widehat{J}J\Lambda_\varphi(R(\tau_{i/2}(y_2^*)))=J\lambda^{i/4}\Lambda_\varphi(\tau_{-i/2}(y_2)\delta^{1/2})=\delta^{1/2}J\Lambda_\varphi(\tau_{-i/2}(y_2))=\lambda^{-1/4}\delta^{1/2}JP^{1/2}\Lambda_\varphi(y_2)\]
and, using \cite{E3}3.8(vii) :
\[\widehat{J}J\Lambda_\varphi(\tau_{-i/2}(x_2))=\lambda^{1/2}\widehat{J}JP^{1/2}\Lambda_\varphi(x_2)\]
Using then \ref{prop2hat} (v), we get :
\[\sigma_{-i/2}^{\widehat{\varphi}}((\omega_{\widehat{J}J\Lambda_\varphi(x_2), J\Lambda_{\varphi}(y_2)}*id)(W)^*)=\lambda^{-i/4}(\omega_{\widehat{J}J\Lambda_\varphi(\tau_{-i/2}(x_2), \widehat{J}J\Lambda_\varphi(R(\tau_{i/2}(y^*)))}*id)(W)^*\]
and, using \cite{E3} 3.8 (ii) :
\begin{eqnarray*}
\Lambda_{\widehat{\varphi}}(X(\omega_{\widehat{J}J\Lambda_\varphi(x_2), J\Lambda_\varphi (y_2)}*id)(W))
&=&
\widehat{J}\sigma_{-i/2}^{\widehat{\varphi}}((\omega_{\widehat{J}J\Lambda_\varphi(x_2), J\Lambda_{\varphi}(y_2)}*id)(W)^*)\widehat{J}\Lambda_{\widehat{\varphi}}(X)\\
&=&
\widehat{J}\lambda^{-i/4}(\omega_{\widehat{J}J\Lambda_\varphi(\tau_{-i/2}(x_2)), \widehat{J}J\Lambda_\varphi(R(\tau_{i/2}(y_2^*)))}*id)(W)^*J\Lambda_\varphi(Y)\\
&=&
\widehat{J}\lambda^{-i/4}\Lambda_\varphi(\omega_{\widehat{J}J\Lambda_\varphi(R(\tau_{i/2}(y_2^*))), \widehat{J}J\Lambda_\varphi(\tau_{-i/2}(x_2)))}\underset{N}{_\beta*_\alpha}id)\Gamma(\sigma^\varphi_{-i/2}(Y^*))\\
&=&
\widehat{J}J\Lambda_\varphi((\omega_{\widehat{J}J\Lambda_\varphi(R(y^*_2)), \widehat{J}J\Lambda_\varphi(x_2)}\underset{N}{_\beta*_\alpha}id)\Gamma(Y))
\end{eqnarray*}
and, thanks to \ref{defLCQG} (v) and \ref{thW}(iv), we get the result. \end{proof}

%%%%Ahat
\subsection{Notations}
\label{Ahat}
Let us denote by $\widehat{\mathcal A}$ the subalgebra of $\widehat{M}$ generated by elements of the form $(\omega_{\widehat{J}J\Lambda_\varphi(x_1), J\Lambda_\varphi (y_1)}*id)(W)$ with $x_1$, $y_1$, in $\mathcal A_0$. Then, using \ref{propA2}, we get that $\widehat{\mathcal A}$ is invariant by $*$. Using now \ref{lemW}(ii), we get that $\widehat{\mathcal A}$ is invariant by $\widehat{R}$. Let us now denote by $\widehat{A}$ the norm closure of $\widehat{\mathcal A}$. It is clear that $\widehat{A}$ is a sub-$\bf C^*$-algebra of $\widehat{M}$, weakly dense in $\widehat{M}$, invariant by $\widehat{R}$.
%%%%thAhat
\subsection{Theorem}
\label{thAhat}
{\it (i) We have $\alpha(B)\subset M(\widehat{A})$ and $\widehat{\beta}(B)\subset M(\widehat{A})$.  Moreover, we have $\widehat{A}\subset \bf C^*$$(\widehat{\varphi}, \widehat{\varphi}\circ \widehat{R}, \widehat{\tau})$.
\newline
(ii) For any $X\in \widehat{A}\cap\gM_{\widehat{T}}$, $\widehat{T}(X)$ belongs to $\alpha(M(B))$; for any $Y\in \widehat{A}\cap\gM_{\widehat{R}\widehat{T}\widehat{R}}$, $\widehat{R}\widehat{T}\widehat{R}(Y)$ belongs to $\widehat{\beta}(M(B))$. 
\newline
(iii) for any $x\in\widehat{A}\cap\gN_{\widehat{T}}$, there exists $x_n$ in $\widehat{\mathcal A}$ (as defined in \ref{Ahat}), such that, for any $b\in B$, $\Lambda_{\widehat{T}}(x_n)b$ is norm converging to $\Lambda_{\widehat{T}}(x)b$.}
\begin{proof}
We have got in \ref{prophat} that, for any $x$, $y$ in $A\cap\gN_T\cap\gN_\varphi$, and $n\in B\cap\gN_\nu$, we have :
\[(\omega_{\widehat{J}J\Lambda_\varphi (x), J\Lambda_\varphi (y)}*id)(W)\alpha(n)=(\omega_{\widehat{J}J\Lambda_\varphi (x\alpha(n)), J\Lambda_\varphi (y)}*id)(W)\]
If now $x\in A\cap\gN_T\cap\gN_\varphi\cap\gN_{RTR}\cap\gN_{\varphi\circ R}$, we get that $R(x\alpha(n))=\beta(n)R(x)$ belongs to $A\cap\gN_T\cap\gN_\varphi$, and, therefore, that $x\alpha(n)$ belongs to $A\cap\gN_T\cap\gN_\varphi\cap\gN_{RTR}\cap\gN_{\varphi\circ R}$. Moreover, if $x$ belongs to $\mathcal A$ (as defined in \ref{propA}), and if $n$ is analytic with respect to $\sigma_t^\nu$ and $\gamma_t$, then $x\alpha(n)$ belongs to $\mathcal A$; therefore, if $X$ belong to $\widehat{\mathcal A}$, then $X\alpha(n)$ belongs to $\widehat{\mathcal A}$. Using now the invariance of $\widehat{\mathcal A}$ by $*$ and $\widehat{R}$, we get that $\alpha(n)X$, $\gamma(n)X$, $X\gamma(n)$ belong to $\widehat{\mathcal A}$. Then, by continuity, we get the first results of (i). Using now \ref{prop2hat}(v), we get that $\widehat{\mathcal A}\subset \bf C^*$$(\widehat{\varphi}, \widehat{\varphi}\circ \widehat{R}, \widehat{\tau})$, which, by density, finishes the proof of (i). 
\newline
Using \ref{lemproduct}, we get that, for any $X\in\widehat{\mathcal A}$, $X$ belongs to $\gN_{\widehat{T}}\cap\gN_{\widehat{\varphi}}$, and $\widehat{T}(X^*X)\in \alpha(M(B))$. Moreover, as $\widehat{\mathcal A}^+$ is norm dense in $\widehat{A}$, for any $x\in A\cap\gM_{\widehat{T}}^+$, it is possible to construct a sequence $x_n$ in $\widehat{\mathcal A}^+$, such that $x_n\leq x$ and $x_n$ is norm converging to $x$. Then, using (\cite{K2}, 3.5) to the restriction of $\widehat{T}$ to $\bf C^*$$(\widehat{\varphi}, \widehat{\varphi}\circ \widehat{R}, \widehat{\tau})$, we get that $\widehat{T}(x_n)$ is strictly converging to $\widehat{T}(x)$; therefore, we get that $\widehat{T}(x)$ belongs to $\alpha(M(B))$. Using the invariance of $\widehat{A}$ by $\widehat{R}$, we finish the proof of (ii). The proof of (iii) is obtained by similar arguments. 

\end{proof}

%%%%thAhat2
\subsection{Theorem}
\label{thAhat2}
{\it For any $x\in\widehat{A}$, $\widehat{\Gamma}(x)$ belongs to $\widehat{A}\underset{B}{_{\widehat{\beta}_{|B}}*_{\alpha_{|A}}}\widehat{A}$.}

\begin{proof}
Let $x$, $y$ in $\widehat{\mathcal A}$. Using \ref{proprho}(i) applied to $\widehat{\gG}$, and using the elements $a_n$ and $b_n$ defined such that $\Sigma_{\bf N} \Lambda_{\widehat{T}}(a_n)\Lambda_{\widehat{T}}(b_n)$ is an approximate unit of $\widehat{A}$, we get that :
\[(1\underset{N}{_{\widehat{\beta}}\otimes_\alpha}\Lambda_{\widehat{T}(x_1)}\Lambda_{\widehat{T}(x_2}^*)\widehat{\Gamma}(x)\rho^{\widehat{\beta}, \alpha}_{\widehat{J}\Lambda_{\widehat{\varphi}}(y)}
=
 \sum_N \rho^{\widehat{\beta}, \alpha}_{\widehat{J}\Lambda_{\widehat{\varphi}}(x_1\widehat{T}(x_2^* ya_n))}(\omega_{\Lambda_{\varphi}(x), \widehat{J}\Lambda_{\widehat{\varphi}}(b_n)}*id)(W)\]
As $x$ and $b_n$ belong to $\widehat{\mathcal A}$, we get that $\Lambda_{\widehat{\varphi}}(x)$ belongs to $\widehat{J}J\Lambda_\varphi(\mathcal A)$ and that $\widehat{J}\Lambda_{\varphi}(b_n)$ belongs to $J\Lambda_\varphi(\mathcal A)$; therefore, $(\omega_{\Lambda_{\varphi}(x), \widehat{J}\Lambda_{\widehat{\varphi}}(b_n)}*id)(W)$ belongs to $\widehat{\mathcal A}$, and we get, by continuity, that $(1\underset{N}{_{\widehat{\beta}}\otimes_\alpha}\Lambda_{\widehat{T}(x_1)}\Lambda_{\widehat{T}(x_2}^*)\widehat{\Gamma}(x)$ belongs to $\widehat{A}\underset{B}{_{\widehat{\beta}_{|B}}*_{\alpha_{|A}}}B(H)$. Which, by continuity, remains true for any $x$ in $\widehat{A}$. Using now similar arguments as in \ref{proprho}, we obtain the result. 
\end{proof}

%%%%thduality
\subsection{Theorem}
\label{thduality}
{\it The octuple $\widehat{\bfG}=(B, \widehat{A}, \alpha_{|B}, \widehat{\beta}_{|B}, \nu_{|B}, \widehat{T}_{|\widehat{A}}, \widehat{R}\widehat{T}\widehat{R}_{|\widehat{A}}, \widehat{\Gamma}_{|\widehat{A}})$ is a locally compact quantum groupoid, we shall call the dual of $\bfG$.}

\begin{proof} This is clear, using \ref{thAhat} and \ref{thAhat2}. \end{proof}

%%%propbidual
\subsection{Proposition}
\label{propbidual}
{\it The octuple $\widehat{\widehat{\bfG}}$, defined as the bidual of $\bfG$, is equal to $\bf G$. }

\begin{proof} Using \ref{Ahat}, we see that $\widehat{\widehat{A}}$ is the norm closure of the algebra $\widehat{\widehat{\mathcal A}}$ generated by all elements of the form 
$(id*\omega_{\widehat{J}\Lambda_{\widehat{\varphi}}(y), J\widehat{J}\Lambda_{\widehat{\varphi}}(x)})(W)$, where $x$, $y$ belong to $\widehat{\mathcal A}_0$, where $\widehat{\mathcal A}_0$ is, by \ref{propA}, the algebra of elements in $\widehat{A}\cap\gM_{\widehat{T}}\cap\gM_{\widehat{R}\widehat{T}\widehat{R}}\cap\gM_{\widehat{\varphi}}\cap\gM_{\widehat{\varphi}\circ\widehat{R}}$ which are analytic wit respect to $\sigma_t^{\widehat{\varphi}}$, $\sigma_t^{\widehat{\varphi}\circ\widehat{R}}$ and $\widehat{\tau_t}$. As $\widehat{\mathcal A}_0$ contains $\widehat{\mathcal A}$, the linear space $\Lambda_{\widehat{\varphi}}(\widehat{\mathcal A}_0)$ contains $\Lambda_{\widehat{\varphi}}(\widehat{\mathcal A})$, which is equal, using \ref{prophat}(i), to $\widehat{J}J\Lambda_\varphi(\mathcal A_0^2)$. Therefore, $\widehat{\widehat{\mathcal A}}$ contains all elements of the form $(id*\omega_{J\Lambda_{\varphi}(y), \Lambda_{\varphi}(x)})(W)$, where $x$, $y$ belong to $\mathcal A_0^2$. By continuity, we get that $\widehat{\widehat{A}}$ contains all elements of the form $(id*\omega_{J\Lambda_{\varphi}(y), \Lambda_{\varphi}(x)})(W)$, where $x$, $y$ belong to $\mathcal A_0$, from which we get that $A\subset \widehat{\widehat{A}}$. But, using \ref{thAhat}(iii), we get that $\widehat{\widehat{A}}$ is the norm closure of the set of all elements of the form $(id*\omega_{\widehat{J}\Lambda_{\widehat{\varphi}}(z_2\alpha(n_2)), J\widehat{J}\Lambda_{\widehat{\varphi}}(z_1\alpha(n_1))})(W)$, where $z_1$, $z_2$ belong to $\widehat{\mathcal A}$ and $n_1$, $n_2$ belong to $B$ and are analytc with respect to $\sigma_t^\nu$ and $\gamma_t$. Therefore, it is the norm closure of the set of all elements of the form $(id*\omega_{\widehat{J}\Lambda_{\widehat{\varphi}}(z_2), J\widehat{J}\Lambda_{\widehat{\varphi}}(z_1)})(W)$, where $z_1$, $z_2$ belong to $\widehat{\mathcal A}$, which is, as we have seen, the norm closure of the set of all elements of the form $(id*\omega_{J\Lambda_{\varphi}(y), \Lambda_{\varphi}(x)})(W)$, where $x$, $y$ belong to $\mathcal A_0^2$, which is $A$. \end{proof}

%%%%action
\section{Action of a locally compact quantum groupoid}
\label{action}

In this chapter, using \cite{Ti4}, we define an action of a locally compact quantum groupoid (\ref{defaction}) and obtain theorems about this notion (\ref{thA}, \ref{thS}), which are similar to some results about actions of measured quantum groupoids obtained in \cite{E4}. 

%%%%defaction
\subsection{\bf Definition}
\label{defaction}
Let $\bfG=(B, A, \alpha, \beta, \nu, T, T', \Gamma)$ be a locally compact quantum groupoid and let $C$ be a $\bf C^*$-algebra. An action of $\bfG$ on $C$ is a couple $(b,\ga)$ where :

(i) $b$ is an injective $*$-antihomomorphism from $B$ to $C$;

(ii) $\ga$ is an injective $*$-homomorphism from $C$ to $C\underset{B}{_b*_\alpha}A$;

(iii) $b$ and $\ga$ are such that, for any $n\in C$ :
\[\ga(b(n))=1\underset{B}{_b\otimes_\alpha}\beta(n)\]
(which allow us to define $\ga\underset{B}{_b*_\alpha}id$ from $C\underset{B}{_b*_\alpha}A$ into $C\underset{B}{_b*_\alpha}A\underset{B}{_\beta*_\alpha}A$) and such that :
\[(\ga\underset{B}{_b*_\alpha}id)\ga=(id\underset{B}{_b*_\alpha}\Gamma)\ga\]

%%%%implementation
\subsection{\bf Definition}
\label{implementation}
Let $\bfG=(B, A, \alpha, \beta, \nu, T, T', \Gamma)$ be a locally compact quantum groupoid and $(b,\ga)$ an action of $\bfG$ on a $\bf C^*$-algebra $C$. Let $\mathcal H$ be an Hilbert space, and $a$ an injective $*$-homomorphism of $C$ into $B(\mathcal H)$. An implementation of $\ga$ is a unitary $V$ in $\mathcal H\underset{B}{_{a\circ b}\otimes_\alpha}H_\varphi$ such that, for all $x\in A$ :
\[(a\underset{B}{_{a\circ b}\otimes id})\ga(x) =V(a(x)\underset{B}{_{a\circ b}\otimes_\alpha}1)V^*\]

%%%thimplemenattion
\subsection{\bf Theorem}
\label{thimplementation}
Let $\bfG=(B, A, \alpha, \beta, \nu, T, T', \Gamma)$ be a locally compact quantum groupoid and $(b,\ga)$ an action of $\bfG$ on a $\bf C^*$-algebra $C$, and $\widehat{\bfG}=(B, \widehat{A}, \alpha_{|B}, \widehat{\beta}_{|B}, \nu_{|B}, \widehat{T}_{|\widehat{A}}, \widehat{R}\widehat{T}\widehat{R}_{|\widehat{A}}, \widehat{\Gamma}_{|\widehat{A}})$ the dual of $\bfG$ constructed in \ref{thduality}. There exists an Hilbert space $\mathcal H$ and an injective 
an injective $*$-homomorphism $a$ of $C$ into $B(\mathcal H)$, such that $\ga$ has an implemantation. 

\begin{proof} Let $\pi$ be a faithful representation of $C$ on a HiIbert space $K$, and let us write $\mathcal K=K\underset{B}{_{\pi\circ b}\otimes_\alpha}H_{\nu\circ T}$, and $a(b)=1_K\underset{B}{_{b\circ\pi}\otimes_\alpha}\widehat{\alpha}(b)$, where $\widehat{\alpha}(b)=J_{\nu\circ\widehat{T}}\widehat{\beta}(b)J_{\nu\circ\widehat{T}}$. Let now $V=1_K\underset{B}{_{b\otimes\pi}\otimes_\alpha}(\sigma W^o\sigma)^*$; then $V$ is an implementation of $\ga$, in the sense of \ref{implementation} \end{proof}

%%%exemple
\subsection{\bf Example}
\label{exaction}
Let $(A,\Gamma, T, T')$ be a locally compact quantum group (\ref{trivial}), and $\ga$ an action of this locally compact quantum on a $\bf C^*$-algebra $C$, i.e. $\ga$ is an injective $*$-homomorphism from $C$ into $M(C\otimes A))$ such that $(\ga\otimes id)\ga=(id\otimes\Gamma)\ga$. (\cite{BS}, 0.2)

%%%%crossedproduct
\subsection{\bf Definition}
\label{crossedproduct}
Let $\bfG=(B, A, \alpha, \beta, \nu, T, T', \Gamma)$ be a locally compact quantum groupoid, and let $(b,\ga)$ an action of $\bfG$ on a $\bf C^*$-algebra $C$; the crossed product of this action is the $\bf C^*$-algebra $C\rtimes_\mathfrak a\ {\bfG}$ generated by the elements $\ga (c)(1\underset{B}{_b\otimes_\alpha} \widehat{J}x\widehat{J})$, with $c\in C$ and $x\in\widehat{A}$. 
%%%exemple
\subsection{\bf Example}
\label{excorossedproduct}
If $(A,\Gamma, T, T')$ be a locally compact quantum group (\ref{trivial}), and $\ga$ an action of this locally compact quantum group on a $\bf C^*$-algebra $C$, the crossed-product is then generated by $\ga(c)(1\otimes \widehat{J}x\widehat{J})$ (\cite{BS}, 7.1).

%%%thA
\subsection{\bf Theorem}
\label{thA}
Let $\bfG=(B, A, \alpha, \beta, \nu, T, T', \Gamma)$ be a locally compact quantum groupoid, and let $(b,\ga)$ an action of $\bfG$ on a $\bf C^*$-algebra $C$; then, the $\bf C^*$-algebra generated by all elements of the form $\ga(c)(1\underset{B}{_b \otimes_\alpha}x)$, where $c\in C$ and $x\in B(H_{\nu\circ T})$ is equal to 
$C\underset{B}{_b *_\alpha}B(H_{\nu\circ T})$. 

\begin{proof}
As $\ga(C)\subset C\underset{B}{_b \otimes_\alpha}A$, we get that the $\bf C^*$-algebra generated by all elements of the form $\ga(c)(1\underset{B}{_b \otimes_\alpha}x)$, where $c\in C$ and $x\in B(H_{\nu\circ T})$ is included in $C\underset{B}{_b *_\alpha}B(H_{\nu\circ T})$. Let us define $\underline{C}$ as the $\bf C^*$-algebra generated by all elements of the form $(id\underset{B}{_b *_\alpha}\omega)\ga(c)$, where $c\in C$ and $\omega$ is a state on $B(H_{\nu\circ T})$. Then, we get that the $\bf C^*$-algebra generated by all elements of the form $\ga(c)(1\underset{B}{_b \otimes_\alpha}x)$, where $c\in C$ and $x\in B(H_{\nu\circ T})$ is equal to $\underline{C}\underset{B}{_b \otimes_\alpha}B(H_{\nu\circ T})$. 

Then , we get that $\ga(C)\underset{B}{_b *_\alpha}B(H_{\nu\circ T})=(\ga\underset{B}{_b \otimes_\alpha}id)(\underline{C}\underset{B}{_b \otimes_\alpha}B(H_{\nu\circ T}))=\ga(\underline{C})\underset{B}{_b *_\alpha}B(H_{\nu\circ T})$, from which we get that $\underline{C}=C$ and the result. 
\end{proof}

%%%thS
\subsection{\bf Theorem}
\label{thS}
Let $\bfG=(B, A, \alpha, \beta, \nu, T, T', \Gamma)$ be a locally compact quantum groupoid, and let $(b,\ga)$ an action of $\bfG$ on a $\bf C^*$-algebra $C$; then, we have :
\[\ga(C)=\{X\in C\underset{B}{_b*_\alpha}A, (\ga\underset{B}{_b*_\alpha}id)(X)=(id\underset{B}{_b*_\alpha}\Gamma)(X)\}\]
\begin{proof}
Using \ref{thimplementation}, we get that there exists an Hilbert space $K$, a faithful representation $\pi$ of $C$ on $K$ and a unitary $V$ on $K\underset{B}{_{\pi\circ b}\otimes_\alpha}H_{\nu\circ T}$ which implements the action $\ga$ on $\pi(C)$. Let us consider now the measured quantum groupoid $\mathfrak{G}=(N, M, \underline{\alpha}, \underline{\beta}, \underline{\Gamma}, \underline{T}, \underline{T'}, \underline{\nu})$ constructed from $\bfG$ (\ref{thMQG}). Then, we get that $V$ implements an action $\underline{\ga}$ of $\mathfrak {G}$ on the von Neumann algebra $D$ generated by $\pi(C)$. Using (\cite{E4}11.8) applied to $D$, we get that :
\[\underline{\ga}(D)=\{Y\in D\underset{N}{_b*_\alpha}M, \underline{\ga}\underset{N}{_b*_\alpha}id)(Y)=(id\underset{N}{_b*_\alpha}\underline{\Gamma})(Y)\}\]
from which we get that :
\[\{X\in C\underset{B}{_b*_\alpha}A, (\ga\underset{B}{_b*_\alpha}id)(X)=(id\underset{B}{_b*_\alpha}\Gamma)(X)\}\subset \underline{\ga}(D)\]
More precisely, there exists a $\bf C^*$-algebra $\widetilde{C}$, with $C\subset \widetilde{C}\subset D$, such that :
\[\{X\in C\underset{B}{_b*_\alpha}A, (\ga\underset{B}{_b*_\alpha}id)(X)=(id\underset{B}{_b*_\alpha}\Gamma)(X)\}=\underline{\ga}(\widetilde{C})\]
Then, it is easy to see that $\underline{\ga}_{|{\widetilde{C}}}$ is an action of $\bfG$ on the $\bf C^*$-algebra $\widetilde{C}$, such that $\underline{\ga}(\widetilde{C})\subset C\underset{B}{_b*_\alpha}A$. Applying now \ref{thA} to $\underline{\ga}_{|{\widetilde{C}}}$, we get that $\widetilde{C}\underset{B}{_b*_\alpha}B(H_{\nu\circ T})=C\underset{B}{_b*_\alpha}B(H_{\nu\circ T})$, from which we get that $\widetilde{C}=C$, which gives the result.

\end{proof}

%%%%%dualaction

\subsection{\bf Definition}
\label{defdualaction}
Let $\bfG=(B, A, \alpha, \beta, \nu, T, T', \Gamma)$ be a locally compact quantum group (\ref{trivial}), and $(b, \ga)$ an action of $\bfG$ on a $\bf C^*$-algebra $C$, as defined in \ref{defaction}, and $C\rtimes_\mathfrak a\ {\bfG}$  the crossed-product of this action (\ref{crossedproduct}). Let us define the dual action $(1\underset{B}{_b\otimes_\alpha}\widehat{\alpha}, \tilde{\ga})$ of $\widehat{\bfG}^c$ on $C\rtimes_\mathfrak a\ {\bfG}$ defined by, for any $c\in C$ and $X\in \widehat{A}$:
\[\tilde{\ga}(\ga(c)(1\underset{B}{_b\otimes_\alpha}\widehat{J}X\widehat{J}))=(\ga(c)\underset{B^o}{_{\widehat{\alpha}}\otimes_\beta}1)(1\underset{B}{_b\otimes_\alpha}\widehat{J}\widehat{\Gamma}(X)\widehat{J})\]
Then, for any Hilbert space $\mathcal H$ such that there exists an injective $*$-homomorphism $a$ of $C\rtimes_\mathfrak a\ {\bfG}$ into $B(\mathcal H)$, we verify that, for any $Y$ in $C\rtimes_\mathfrak a\ {\bfG}$, we have :
\[\tilde{\ga}(Y)=(1\underset{B}{_b\otimes_\alpha}(J\underset{B^o}{_{\alpha}\otimes_{\widehat{\beta}}}\widehat{J})W(\widehat{J}\underset{B}{_{\widehat{\beta}}\otimes_\alpha}J))Y(1\underset{B}{_b\otimes_\alpha}(\widehat{J}\underset{B}{_{\widehat{\beta}}\otimes_\alpha}J)W^*(J\underset{B^o}{_{\alpha}\otimes_{\widehat{\beta}}}\widehat{J}))\]
which proves that $1\underset{B}{_b\otimes_\alpha}(J\underset{B^o}{_{\alpha}\otimes_{\widehat{\beta}}}\widehat{J})W(\widehat{J}\underset{B}{_{\widehat{\beta}}\otimes_\alpha}J))Y(1\underset{B}{_b\otimes_\alpha}(\widehat{J}\underset{B}{_{\widehat{\beta}}\otimes_\alpha}J)W^*(J\underset{B^o}{_{\alpha}\otimes_{\widehat{\beta}}}\widehat{J})$ is an implementation of $\tilde{\ga}$.

%%%%%%%exemples
\section{Examples}
\label{Examples}

%%%finitebasis
 \subsection{Locally compact sub-quantum groupoid with finite dimensional basis}
 \label{finite}
 
When the basis is finite dimensional, the construction of locally compact quantum groupoids had already been studied by De Commer (\cite{DC1}, chapter 11) and Baaj and Crespo (\cite {BC} 2.2, \cite{Cr}).

%%%%%%fieldsqg
\subsection{Continuous field of locally compact quantum groups \cite{B2}}
\label{Bl}
In this chapter, we define a notion of continuous field of locally compact quantum groups (\ref{deffieldqg}), which was underlying in \cite{B2}. We show that these are exactly locally compact quantum groupoids with central basis, and that the dual object is of the same kind (\ref{bicentral}). We finish by recalling concrete examples (\ref{SU}, \ref{ax+b}, \ref{Emu}) given by Blanchard, which are examples of locally compact quantum groupoids.

%%%%%%deffieldqg
\subsubsection{\bf{Definition}}
\label{deffieldqg}
A 6-uple $(X, \alpha, A, \Gamma^x, \varphi^x, \psi^x)$ will be called a \emph{continuous field of locally compact quantum groups} if :
\newline
(i) $(X, \alpha, A, \varphi^x)$ and $(X, \alpha, A, \psi^x)$ are measurable continuous fields of $\bf C^*$-algebras; 
\newline
(ii) for any $x\in X$, there exists a simplifiable coproduct $\Gamma^x$ from $A^x$ to $M(A^x\otimes_m A^x)$ such that $(A^x, \Gamma^x, \varphi^x, \psi^x)$ is a locally compact quantum group, in the sense of [KV1]. 
\newline
Simplifiable means that the closed linear set generated by $\Gamma^x(A^x)(A^x\otimes_m 1)$ is equal to $A^x\otimes_m A^x$ (resp. $\Gamma^x(A^x)(1\otimes_m A^x)$).

%%%%deffieldqg
\subsubsection{\bf{Definition }}(\cite{L}, 17.1.3)
\label{deffieldqg}
Let $(X, \nu)$ be a $\sigma$-finite standard measure space; let us take $\{M^x, x\in X\}$ a measurable field of von Neumann algebras over $(X, \nu)$ and $\{\varphi^x, x\in X\}$ (resp. $\{\psi^x\}$) a measurable field of normal semi-finite faithful weights on $\{M^x\}$ (\cite{T}, 4.4). Moreover, let us suppose that :
\newline
(i) there exits a measurable field of injective $*$-homomorphisms $\Gamma^x$ from $M^x$ into $M^x\otimes M^x$ (which is also a measurable field of von Neumann algebras, on the measurable field of Hilbert spaces $H_{\varphi^x}\otimes H_{\varphi^x}$). 
\newline
(ii) for almost all $x\in X$, $\textbf{G}^x=(M^x, \Gamma^x, \varphi^x, \psi^x)$ is a locally compact quantum group (in the von Neumann sense (\cite{KV2}). 
\newline
In that situation, we shall say that $(M^x, \Gamma^x, \varphi^x, \psi^x, x\in X)$ is a \emph{measurable field of locally compact quantum groups} over $(X, \nu)$. 

%%%%mqgfield
\subsubsection{\bf{Theorem} }(\cite{L}, 17.1.3, \cite{E5}, 8.2)
\label{mqgfield}
{\it Let $\textbf{G}^x=(M^x, \Gamma^x, \varphi^x, \psi^x, x\in X)$ be a measurable field of locally compact quantum groups over $(X, \nu)$. Let us define :
\newline
(i) $M$ as the von Neumann algebra made of decomposable operators $\int_X^{\oplus}M^xd\nu(x)$, and $\alpha$ the $*$-isomorphism which sends $L^\infty(X, \nu)$ into the algebra of diagonalizable operators, which is included in $Z(M)$.
\newline
(ii) $\Phi$ (resp. $\Psi$) as the direct integral $\int_X^{\oplus}\varphi^x d\nu(x)$ (resp. $\int_X^{\oplus}\psi^xd\nu(x)$). Then, the Hilbert space $H_\Phi$ is equal to the direct integral $\int_X^{\oplus}H_{\varphi^x} d\nu(x)$, the relative tensor product $H_\Phi\underset{\nu}{_\alpha\otimes_\alpha}H_\Phi$ is equal to the direct integral $\int_X^{\oplus}(H_{\varphi^x}\otimes H_{\varphi^x})d\nu(x)$, and the product $M\underset{N}{_\alpha*_\alpha}M$ is equal to the direct integral $\int_X^{\oplus}(M^x\otimes M^x)d\nu(x)$.
\newline
(iii) $\Gamma$ as the decomposable $*$-homomorphism $\int_X^{\oplus}\Gamma^xd\nu(x)$, which sends $M$ into $M\underset{N}{_\alpha*_\alpha}M$.
\newline
(iv) $T$ (resp. $T'$) as an operator-valued weight from $M$ into $\alpha (L^\infty(X, \nu))$ defined the following way : $a\in M^+$ represented by the field $\{a^x\}$ belongs to $\gM_T^+$ if, for almost all $x\in X$, $a^x$ belongs to $\gM_{\varphi^x}$ (resp. $\gM_{\psi^x}$), and the function $x\mapsto \varphi^x(a^x)$ (resp. $x\mapsto \psi^x(a^x)$) is essentiallly bounded; then $T(a)$ (resp. $T'(a)$) is defined as the image under $\alpha$ of this function. 
\newline
Then, $(L^\infty(X, \nu), M, \alpha, \alpha, \Gamma, T, T', \nu)$ is a measured quantum groupoid, we shall denote by $\int_X^{\oplus}\textbf{G}^xd\nu(x)$.}

%%%%propmqgfield
\subsubsection{\bf{Proposition}}(\cite{E5}, 8.3)
\label{propmqgfield}
{\it Let $(X,\nu)$ be a $\sigma$-finite standard measure space, and $\{\textbf{G}^x, x\in X\}$ a measurable field of locally compact quantum groups, as defined in \ref{deffieldqg}; let $\int_X^{\oplus}\textbf{G}^xd\nu(x)$ be the measured quantum groupoid constructed in \ref{mqgfield}; then :
\newline
(i) we have $\alpha=\beta=\hat{\beta}$;
\newline
(ii) the pseudo-multiplicative unitary of the the measured quantum groupoid is a unitary on $H_\Phi\underset{\nu}{_\alpha\otimes_\alpha}H_\Phi$, which is equal to the decomposable operator $\int_X^{\oplus}W^xd\nu(x)$, where $W^x$ is the multiplicative unitary associated to the locally compact quantum group $\textbf{G}^x$. 
\newline
(iii) we have :}
\[\widehat{\int_X^{\oplus}\textbf{G}^xd\nu(x)}=\int_X^{\oplus}\widehat{\textbf{G}^x}d\nu(x)\]

%%%%propbicentral
\subsubsection{\bf{Proposition}}(\cite{E5}, 8.4)
\label{propbicentral}
{\it Let $\gG=(N, M, \alpha, \beta, \Gamma, T, T', \nu)$ be a measured quantum groupoid and $\widehat{\gG}=(N, \widehat{M}, \alpha, \hat{\beta}, \widehat{\Gamma}, \hat{T}, \widehat{R}\hat{T}\widehat{R}, \nu)$ its dual measured quantum groupoid. Then, are equivalent:
\newline
(i) $\alpha (N)\subset Z(M)\cap Z(\widehat{M})$.
\newline
(ii) $\alpha=\beta=\hat{\beta}$. }

%%%%%%thbicentral
\subsubsection{\bf{Theorem}}(\cite{E5}, 8.5)
\label{bicentral}
{\it Let $\gG=(N, M, \alpha, \beta, \Gamma, T, T', \nu)$ be a measured quantum groupoid and $\widehat{\gG}=(N, \hat{M}, \alpha, \hat{\beta}, \widehat{\Gamma}, \hat{T}, \widehat{R}\hat{T}\widehat{R}, \nu)$ its dual measured quantum groupoid; let $W$ and $\widehat{W}$ be the pseudo-multiplicative unitaries associated, and $\Phi=\nu\circ\alpha^{-1}\circ T$ (resp. $\widehat{\Phi}=\nu\circ\alpha^{-1}\circ\hat{T}$); let us suppose that $\alpha (N)$ is central in both $M$ and $\widehat{M}$; let $X$ be the spectrum of ${\bf C}^*$algebra ${\bf C}^*(\nu)$, we shall therefore identify with $C_0(X)$; for any $x\in X$, let $C_x(X)$ be the subalgebra of $C_0(X)$ made of functions which vanish at $x$; let $A_n(W)$ be the sub-${\bf C}^*$-algebra of $M$ introduced in \ref{AnW}, which is, thanks to $\alpha_{|C_0(X)}$, a continuous field over $X$ of ${\bf C}^*$-algebras ([T], 4.11); let $\varphi^x$ be the desintegration of $\Phi_{|A_n(W)}$ over $X$; $\varphi^x$ is a lower semi-continuous weight  on $A_n(W)$, faithful when considered on $A_n(W)/\alpha(C_x(X))A_n(W)$, and the representation $\pi_{\varphi^x}$ form a continuous field of faithful representation of $A_n(W)$. Then :
\newline
(i) the Hilbert space $H_\Phi\underset{\nu}{_\alpha\otimes_\alpha}H_\Phi$ is equal to $\int_X^{\oplus}H_{\varphi^x}\otimes H_{\varphi^x}d\nu(x)$. 
\newline
(ii) the von Neumann algebra $M\underset{N}{_\alpha*_\alpha}M$ is equal to :
\[\int_X^{\oplus}\pi_{\varphi^x}(A_n(W)/\alpha(C_x(X))A_n(W))"\otimes\pi_{\varphi^x}(A_n(W)/\alpha(C_x(X))A_n(W))"d\nu(x)\]
\newline
(ii) the coproduct $\Gamma_{|A_n(W)}$ can be desintegrated in $\Gamma_{|A_n(W)}=\int_X^{\oplus}\Gamma^xd\nu(x)$, where $\Gamma^x$ is a continuous field of coassociative coproducts on $A_n(W)/\alpha(C_x(X))A_n(W)$. 
\newline
(iii) $R^x$ is a anti-$*$-automorphism of $A_n(W)/\alpha(C_x(X))A_n(W)$, and, for all $x\in X$, $(A_n(W)/\alpha(C_x(X))A_n(W), \Gamma^x, \varphi^x, \varphi^x\circ R^x)$ is a locally compact quantum group (in the ${\bf C}^*$-sense), we shall denote ${\bf G}^x$. We shall denote also ${\bf G}^x$ its von Neumann version. 
\newline
(iv) we have, with the notations of \ref{mqgfield}, $\gG=\int_X^{\oplus}{\bf G}^x d\nu(x)$. }

%%%cfield
\subsubsection{\bf{Theorem}}(\cite{E5}, 8.6)
\label{cfield}
{\it Let $(X,\nu)$ be a $\sigma$-finite standard measure space, $\textbf{G}^x$ be a measurable field of locally compact quantum groups over $(X, \nu)$, ad defined in \ref{deffieldqg}, and $\int_X^{\oplus}\textbf{G}^xd\nu(x)$ be the measured quantum groupoid constructed in \ref{mqgfield}. Then :
\newline
(i) there exists a locally compact set $\tilde{X}$, and a positive Radon measure $\tilde{\nu}$ on $\tilde{X}$, such that $L^\infty(X,\nu)$ and $L^\infty(\tilde{X}, \tilde{\nu})$ are isomorphic, and such that this isomorphism sends $\nu$ on $\tilde{\nu}$. 
\newline
(ii) there exists a continuous field $(A^x)_{x\in\tilde{X}}$ of ${\bf C}^*$-algebras, and a continuous field of coassociative coproducts $\tilde{\Gamma}^x : A^x\rightarrow A^x\otimes^mA^x$;
\newline
(iii) there exists left-invariant (resp. right-invariant) weights $\tilde{\varphi}^x$ (resp. $\tilde{\psi}^x$), such that $(A^x, \tilde{\Gamma}^x, \tilde{\varphi}^x, \tilde{\psi}^x)$ is a locally compact quantum group $\tilde{\textbf{G}}^x$ (in the ${\bf C}^*$ sense). 
\newline
(iv) we have : $\int_X^{\oplus}\textbf{G}^xd\nu(x)=\int_{\tilde{X}}^{\oplus}\tilde{\textbf{G}}^xd\tilde{\nu}(x)$.}

%%%%lcquantumgroupoid
\subsubsection{\bf{Theorem}}
\label{lcquantumgroupoid}
{\it With the notations of \ref{cfield}, let us define $A=\int_{\tilde{X}}^{\oplus}A^xd\tilde{\nu}(x)$. Then, $(C_0(\tilde{X}), A, \alpha_{|C_0(\tilde{X})}, \alpha_{|C_0(\tilde{X})},, \tilde{\nu}, \Gamma_{|A})$ is a locally compact quantum groupoid. }

\begin{proof}
Let $a=\int_{\tilde{X}}^{\oplus}a^xd\tilde{\nu}(x)$ in $A$, and $b=\int_{\tilde{X}}^{\oplus}b^xd\tilde{\nu}(x)$ in $A$, with $b^x\in\gN_{T^x}\cap\gN_{\varphi^x}$; then, we get that $\rho^{\alpha_{|A}, \alpha_{|A}}_{J_\varphi\Lambda_\varphi(b)}\Gamma(a)=\int_{\tilde{X}}^{\oplus}\rho^{\alpha_x, \alpha_x}_{J_{\varphi_x}(b^x)}\Gamma_x(a^x)d\tilde{\nu}(x)$.\end{proof}

%%%%%SU
\subsubsection{\bf{Example}}
\label{SU}
As in (\cite{B1}, 7.1), let us consider the ${\bf C}^*$-algebra $A$ whose generators $\alpha$, $\gamma$ and $f$ verify :
\newline
(i) $f$ commutes with $\alpha$ and $\gamma$;
\newline
(ii) the spectrum of $f$ is $[0,1]$;
\newline
(iii) the matrix 
$\left(\begin{array}{cc}
\alpha&-f\gamma\\
\gamma & \alpha^*
\end{array}
\right)$ is unitary in $M_2(A)$. 
Then, using the sub ${\bf C}^*$-algebra generated by $f$, $A$ is a $C([0,1])$-algebra; let us consider now $A$ as a $C_0(]0,1])$-algebra.  Then, Blanchard had proved (\cite{B2} 7.1) that $A$ is a continuous field over $]0,1]$ of ${\bf C}^*$-algebras, and that, for all $q\in]0,1]$, we have $A^q=SU_q(2)$, where the $SU_q(2)$ are the compact quantum groups constructed by Woronowicz and $A^1=C(SU(2))$. 
\newline
Moreover, using the coproducts $\Gamma^q$ defined by Woronowicz as 
\[\Gamma^q(\alpha)=\alpha\otimes\alpha-q\gamma^*\otimes\gamma\]
\[\Gamma^q(\gamma)=\gamma\otimes\alpha+\alpha^*\otimes\gamma\]
and the (left and right-invariant) Haar state $\varphi^q$, which verifies : 
\newline
$\varphi^q(\alpha^k\gamma^{*m}\gamma^n)=0$, for all $k\geq 0$, and $m\not=n$, 
\newline
$\varphi^q(\alpha^{*|k|}\gamma^{*m}\gamma^n)=0$, for all $k<0$, and $m\not=n$, 
\newline
and $\varphi^q((\gamma^*\gamma)^m)=\frac{1-q^2}{1-q^{2m+2}}$, 
\newline
we obtain this way a continuous field of compact quantum groups (see \cite{B2}, 6.6 for a definition); this leads to put on $A$ a structure of locally compact quantum groupoid (of compact type, in the sense of \cite{E2}, because $1\in A$). 
\newline
This structure is given by a coproduct $\Gamma$ which is $C_0(]0,1])$-linear from $A$ to $A\underset{C_0(]0,1])}{\otimes^m}A$, and given by :
\[\Gamma(\alpha)=\alpha\underset{C_0(]0,1])}{\otimes^m}\alpha-f\gamma^*\underset{C_0(]0,1])}{\otimes^m}\gamma\]
\[\Gamma(\gamma)=\gamma\underset{C_0(]0,1])}{\otimes^m}\alpha+\alpha^*\underset{C_0(]0,1])}{\otimes^m}\gamma\]
and by a conditional expectation $E$ from $A$ on $M(C_0(]0,1]))$ given by :
\newline
$E(\alpha^k\gamma^{*m}\gamma^n)=0$, for all $k\geq 0$, and $m\not=n$
\newline
$E(\alpha^{*|k|}\gamma^{*m}\gamma^n)=0$, for all $k<0$, and $m\not=n$
\newline
$E((\gamma^*\gamma)^m)$ is the bounded function $x\mapsto\frac{1-q^2}{1-q^{2m+2}}$. 
Then $E$ is both left and right-invariant with respect to $\Gamma$. This example gives results at the level of ${\bf C}^*$-algebras, which are more precise than theorem \ref{mqgfield}. 

%%%%%ax+b
\subsubsection{\bf{Example}}
\label{ax+b}
One can find in \cite{B2} another example of a continuous field of locally compact quantum group. Namely, in (\cite{B2}, 7.2), Blanchard constructs a ${\bf C}^*$-algebra $A$ which is a continuous field of ${\bf C}^*$-algebras over $X$, where $X$ is a compact included in $]0, 1]$, with $1\in X$. For any $q\in X$, $q\not=1$, we have $A^q=SU_q(2)$, and $A^1={\bf C}^*_r(G)$, where $G$ is the "$ax+b$" group. (\cite{B2}, 7.6). 
\newline
Moreover, he constructs a coproduct (denoted $\delta$) (\cite{B2} 7.7(c)), and "the system of Haar weights" $\Phi$ ([B2] 7.2.3), which bear left-invariant-like properties (end of remark after \cite{B2} 7.2.3). 
\newline
Finally, he constructs a unitary $U$ in $\mathcal L(\mathcal E_\Phi)$ (\cite{B2} 7.10), with which it is possible to construct a co-inverse $R$ of $(A, \delta)$, which leads to the fact that $\Phi\circ R$ is right-invariant. 
\newline
Clearly, the fact that we are here dealing with non-compact locally compact quantum groups made the results more problematic at the level of ${\bf C}^*$-algebra; at the level of von Neumann algebra, \ref{mqgfield} allow us to construct an example of measured quantum groupoid from these data. 

%%%%%Emu
\subsubsection{\bf{Example}}
\label{Emu}
Let us finish by quoting a last example given by Blanchard in (\cite{B2}, 7.4): for $X$ compact in $[1, \infty[$, with $1\in X$, he constructs a ${\bf C}^*$-algebra which is a continuous field over $X$ of ${\bf C}^*$ algebras, whose fibers, for $\mu\in X$, are $A^\mu=E_\mu(2)$, with a coproduct $\delta$ and a continuous field of weights $\Phi$, which is left-invariant. As in \ref{ax+b}, he then constructs a unitary $U$ on $\mathcal L(\mathcal E_\Phi)$, which will lead to a co-inverse, and, therefore, to a right-invariant ${\bf C}^*$-weight. 

%%%%%sum
\subsubsection{\bf{Example}} (\cite{L}, 17.1)
\label{sum}
Let us return to \ref{deffieldqg}; let $I$ be a (discrete) set, and, for all $i$ in $I$, let $\textbf{G}_i=(M_i, \Gamma_i, \varphi_i, \psi_i)$ be a locally compact quantum group; then the product $\Pi_i \textbf{G}_i$ is a measured field of locally compact quantum groups, and can be given a natural structure of measured quantum groupoid, described in (\cite{L}, 17.1).  

%%%%C2
\subsubsection{\bf{Example}} (\cite{DC2}3.17)(\cite{E5}, 9)
\label{C2}
Let $({\bf C}^2, M, \alpha, \beta, \Gamma, T, RTR, \nu)$ be a measured quantum groupoid; let $(e_1, e_2)$ be the canonical basis of $\bf C^2$ et let us suppose that $\nu(e_1)=\nu(e_2)=1/2$; let us define $M_{i,j}=M\alpha(e_i)\beta(e_j)$; then, the fiber product $M\underline{\nu}{_\beta\*_\alpha}M$ can be idetified with the reduces von Neumann algebra $(M\otimes M)_{\beta(e_1)\otimes\alpha(e_2)+\beta(e_2)\alpha(e_1)}$, and $\Gamma$ can be identified with an injective $*$-homomorphism from $M$ to $M\otimes M$, such that $(\Gamma\otimes id)\Gamma=(id\otimes\Gamma)\Gamma$ and :
\[\Gamma(1)=\beta(e_1)\otimes\alpha(e_1)+\beta(e_2)\otimes\alpha(e_2)\]
\[\Gamma(\alpha(e_i)\beta(e_j)=\alpha(e_i)\beta(e_1)\otimes\alpha(e_1)\beta(e_j)+\alpha(e_i)\beta(e_2)\otimes\alpha(e_2)\beta(e_j)\]
Let $M_{i,j}=M\alpha(e_i)\beta(e_j)$, we have $M_{1,1}\neq\{0\}$, $M_{2,2}\neq\{0\}$, and 
\[\Gamma(M_{i,j})\subset (M_{i,1}\otimes M_{1,j})\oplus (M_{i,2}\otimes M_{2,j})\]
There exist normal semi-finite faithful weights $\varphi_{i,j}$ on $M_{i,j}$ such that, for any $X\in\gM_T^+$, $X=\sum_{i,j}x_{i,j}$, with $x_{i,j}\in M_{i,j}^+$, we get that $T(X)$ is the image under $\alpha$ of 
\[(\varphi_{1,1}(x_{1,1}+\varphi_{1,2}(x_{1,2})e_1+(\varphi_{2,1}(x_{2,1})+\varphi_{2,2}(x_{2,2})e_2\]
For $x_{i,j}\in M_{i,j}$, and $k=1,2$, let us define $\Gamma_{i,j}^k(x_{i,j})=\Gamma(x_{i,j})[\alpha((e_i)\beta(e_k)\otimes\alpha(e_k)\beta(e_j)]$
Then,  $\textbf{G}^1=(M_{1,1}, \Gamma_{1,1}^1, \varphi_{1,1}, \psi_{1,1})$ and $\textbf{G}^2=(M_{2,2}, \Gamma_{2,2}^2, \varphi_{2,2}, \psi_{2,2})$ are two locally compact quantum groups; 
\newline
if $\alpha({\bf C}^2)\subset Z(\widehat{M})$, then $\beta=\alpha$ and $\gG=\textbf{G}^1\oplus\textbf{G}^2$ (i.e. $M_{1,2}=M_{2,1}=\{0\}$).
\newline
If $\alpha(e_1)\notin Z(\widehat{M})$, then $M_{1,2}\neq \{0\}$, $M_{2,1}=R(M_{1,2})\neq \{0\}$, $\Gamma_{1,2}^2 : M_{1,2}\rightarrow M_{1,2}\otimes M_{2,2}$ is a right action of $\textbf{G}^2$ on $M_{1,2}$, $\Gamma_{1,2}^1: M_{1,2}\rightarrow M_{1,1}\otimes M_{1,2}$ is a left action of $\textbf{G}^1$ on $M_{1,2}$, $\Gamma_{2,1}^1: M_{2,1}\rightarrow M_{2,1}\otimes M_{1,1}$ verify $\Gamma_{2,1}^1=\varsigma(R^1\otimes R)\Gamma_{1,2}^1R$, and $\Gamma_{2,1}^2: M_{2,1}\rightarrow M_{2,2}\otimes M_{2,1}$ verify $\Gamma_{2,1}^2=\varsigma(R\otimes R^2)\Gamma_{1,2}^2R$. Moreover, these actions are ergodic and integrable.

%%abelian
\subsection{Abelian measured quantum groupoid}
\label{abelian}
We consider now the case of an "abelian" measured quantum groupoid (i.e. a measured quantum groupoid $\gG=(N, M, \alpha, \beta, \Gamma, T, T', \nu)$ where the underlying von Neuman algebra $M$ is abelian); then we prove that it is possible to put on the spectrum of the ${\bf C}^*$-algebra $A_n(W)$ a structure of a locally compact groupoid, whose basis is the spectrum of ${\bf C}^*(\nu)$ ([E5], 7.1). Starting from a measured groupoid equipped with a left-invariant Haar system, we recover Ramsay's theorem which says that this groupoid is measure-equivalent to a locally compact one (\ref{ramsay}).

%%%%%abelian

\subsubsection{\bf{Theorem}}
\label{thabelian}
{\it Let $\bfG=(B, A, \alpha, \beta, \nu, T, T', \Gamma)$ be a locally compact quantum groupoid, such that the $\bf C^*$-algebra $A$ (and, therefore, the $\bf C^*$-algebra $B$) is abelian. Let $\mathcal G$ be the spectrum of $A$, and $\mathcal G^{(0)}$ be the spectrum of $B$. Then :

(i) there are two continuous open applications $r$, $s$ from $\mathcal G$ to $\mathcal G^{(0)}$ such that, for all $f\in \bf C$$_0(\mathcal G^{(0)})$, we have $\alpha(f)=f\circ r$ and $\beta(f)=f\circ s$. 

(ii) there is a partially defined multiplication on $\mathcal G$, which gives to $\mathcal G$ a structure of locally compact groupoid, with $\mathcal G^{(0)}$ as set of units. Then, for any $f\in \bf C$$_0(\mathcal G)$, $\Gamma(f)\in\bf C$$_b(\mathcal G^{(2)})$.

(iii) for all $f\in\mathcal K(\mathcal G)$ and $u\in\mathcal G^{(0)}$, the application $f\rightarrow \alpha^{-1}(T(f)))(u)$ defines a positive Radon measure $\lambda^u$ on $\mathcal G$, whose support is $\mathcal G^u=\{x\in\mathcal G, r(x)=u\}$. The measures $(\lambda^u)_{u\in\mathcal G^{(0)}}$ are a Haar system on $\mathcal G$

(iv) the trace $\nu$ on $\bf C$$_0(\mathcal G^{(0)})$ is a measure on $\mathcal G^{(0)}$, which is quasi-invariant with respect to this Haar system. 

(v) we have $\bfG=(\bf C$$_0(\mathcal G^{(0)}), \bf C$$_0(\mathcal G), r_{\mathcal G}, s_{\mathcal G}, \Gamma_{\mathcal G}, (\lambda^u)_{u\in\mathcal G^{(0)}}, (\lambda_u)_{u\in\mathcal G^{(0)}}, \nu)$, where $r_{\mathcal G}(f)=f\circ r$, $s_{\mathcal G}=f\circ s$, $\Gamma_{\mathcal G}(f)(x_1, x_2)=f(x_1x_2)$, for any $f$ in $\bf C$$_0(\mathcal G)$, $(x_1, x_2)\in\mathcal G^{2)}$, and $\lambda_u$ is the image of $\lambda^u$ by $x\rightarrow x^{-1}$.}

\begin{proof} Completely similar to (\cite{E5}, 7.1). \end{proof}

 %%%exgroupoid
\subsubsection{\bf Example}
\label{Examplegroupoid}
Let $\mathcal G$ be a locally compact groupoid, equipped with a left Haar system $(\lambda^u)_{u\in {\mathcal G}^o}$ and a quasi-invariant Radon measure $\nu$ on $\mathcal G^o$. Then, using the same notations as in \ref{defMQG}, we get that $(\bf C_0$$(\mathcal G^{(o)}), \bf C_0$$(\mathcal G), r_\mathcal G, s_\mathcal G, \nu, T_\mathcal G, T^{-1}_\mathcal G, \Gamma_\mathcal G)$ is a locally compact quantum groupoid, and the measured quantum groupoid constructed by \ref{thMQG} is then $(L^\infty(\mathcal G^{(o)}, \nu), L^\infty(\mathcal G, \mu_l), r_\mathcal G, s_\mathcal G, \Gamma_\mathcal G, T_\mathcal G, T^{-1}_\mathcal G, \nu)$ where $\mu_l=\int_{\mathcal G^{(o)}}\lambda^ud\nu$ and $\mu_r=\int_{\mathcal G^{(0)}}\lambda_ud\nu$. 

The canonical locally compact quantum groupoid constructed by \ref{ThLCQG}, is then :

$(\bf C^*(\nu), \bf C^*(\mu_l)$$, r_\mathcal G, s_\mathcal G, \nu, T_\mathcal G, T^{-1}_\mathcal G, \Gamma_\mathcal G)$, where 
$\bf C^*(\nu)$ (resp. $\bf C^*(\mu_l)$) is the norm closure of $L^1(\mathcal G^{(0)}, \nu)\cap L^\infty (\mathcal G^{(0)}, \nu)$ (resp. $L^1(\mathcal G, \mu_1)\cap L^1(\mathcal G, \mu_r)\cap L^\infty (\mathcal G, \mu_l)$) in $L^\infty (\mathcal G^{(0)}, \nu)$ (resp. in $L^\infty (\mathcal G, \mu_l)$), which is not, in general, the initial one.

Moreover, this example shows that there is no hope for a unicity theorem for a locally compact sub-quantum groupoid of a measured quantum groupoid. In particular, if $X$ is a locally compact space, and $\mu$ a Radon measure on $X$ with full support, then $X$ (with $X^{(0)}=X$ and then $X^{(2)}=X$) can be considered as a locally compact groupoid (called space groupoid), which has $id$ as left Haar system, and $\mu$ as quasi-invariant measure. Therefore, $(\bf C$$_0(X), \bf C$$_0(X), id, id, \mu, id, id, id)$ is a locally compact quantum groupoid. The measured quantum groupoid associated by \ref{thMQG} is then $(L^\infty (X, \mu), L^\infty (X, \mu), id, id, id, id, id, \mu)$. But then, any dense sub-$\bf C^*$-algebra $A$ of $L^\infty (X, \mu)$ such that $A\cap\L^1(, \mu)$ is dense in $A$ gives another locally compact quantum groupoid.

%%%%%%%Ramsay
\subsubsection{\bf{Ramsay's theorem}} [Ra]
\label{ramsay}
{\it Let $\mathcal G$ be a measured groupoid, with $\mathcal G^{(0)}$ as space
of units, and $r$ and $s$ the range and source functions from $\mathcal G$ to $\mathcal G^{(0)}$, with a Haar system $(\lambda^u)_{u\in \mathcal G^{(0)}}$ and a quasi-invariant measure $\nu$ on $\mathcal G^{(0)}$. Let us write $\mu=\int_{\mathcal G^{(0)}}\lambda^ud\nu$. Let $\Gamma_{\mathcal G}$, $r_{\mathcal G}$, $s_{\mathcal G}$ be the morphisms associated in \ref{thabelian}. Then, there exists a locally compact groupoid $\tilde{\mathcal G}$, with set of units $\tilde{\mathcal G}^{(0)}$, with a Haar system $(\tilde{\lambda}^u)_{u\in \tilde{\mathcal G}^{(0)}}$, and a quasi-invariant measure $\tilde{\nu}$ on $\tilde{\mathcal G}^{(0)}$, such that, if $\tilde{\mu}=\int_{\tilde{\mathcal G}^{(0)}}\tilde{\lambda}^ud\tilde{\nu}$, we get that the abelian measured quantum groupoids $\gG(\mathcal G)$ and $\gG(\tilde{\mathcal G})$ are isomorphic. }

%%%%%%transformation
\subsection{Locally compact transformation groupoid}
\label{transformation}
In \cite{ET} had been defined the "measured quantum transformation groupoids", which are quantum groupoids constructed on a crossed product of a von Neumann algebra $N$ by a specific action of a quantum group $\bf G$; more precisely, $N$ must be a braided commutative Yetter-Drinfel'd algebra. Frank Taipe (\cite{Ta}) had proved that, if the quantum group is constructed from an algebraic quantum group, as defined by Van Daele in  (\cite{VD2}), acting on a $*$-algebra, with appropriate axioms, then the crossed product of the action of the quantum group on the norm closure of the algebra is a locally compact quantum transformation groupoid. Here, mimicking this thesis, we obtain the same result for appropriate action of a locally compact quantum group on a $\bf C^*$-algebra. I must thank Frank Taipe who gave me new ideas on the subject. 

%%%YD
\subsubsection{{\bf Definition}} \cite{ET}(2.4)
\label{YD}
Let $\bf G$$=(M, \Gamma, \varphi, \varphi\circ R)$ be (the von Neumann version of) a locally compact quantum group, and $\widehat{\bf G}=(\widehat{M}, \widehat{\Gamma}, \widehat{\varphi}, \widehat{\varphi}\circ\widehat{R})$ its dual. A $\bf G$-Yetter-Drinfeld algebra is a von Neumann algebra $N$, with a left action $\alpha$ of $\bf G$ on $N$, and a left action $\widehat{\alpha}$ of $\widehat{\bf G}$ on $N$, such that, for all $x\in N$:
\[(id\otimes\alpha)\widehat{\alpha}(x)=Ad(\sigma W\otimes 1)(id\otimes\widehat{\alpha})\alpha(x)\]

%%%braided-commutative
\subsubsection{{\bf Definition}} \cite{ET}(2.5.3)
\label{braided-commutative}
For any action $\alpha$ of $\bf G$ on a von Neumann algebra $N$, we define the action $\alpha^c$ of $\bf G^c$ on $N^o$ and the action $\alpha^o$ of $\bf G^o$ on $N^o$, where $\alpha^c=(j\otimes .^o)\alpha$, and $\alpha^o=(R\otimes .^o)\alpha$, where, for any $x\in M$ $j(x)=Jx^*J$. 

If $(N, \alpha, \widehat{\alpha})$ is a $\bf G$-Yetter-Drinfel'd algebra, then, are equivalent :

(i) $\alpha^c(N^o)$ and $\widehat{\alpha}^c(N^o)$ commute;

(ii)$ \alpha^o(N^o)$ and $\widehat{\alpha}^o(N^o)$ commute.

In that case, we shall say that $(N, \alpha, \widehat{\alpha})$ is braided-commutative. 
Let $\nu$ be a normal faithful semi-finite weight on $N$, and let $U^\alpha_\nu=J_{\widetilde{\nu}}(\widehat{J}\otimes J_\nu)$ be the standard implememtation of $\alpha$ on $H\otimes H_\nu$ (\cite{ET}2.2). Let us define an injective anti-$*$-homomorphism $\beta$ by $\beta(x)=U^\alpha_\nu \widehat{\alpha}^o(x^o) (U^\alpha_\nu)^*$; then $(N, \alpha, \widehat{\alpha})$ is braided-commutative if and only if $\beta(N)$ is included in the crossed product $\bf G$$\ltimes_\ga N$ of $N$ by $\alpha$. 

%%%%transgroupoid
\subsubsection{{\bf Definition} }\cite{ET}
\label{transgroupoid}
If $(N, \alpha, \widehat{\alpha})$ is a braided commutative $\bf G$-Yetter-Drinfel'd algebra, we can construct on the crossed product $\bf G$$\ltimes_\ga N$ a structure of Hopf bimodule (\cite{ET}, 4) by the following way :

On the crossed-product $\bf G$$\ltimes_\ga N$, let $\tilde{\alpha}$ be the dual action of $\widehat{\bf G}^o$ given, for $X\in \bfG\ltimes_\ga N$, by 
\[\widehat{\alpha}(X)=(\widehat{W}^{o*}\otimes 1)(1\otimes X)(\widehat{W}^o\otimes 1)\]
For any $\eta\in H$ and $p\in\gN_\nu$, the vector $U_\alpha^\nu(\eta\otimes J_\nu\Lambda_\nu(p))$ belongs to $D_\alpha (H\otimes H_\nu)$ (\cite{ET}4.3(i))and we can define a unitary $V_1$ from $(H\otimes H_\nu)\underset{\nu}{_\beta\otimes_\alpha} (H\otimes H_\nu)$ onto $H\otimes H\otimes H_\nu$ by (for all $\Xi$ in $H\otimes H_\nu$) :
\[V_1(\Xi\underset{\nu}{_\beta\otimes_\alpha}U_\alpha^\nu(\eta\otimes J_\nu\Lambda_\nu(p)))=\eta\otimes\beta(p^*)\Xi\]
By (\cite{ET}, 2.2), we have that $U_\alpha^\nu=J_{\widetilde{\nu}}(\widehat{J}\otimes J_\nu)$ and, for $x\in\gN_\nu$ and $y\in\gN_{\widehat{\varphi}}$, $(y\otimes 1)\alpha(x)$ belongs to $\gN_{\widetilde{\nu}}$, and $J_{\widetilde{\nu}}\Lambda_{\widetilde{\nu}}((y\otimes 1)\alpha(x))=U^\alpha_\nu(\widehat{J}\Lambda_{\widehat{\varphi}}(y)\otimes J_\nu\Lambda_\nu(x))$.

For any $X\in \bf G$$\ltimes_\ga N$, let us define $\widetilde{\Gamma}(X)=V_1^*\widetilde{\alpha}(X)V_1$; then $(N, \bfG\ltimes_\ga N, \alpha, \beta, \widetilde{\Gamma})$ is a Hopf-bimodule (\cite{ET}, 4.4 ). This Hopf-bimodule can be equipped with a co-inverse $\widetilde{R}$, a left-invariant operator-valued weight $T_{\widetilde{\alpha}}$ and a right-invariant operator-valued weight $\widetilde{R}\circ T_{\widetilde{\alpha}}\circ\widetilde{R}$ (\cite{ET} 5.4) ; 
if the automorphism groups $\sigma_t^{\widetilde{\nu}}$ and $\sigma_s^{\widetilde{\nu}\circ\widetilde{R}}$ commute, 
for all $s$, $t$ in $R$, 
then $(N, \bf G$$\ltimes_\ga N, \alpha, \beta, \widetilde{\Gamma}, T_{\widetilde{\alpha}}, \widetilde{R}\circ T_{\widetilde{\alpha}}\circ\widetilde{R}, \nu)$ is a measured quantum groupoid (\cite{ET}, 5.9). 
If $z\in\widehat{M}$ and $x\in N$, we have 
\[\widetilde{\Gamma}[(z\otimes 1)\alpha(x)]=V_1^*(\widehat{\Gamma}^o(z)\otimes 1)V_1(\alpha(x)\underset{\nu}{_\beta\otimes_\alpha}1)\]

%%%C*algebra
\subsubsection{{\bf Definition}}
\label{C*algebra}
Let $\alpha$ an action of a (von Neumann version of a) locally compact quantum group $\bf G$ on a von Neumann algebra $N$; let $\bf G$$=(A, \Gamma_{|A}, \varphi_{|A}, \varphi\circ R_{|A})$ be the $\bf C^*$ version of $\bf G$, and $\widehat{\bf G}=(\widehat{A}, \widehat{\Gamma}_{|\widehat{A}}, \widehat{\varphi}|_{\widehat{A}}, \widehat{\varphi}\circ\widehat{R}|_{\widehat{A}})$ its dual; let $B$ a sub-$\bf C^*$-algebra of $N$, weakly dense in $N$; then $\alpha_{|N}$ is an action of $\bf G$ on $B$ if $\alpha(B)\subset M(A\otimes B)$. 

The crossed-product $\bf G\ltimes_\ga B$ is equal to norm closure of the linear space generated by elements of the form $(y\otimes 1)\alpha(b)$, for any $y\in \widehat{A}$, and $b\in B$ (\cite{BS}, 7.2).

Let us suppose that the weight $\nu_{|B}$ is semi-finite; then $T_{\widetilde{\alpha}|\bfG\ltimes_\ga B}$ is a semi-finite operator-valued weight from $\bf G$$\ltimes_\ga B$ to $\alpha(M(B))$.

%%%C*YD
\subsubsection{{\bf Definition}}
\label{C*YD}
Let $\bf G $$ =(M, \Gamma, \varphi, \varphi\circ R)$ be (the von Neumann version of) a locally compact quantum group, $(A, \Gamma_{|A}, \varphi_{|A}, \varphi_{|A}\circ R_{|A})$ its $\bf C^*$-version, and $\widehat{\bf G}=(\widehat{M}, \widehat{\Gamma}, \widehat{\varphi}, \widehat{\varphi}\circ\widehat{R})$ its dual (von Neumann version) or 
$(\widehat{A}, \widehat{\varphi}_{|\widehat{A}}, \widehat{\varphi}_{|\widehat{A}}\circ \widehat{R}_{|\widehat{A}}, \widehat{\Gamma}_{|\widehat{A}})$ ($\bf C^*$ version). Let $(N, \alpha, \widehat{\alpha})$ 
be a $\bf G$-Yetter-Drinfeld algebra as defined in \ref{YD}. 
Let $B$ be sub-$\bf C^*$-algebra of $N$, such that $\alpha(B)\subset M(A\otimes B)$, and $\widehat{\alpha}(B)\subset M(\widehat{A}\otimes B)$; then, we say that $(B, \alpha_{|B}, \widehat{\alpha}_{|B})$ is a $\bf C^*$-$\bf G$-Yetter- Drinfel'd algebra. 
If $(N, \alpha, \widehat{\alpha})$ is braided-commutative (as defined in \ref{braided-commutative}) , we shall say that $(B, \alpha_{|B}, \widehat{\alpha}_{|B})$ is braided-commutative if $\beta(B)\subset \bfG\ltimes_\alpha B$ (as defined in \ref{C*algebra}). Then $R(\bf G$$\ltimes_\alpha B)=\bf G$$\ltimes_\alpha B$. 

%%%%lemGamma
\subsubsection{{\bf Lemma}}
\label{lemGamma}
{\it Let $x\in A$, and $y\in A\cap\gM_\varphi$; then, there exist $x_n\in A\cap\gN_\varphi$ and $y_n\in A$ such that, for all $\xi\in H$, $\Gamma(x)(\xi\otimes \Lambda_\varphi(y))$ is equal to the norm limit of $\Sigma_n(y_n\xi\otimes\Lambda_\varphi(x_n))$. }

\begin{proof} 
Let us write $y=y_1^*y_2$, with $y_1$ and $y_2$ in $A\cap\gN_\varphi$. We have then 
\[\Gamma(x)(\xi\otimes\Lambda_\varphi(y))=[\Gamma(x)(1\otimes y_1^*)](\xi\otimes\Lambda_\varphi(y_2))\]
By  \cite{KV1}4.4, we know that $\Gamma(x)(1\otimes y_1^*)$ belongs to $A\otimes A$, and, therefore is a norm limit of a sum $\Sigma_n(y_n\otimes x^1_n)$, with $x^1_n\in A$ and $y_n\in A$; then, we get that $\Gamma(x)(\xi\otimes \Lambda_\varphi(y))$ is the norm limit of $\Sigma_n ( y_n\xi\otimes \Lambda_\varphi (x^1_ny_2))$; writing $x_n=x^1_ny_2$, we get the result. \end{proof}
%%%%propGamma
\subsubsection{{\bf Proposition}}
\label{propGamma}
{\it Let $X\in\bf G$$\ltimes_\alpha B$, $a$ in $\gN_\nu$, $y\in\gN_{\widehat{\varphi}}$. Then $\widetilde{\Gamma}(X)\rho^{\beta, \alpha}_{J_{\widetilde{\nu}}\Lambda_{\widetilde{\nu}}((y\otimes 1)\alpha(a))}$ is equal to the norm limit of elements of the form $\rho^{\beta, \alpha}_{J_{\widetilde{\nu}}\Lambda_{\widetilde{\nu}}((z_p\otimes 1)\alpha(x_p))}(y_p\otimes 1)\alpha (a_p)$ (with $y_p\in\widehat{A}$, $z_p\in\widehat{A}\cap\gN_{\widehat{\varphi}}$, $a_p\in B$, $x_p\in B\cap\gN_\nu$)}.

\begin{proof}
Using \ref{transgroupoid} , we get that $\widetilde{\Gamma}((z\otimes 1)\alpha(x))[\Xi\underset{\nu}{_\beta\otimes_\alpha}J_{\widetilde{\nu}}\Lambda_{\widetilde{\nu}}((y\otimes 1)\alpha(a))]$ is equal to 

\[V_1^*(\widehat{\Gamma}^o(z)\otimes 1)V_1[\alpha(x)\Xi\underset{\nu}{_\beta\otimes_\alpha}U^\alpha_\nu(\widehat{J}\Lambda_{\widehat{\varphi}}(y)\otimes J_\nu\Lambda_\nu(a))]\]
and to :
\[V_1^*(\widehat{\Gamma}^o(z)\otimes 1)(\widehat{J}\Lambda_{\widehat{\varphi}}(y)\otimes \beta(a^*)\alpha(x)\Xi)
=
V_1^*(\widehat{\Gamma}^o(z)\otimes 1)(\Lambda_{\widehat{\varphi}}(\sigma_{-i/2}^{\widehat{\varphi}}(y^*))\otimes\beta(a^*)\alpha(x)\Xi)\]

Using \ref{lemGamma} applied to the $\bf C^*$-version of $\widehat{\bf G}$ and \ref{transgroupoid},  there exist $z_p\in \widehat{A}\cap\gN_{\widehat{\varphi}}$ and $y_p\in\widehat{A}$, such that it is equal to the norm limit of :
\[\Sigma_p V_1^*\Lambda_{\widehat{\varphi}}(z_p)\otimes (y_p\otimes 1)\beta(a^*)\alpha(x)\Xi\]
Using \ref{C*algebra}, we get that $(y_p\otimes 1)\beta(a^*)=\widetilde{R}[\alpha(a^*)(\widehat{R}(y_p)\otimes 1)]$ is the norm limit of a sum $\Sigma_n\widetilde{R}[(y_{n,p}\otimes 1)\alpha(a_n)]=\Sigma_n\beta(a_n)(\widehat{R}(y_{n,p})\otimes 1)$, with $y_{n,p}\in \widehat{A}$ and $a_n\in A\cap\gN_\nu$ such that $a_n$ is analytic and $\sigma_{-i/2}^\nu(a_n^*)\in \gN_\nu$ and, therefore, we get that 
\[\widetilde{\Gamma}((z\otimes 1)\alpha(x))[\Xi\underset{\nu}{_\beta\otimes_\alpha}J_{\widetilde{\nu}}\Lambda_{\widetilde{\nu}}((y\otimes 1)\alpha(a))]\] is equal to the norm limit of 
\[\Sigma_{p,n}V_1^*\Lambda_{\widehat{\varphi}}(z_p)\otimes \beta(a_n)(\widehat{R}(y_{n,p})\otimes 1)\alpha(x)\Xi=
\Sigma_{p,n}(\widehat{R}(y_{n,p})\otimes 1)\alpha(x)\Xi\underset{\nu}{_\beta\otimes_\alpha}U^\alpha_\nu(\Lambda_{\widehat{\varphi}}(z_p)\otimes \Lambda_\nu(a_n))\]
Therefore, we get that $\widetilde{\Gamma}((z\otimes 1)\alpha(x))\rho^{\beta, \alpha}_{J_{\widetilde{\nu}}\Lambda_{\widetilde{\nu}}((y\otimes 1)\alpha(a))}$ is the norm limit of 
\[\Sigma_{p,n}\rho^{\beta, \alpha}_{J_{\widetilde{\nu}}\Lambda_{\widetilde{\nu}}((z_p\otimes 1)\alpha(\sigma_{-i/2}^\nu(a_n^*))}(\widehat{R}(y_{n,p})\otimes 1)\alpha(x)\] 
By continuity of $\widetilde{\Gamma}$, we obtain the result. 
\end{proof}

%%%%lemX
\subsubsection{{\bf Lemma}}
\label{lemX}
{\it Let $y_n\in\widehat{A}\cap\gM_{\widehat{\varphi}}$, such that $(\Lambda_{\widehat{\varphi}}(y_n))_n$ is an orthonormal basis of $H$. Let $X\in\bfG\ltimes_\alpha B\cap \gN_{T_{\widetilde{\alpha}}}\cap\gN_{\widetilde{\nu}}$ (resp. $X'\in \widetilde{A}\cap \gN_{T_{\widetilde{\alpha}}}\cap\gN_{\widetilde{\nu}}$) and $\epsilon>0$. 

(i) There exists $N$ in $\bf N$, and, for any $n\in \bf N$, such that $1\leq n\leq N$, there exists  $a_{n,N}\in B$ such that 
\[||X-\Sigma_{n=1}^N(y_n\otimes 1) \alpha(a_{n,N})||<\epsilon\]

(ii) There exists $N$ in $\bf N$, and, for any $n\in \bf N$, such that $1\leq n\leq N$, there exist $a_n\in B\cap\gN_\nu$ such that :
\[||\Lambda_{\widetilde{\nu}}(X)-\Sigma_{n=1}^N\Lambda_{\widehat{\varphi}}(y_n)\otimes\Lambda_\nu(a_n)||^2<\epsilon\]

(iii) $\Lambda_{T_{\widetilde{a}}}(X)$ is the weak limit of $\Sigma_n\Lambda_{\widehat{\varphi}}(y_n)\otimes a_n$ and $T_{\widetilde{a}}(X^*X)$ is the weak limit of $\alpha(\Sigma_n a_n^*a_n)$.

(iv) Let us suppose that $T_{\widetilde{a}}(X^*X)$ belongs to $\alpha(B)$; then $T_{\widetilde{a}}(X^*X)$ is the norm limit of $\alpha(a_n^*a_n)$ and $\Lambda_{T_{\widetilde{a}}}(X)$ is the norm limit of $\Sigma_n\Lambda_{\widehat{\varphi}}(y_n)\otimes a_n$. }

\begin{proof}
Using \ref{C*algebra}, we get that, for all $\epsilon > 0$, there exist $z_p\in \widehat{A}\cap\gN_{\widehat{\varphi}}$ and $b_p\in B\cap\gN_\nu$ such that 
\[||X-\Sigma_{p=1}^P(z_p\otimes 1)\alpha(b_p)||< \epsilon\]
We have $\Lambda_{\widetilde{\nu}}(\Sigma_{p=1}^P(z_p\otimes 1)\alpha(b_p))=\Sigma_{p=1}^P\Lambda_{\widetilde{\varphi}}(z_p)\otimes\Lambda_\nu(b_p)$.
Using now the basis $(\Lambda_{\widehat{\varphi}}(y_n))$ we get that there exists $N$ and $\lambda_{p,n}\in \bf C$ such that, for all $1\leq p\leq P$, we have $\Lambda_{\widehat{\varphi}}(z_p)=\Sigma_{n=1}^N\lambda_{p,n}\Lambda_{\widehat{\varphi}}(y_n)$, and, therefore, $\Lambda_{\widetilde{\nu}}(\Sigma_{p=1}^P(z_p\otimes 1)\alpha(b_p))=\Sigma_{n=1}^N\Lambda_{\widehat{\varphi}}(y_n)\otimes\Lambda_\nu(\Sigma_{p=1}^P\lambda_{p,n}b_p)$. From which we get that $\Sigma_{p=1}^P(z_p\otimes 1)\alpha(b_p)=\Sigma_{n=1}^N(y_n\otimes 1)\alpha(\Sigma_{p=1}^P\lambda_{p,n}b_p)$. 

Writing $a_{n,N}=\Sigma_{p=1}^P\lambda_{p,n}b_p$, we get (i). 

There exists a family $\xi_n$ in $H_\nu$ such that $\Lambda_{\widetilde{\nu}}(X)=\Sigma_n \Lambda_{\widehat{\varphi}}(y_n)\otimes\xi_n$. We have $\widetilde{\nu}(X^*X)=\Sigma_n||\xi_n||^2$. Therefore, there exists $N$ such that $\Sigma_{n>N}||\xi_n||^2<\epsilon/2$. And, for all $n\leq N$, there exist $a_n\in B\cap\gN_{\nu}$ such that $||\xi_n-\Lambda_\nu(a_n)||^2<\epsilon/2N$. From which we get that 
\[||\Lambda_{\widetilde{\nu}}(X)-\Sigma_{n=1}^N\Lambda_{\widehat{\varphi}}(y_n)\otimes\Lambda_\nu(a_n)||^2=
||\Sigma_{n=1}^N\Lambda_{\widehat{\varphi}}(y_n)\otimes (\xi_n-\Lambda_\nu(a_n))||^2+||\Sigma_{n>N}\Lambda_{\widehat{\varphi}}(y_n)\otimes \xi_n||^2<\epsilon\]
which is (ii). 

Let now $x\in\gN_\nu$; we have, using (ii)  :
\begin{multline*}
\Lambda_{T_{\widetilde{\alpha}}(X)}\Lambda_\nu(x)=\Lambda_{\widetilde{\nu}}(X\alpha(x))=\Sigma_n\Lambda_{\widetilde{\nu}}((y_n\otimes 1)\alpha(a_nx))=\\\Sigma_n\Lambda_{\widehat{\varphi}}(y_n)\otimes\Lambda_\nu(a_nx)=\Sigma_n\Lambda_{\widehat{\varphi}}(y_n)\otimes a_n\Lambda_\nu(x)
\end{multline*}
from which we get (iii). 
Let us suppose now that $T_{\widetilde{a}}(X^*X)$ belongs to $\alpha(B)$; there exists $b_n\in B\cap\gN_\nu$ such that $T_{\widetilde{a}}(X^*X)$ is the norm limit of $\Sigma_nb_n^*b_n$. Let us define $Y=\Sigma_n\Lambda_{\widehat{\varphi}}(y_n)\otimes b_n\in B(H_\nu, H_{\widetilde{\nu}})$. We get that, for all $\xi\in H_\nu$, we have $||\Lambda_{T_{\widetilde{\alpha}}}(X)\xi||=||Y\xi||$; therefore, there exists an isometry $U\in B(H_{\widetilde{\nu}})$ such that $UY=\Lambda_{T_{\widetilde{\alpha}}}(X)$; so, we get that $a_n=(\omega_{\Lambda_{\widetilde{\varphi}}(y_n)}\otimes id)(U)b_n$; then, we have $\Sigma_{n=N}^\infty a_n^*a_n\leq \Sigma_{n=N}^\infty b_n^*b_n$, and, therefore $\Sigma_n a_n^*a_n$ is norm converging to $T_{\widetilde{\alpha}}(X^*X)$, and, then, we deduce that $\Lambda_{T_{\widetilde{\alpha}}}(X)$ is the norm limit of $\Sigma_n\Lambda_{\widehat{\varphi}}(y_n)\otimes a_n$, which is (iv). \end{proof} 

%%%propGamma2
\subsubsection{{\bf Proposition}}
\label{propGamma2}
{\it Let $X\in\bf G$$\ltimes_\alpha B$, and $Y\in\bf G$$\ltimes_\alpha B\cap\gN_{T_{\widetilde{\alpha}}}\cap\gN_{\widetilde{\nu}}$; lhen :

(i) Let us suppose that $T_{\widetilde{\alpha}}(Y^*Y)\in \alpha(A)$; then there exist $y_p\in\widehat{A}$, $z_p\in\widehat{A}\cap\gN_{\widehat{\varphi}}$, $a_p\in B$, $x_p\in B\cap\gN_\nu$ such that $\widetilde{\Gamma}(X)\rho^{\beta, \alpha}_{J_{\widetilde{\nu}}\Lambda_{\widetilde{\nu}}(Y)}$ is the norm limit of 
$\rho^{\beta, \alpha}_{J_{\widetilde{\nu}}\Lambda_{\widetilde{\nu}}((z_p\otimes 1)\alpha(x_p))}(y_p\otimes 1)\alpha (a_p)$.

(ii) Let $a$ in $\gN_\nu$. Then there exist $y_p\in\widehat{A}$, $z_p\in\widehat{A}\cap\gN_{\widehat{\varphi}}$, $a_p\in B$, $x_p\in B\cap\gN_\nu$ such that $\widetilde{\Gamma}(X)\rho^{\beta, \alpha}_{J_{\widetilde{\nu}}\Lambda_{\widetilde{\nu}}(Y\alpha(a))}$ is the norm limit of $\rho^{\beta, \alpha}_{J_{\widetilde{\nu}}\Lambda_{\widetilde{\nu}}((z_p\otimes 1)\alpha(x_p))}(y_p\otimes 1)\alpha (a_p)$.

(iii) Let $a$ in $\gN_\nu$, analytic with respect to $\sigma_\nu$. There exist $y_p\in\widehat{A}$, $z_p\in\widehat{A}\cap\gN_{\widehat{\varphi}}$, $a_p\in B$, $x_p\in B\cap\gN_\nu$ such that $\widetilde{\Gamma}(X)(\beta(a^*)\underset{\nu}{_\beta\otimes_\alpha} 1)\rho^{\beta, \alpha}_{J_{\widetilde{\nu}}\Lambda_{\widetilde{\nu}}(Y)}$ is the norm limit of 
\[\rho^{\beta, \alpha}_{J_{\widetilde{\nu}}\Lambda_{\widetilde{\nu}}((z_p\otimes 1)\alpha(x_p))}(y_p\otimes 1)\alpha (a_p)\]

(iv) There exist $y^n_p\in\widehat{A}$, $z^n_p\in\widehat{A}\cap\gN_{\widehat{\varphi}}$, $a^n_p\in B$ and $x^n_p\in B\cap\gN_\nu$ such that $\widetilde{\Gamma}(X)\rho^{\beta, \alpha}_{J_{\widetilde{\nu}}\Lambda_{\widetilde{\nu}}(Y)}$ is a norm limit of $\rho^{\beta, \alpha}_{J_{\widetilde{\nu}}\Lambda_{\widetilde{\nu}}((z_p^n\otimes 1)\alpha(x_p^n))}(y_p^n\otimes 1)\alpha (a_p^n)$}

\begin{proof}
Using \ref{lemX}(iv) applied to $Y$ and \ref{propGamma}, we get (i). 
Then, as 
\[T((Y\alpha(a))^*(Y\alpha(a))=\alpha(a^*)T(Y^*Y)\alpha(a)\]
 belongs to $\alpha(A)$, 
using (i) applied to $Y\alpha(a)$, we get (ii). 
As we have :
\[\widetilde{\Gamma}(X)(\beta(a^*)\underset{\nu}{_\beta\otimes_\alpha} 1)\rho^{\beta, \alpha}_{J_{\widetilde{\nu}}\Lambda_{\widetilde{\nu}}(Y)}=
\widetilde{\Gamma}(X)(1\underset{\nu}{_\beta\otimes_\alpha}\alpha(\sigma_{-i/2}(a^*))\rho^{\beta, \alpha}_{J_{\widetilde{\nu}}\Lambda_{\widetilde{\nu}}(Y)}
=\widetilde{\Gamma}(X)\rho^{\beta, \alpha}_{J_{\widetilde{\nu}}\Lambda_{\widetilde{\nu}}(Y\alpha(a))}\]
we get that (ii) implies (iii). 
Let $\Xi\in H_{\widetilde{\nu}}$; then there exist $a_n$ in $B$, analytic with respect to $\sigma_\nu$, such that $\Xi$ is the norm limt of 
$\beta(a_n^*)\Xi$; then $\widetilde{\Gamma}(X)\Xi\underset{\nu}{_\beta\otimes_\alpha}J_{\widetilde{\nu}}\Lambda_{\widetilde{\nu}}(Y)$ is the norm limit of 
\[\widetilde{\Gamma}(X)\beta(a_n^*)\Xi\underset{\nu}{_\beta\otimes_\alpha}J_{\widetilde{\nu}}\Lambda_{\widetilde{\nu}}(Y)\]
and, using (iii), there exist $y_p^n$, $z^n_p$, $a^n_p$ and $x^n_p$ such that it is the norm limit of elements of the form $(y_p^n\otimes 1) \alpha(a_p^n)\Xi\underset{\nu}{_\beta\otimes_\alpha}J_{\widetilde{\nu}}\Lambda_{\widetilde{\nu}}((z_p^n\otimes 1)\alpha(x^n_p))$. So, we get (iv). \end{proof}

%%%th
\subsubsection{{\bf Theorem}}
\label{th}
{\it  Let's use the notations of \ref{C*YD}; then :

$(B, \bf G$$\ltimes_\alpha B, \alpha_{|B}, \beta_{|B}, T_{\widetilde{\alpha}|_{\bf G\ltimes_\alpha B}}, (\widetilde{R}T\widetilde{R})_{|_{\bf G\ltimes_\alpha B}}	, \Gamma_{|_{\bfG\ltimes_\alpha B}})$
 is a locally compact quantum groupoid.}										
\begin{proof}

Using \ref{propGamma2}(iv), we get that $\widetilde{\Gamma}(\bf G$$\ltimes_\alpha B)\subset \bf G$$\ltimes_\alpha B\underset{\nu}{_\beta *_\alpha} B(H_{\widetilde{\nu}})$; as $\widetilde{R}(\bf G$$\ltimes_\alpha B)=\bf G$$\ltimes_\alpha B$ (\ref{C*YD}), we get that $\widetilde{\Gamma}(\bf G$$\ltimes_\alpha B)\subset \bf G$$\ltimes_\alpha B\underset{\nu}{_\beta *_\alpha}\bf G$$\ltimes_\alpha B$. \end{proof}

%%%%%%bibli

\end{document}